\documentclass[11pt]{article}
\usepackage{amsfonts}
\usepackage{amsgen,amsmath,amstext,amsbsy,amsopn,amssymb,subfigure,amsthm}
\usepackage[dvips]{graphicx}
\usepackage{comment}
\usepackage{enumitem}
\usepackage[numbers]{natbib}
\usepackage{booktabs}
\usepackage{blkarray}
\usepackage{algorithm}
\usepackage{algpseudocode}
\usepackage{url}
\usepackage{longtable}
\usepackage{mathtools}

\DeclarePairedDelimiter\floor{\lfloor}{\rfloor}

\usepackage[usenames]{color}

 \allowdisplaybreaks

\newtheorem{Definition}{Definition}
\newtheorem{Theorem}{Theorem}

\newtheorem{Lemma}{Lemma}
\newtheorem{Remark}{Remark}

\newtheorem{Hypothesis}{Hypothesis}
\newtheorem{Proposition}{Proposition}

\newcommand{\be}{\begin{equation}}
\newcommand{\ee}{\end{equation}}
\newcommand{\bea}{\begin{eqnarray}}
\newcommand{\eea}{\end{eqnarray}}
\newcommand{\beas}{\begin{eqnarray*}}
\newcommand{\eeas}{\end{eqnarray*}}

\newcommand{\br}{{\boldsymbol{r}}}
\newcommand{\bp}{{\boldsymbol{p}}}

\newcommand{\X}{{\mathbf{X}}}

\newcommand{\E}{{\mathbf{E}}}

\renewcommand{\S}{{\mathbf{S}}}

\newcommand{\V}{{\rm V}}

\newcommand{\A}{{\mathbf{A}}}
\newcommand{\Y}{{\mathbf{Y}}}
\newcommand{\Z}{{\mathbf{Z}}}

\newcommand{\Proj}{P}

\newcommand{\rank}{{\rm rank}}

\newcommand{\diag}{{\rm diag}}

\newcommand{\SVD}{{\rm SVD}}

\newcommand{\argmin}{\mathop{\rm arg\min}}
\newcommand{\argmax}{\mathop{\rm arg\max}}
\newcommand{\RR}{\mathbb{R}}

\makeatletter
\newcommand*{\rom}[1]{\expandafter\@slowromancap\romannumeral #1@}
\makeatother

\allowdisplaybreaks

\begin{document}

\title{
	Tensor SVD: Statistical and Computational Limits
	\footnote{Anru Zhang is Assistant Professor, Department of Statistics, University of Wisconsin-Madison, Madison, WI 53706, E-mail: anruzhang@stat.wisc.edu; Dong Xia is Visiting Assistant Professor, Department of Statistics, University of Wisconsin-Madison, Madison, WI 53706, E-mail: dongxia@stat.wisc.edu.}}
\author{Anru Zhang ~ and ~ Dong Xia\\
University of Wisconsin-Madison}
\date{}
\maketitle

\bigskip

\begin{abstract}
In this paper, we propose a general framework for tensor singular value decomposition (tensor SVD), which focuses on the methodology and theory for extracting the hidden low-rank structure from high-dimensional tensor data. Comprehensive results are developed on both the statistical and computational limits for tensor SVD. 
This problem exhibits three different phases according to the signal-to-noise ratio (SNR). In particular, with strong SNR, we show that the classical higher-order orthogonal iteration achieves the minimax optimal rate of convergence in estimation; with weak SNR, the information-theoretical lower bound implies that it is impossible to have consistent estimation in general; with moderate SNR, we show that the non-convex maximum likelihood estimation provides optimal solution, but with NP-hard computational cost; moreover, under the hardness hypothesis of hypergraphic planted clique detection, there are no polynomial-time algorithms performing consistently in general. 
\end{abstract}

\section{Introduction}
There is no need to argue the importance of singular value decomposition (SVD) in data analysis. As one of the most important tools in multivariate analysis, SVD along with the closely related formulation, i.e., principal component analysis (PCA), has been a mainstay of data analysis since more than a century ago, and widely used in various subjects. Attributed to the modern high-dimensional data, the popularity of SVD and PCA continues to surge in the recent decades, and many important variations, such as sparse SVD \citep{shen2008sparse,lee2010biclustering,yang2014sparse,yang2016rate}, matrix denoising \cite{candes2013unbiased,shabalin2013reconstruction,donoho2014minimax,gavish2014optimal}, sparse PCA \citep{zou2006sparse,cai2013sparse,birnbaum2013minimax}, robust PCA \citep{candes2011robust}, have been proposed and developed recently. Traditionally, most of the SVD and PCA results focused on exploiting low-rank structures from datasets in the form of matrices.

Motivated by modern scientific research, tensors, or high-order arrays, have been actively studied in machine learning, electrical engineering, and statistics. Some specific scientific applications  involving tensor data include neuroimaging analysis \citep{zhou2013tensor,zhang2016cross}, recommender systems \citep{karatzoglou2010multiverse,rendle2010pairwise}, computer vision \citep{lu2008mpca,liu2013tensor}, topic modeling \citep{anandkumar2014tensor-latent}, community detection \citep{anandkumar2014tensor-membership}, hyperspectral image compression \citep{li2010tensor}, spatiotemporal gene expression \citep{liu2017characterizing}, etc. A common objective in these problems is to dig out the underlying high-order low-rank structure, such as the singular subspaces and the whole low-rank tensors, buried in the noisy observations. To achieve this goal, we are in strong need of a statistical tool for tensor data that is the counterpart of regular singular value decomposition for traditional order-2 datasets. Richard and Montanari \cite{richard2014statistical}, Hopkins et al \cite{hopkins2015tensor}, Perry et al \cite{perry2016statistical} considered a rank-1 spiked tensor SVD statistical model and proposed various methods, including tensor unfolding and sum of square optimization (SOS). However, as far as we know, the statistical framework for general rank-$r$ high-order tensor SVD or PCA was not well established or studied in the literature.

In this paper, we propose a general framework of tensor singular value decomposition (tensor SVD). To be specific, suppose we are interested in a low-rank tensor $\X \in \mathbb{R}^{p_1\times p_2 \times p_3}$, which is observed with entry-wise corruptions as follows,
\begin{equation}\label{eq:tensor_pca_model}
\Y = \X + \Z.
\end{equation}
Here $\Z$ is the $p_1$-by-$p_2$-by-$p_3$ noisy tensor with $\{Z_{ijk}\}_{i,j,k=1}^{p_1, p_2, p_3}\overset{iid}{\sim} N(0, \sigma^2)$; $\X$ is a fixed tensor with low Tucker ranks in the sense that all fibers of $\X$ along three directions (i.e., counterpart of matrix columns and rows for tensors, also see Section \ref{sec:notation_preliminary} for formal definitions) lie in low-dimensional subspaces, say $U_1$, $U_2$, and $U_3$, respectively. Our goal is to estimate $U_1, U_2, U_3$, and $\X$ from the noisy observation $\Y$. 

It is worth mentioning that the analog of this problem when $\X$ is an order-2 tensor, i.e., a matrix, has been previously studied in the context of matrix denoising in \cite{candes2013unbiased,donoho2014minimax,gavish2014optimal,cai2016rate}. For the matrix denoising problem, the best low-rank matrix approximation provides the optimal results, which can be calculated efficiently via singular value decomposition, as guaranteed by the well-regarded Eckart-Young-Mirsky Theorem.

Although there have been significant efforts in developing methodologies and theories for matrix SVD or matrix denoising, there is a paucity of literature on the analogous question for tensors of order three or higher. In fact, SVD for high-order tensors is much more difficult than its counterpart for matrices in various aspects. First, tensors have more involved structures along three or more ways, while the traditional tools for matrices could typically incorporate two ways. As we will see later, one may achieve a sub-optimal result by simply ignoring the structure beyond two ways. Second, many operations for matrices, such as operator norm, singular value decomposition, are either not well defined or computational NP-hard for high-order tensors \citep{hillar2013most}. Third, high-order tensors often bring about high-dimensionality and impose significant computational challenges. For example, a 500-by-500-by-500 tensor contains 12,500,000 entries. All these characteristics make the tensor SVD distinct from the classical matrix setting.
	
The best low-rank tensor approximation, or equivalently the maximum likelihood estimation (MLE), is a straightforward solution for tensor SVD. However, MLE is non-convex and computationally NP-hard in general (see, e.g. Hillar and Lim \cite{hillar2013most}). De Lathauwer, De Moor, and Vandewalle instead introduced the higher order SVD (HOSVD) \cite{de2000multilinear} and higher order orthogonal iteration (HOOI) \cite{de2000best}, which aims to approximate the best low-rank approximation via efficient spectral and power iteration method. Since then, HOSVD and HOOI have been widely studied in the literature (see, e.g. \cite{vasilescu2003multilinear,sheehan2007higher,costantini2008higher,haardt2008higher,liu2014generalized}). However as far as we know, many theoretical properties of these procedures, such as the error bound and the necessary iteration times, still remain unclear.

In this paper, we develop comprehensive results on both the statistical and computational limits for tensor SVD. To be specific, we establish upper bounds on estimation errors for both higher order orthogonal iteration (HOOI) and maximum likelihood estimator (MLE). It is also shown that HOOI converges within a logarithm factor of iterations. Then the matching information-theoretical lower bounds over a large class of low-rank tensors are correspondingly introduced. To the best of our knowledge, we are among the first to develop the statistical guarantees for both HOOI and MLE. Let the Tucker rank of $\X$ be $(r_1,r_2,r_3)$ (see formal definition in Section \ref{sec:theory}). The statistical and computational barriers of tensor SVD problem rely on a key factor $\lambda$, i.e., the smallest non-zero singular values of matricizations of $\X$ (also see formal definition in Section \ref{sec:theory}), which essentially measures the signal strength of the problem. When $p=\min\{p_1, p_2, p_3\}$, $p_k \leq Cp$, $r_k \leq Cp^{1/2}$ for $k=1,2,3$ and a constant $C>0$, our main results can be summarized into the following three phases according to signal-to-noise ratio (SNR): $\lambda/\sigma$.
\begin{enumerate}
	\item When $\lambda/\sigma = p^{\alpha}$ for $\alpha \geq 3/4$, the scenario is referred to as the {\bf strong SNR case}. The fast higher order orthogonal iteration (HOOI) recovers $U_1, U_2, U_3$, and $\X$ with the minimax optimal rate of convergence over a general class of cases. 

	\item When $\lambda/\sigma = p^{\alpha}$ for $\alpha < 1/2$, we refer to this case as the {\bf weak SNR case} and propose the minimax lower bound to show that there are no consistent estimators of $U_1, U_2, U_3$, or $\X$;
	\item When $\lambda/\sigma = p^{\alpha}$ for $1/2 \leq \alpha <3/4$, the scenario is referred to as the {\bf moderate SNR case}. We provide a computational lower bound to show that no polynomial-time algorithm can recover $U_1, U_2, U_3$ consistently based on an assumption of hypergraphic planted clique detection. Meanwhile, the maximum likelihood estimator, although being computationally intractable, achieves optimal rates of convergence over a general class of cases.
\end{enumerate}
Especially when the tensor is rank-1, our results in the strong SNR case confirm a heuristic conjecture recently raised by Richard and Montanari \cite{richard2014statistical} that the tensor unfolding method yields reliable estimates for order-3 spiked tensors provided that $\lambda/\sigma > \Omega(p^{3/4})$. It is also noteworthy that our results can be further generalized to fourth or higher order tensors, or when the noise $\Z$ is i.i.d. sub-Gaussian distributed. 

Our work is also related to several recent results in literature. For example,  \cite{richard2014statistical,hopkins2015tensor,zheng2015interpolating,anandkumar2016homotopy} considered the extraction of rank-1 symmetric tensors from i.i.d. (symmetric) Gaussian noise, which is a rank-1 special case of our tensor SVD model; \cite{anandkumar2014guaranteed,sun2015provable} considered the CP low-rank tensor decomposition based on noisy observations; \cite{perry2016statistical} considered the statistical limit of detecting and estimating a randomly sampled rank-one structure from a symmetric random Gaussian tensor; \cite{allen2012regularized,allen2012sparse} considered the regularized tensor factorizations with/without sparsity; \cite{qi2016uniqueness} and \cite{lu2016tensor} further considered non-negative tensor decomposition and robust tensor principal component analysis; \cite{liu2017characterizing} focused on orthogonal decomposable tensor SVD problem; Lesieur et al \cite{lesieur2017statistical} considered a Bayesian symmetric spiked tensor estimation model -- an approximate message passing algorithm (AMP) was particularly introduced and the rigorous asymptotic analysis for statistical and computational phase transitions were performed on high-order, symmetric, and rank-1 tensor estimation. It should be noted that different from previous works, we perform non-asymptotic analysis for tensor SVD, where the signal tensor $\X$ can be generally Tucker-rank-$r$, non-random, and asymmetric. Also, to the best of our knowledge, we are among the first to provide a comprehensive analysis of both statistical and computational optimality of tensor SVD.

The rest of the article is organized as follows. After a brief explanation of basic notation and tensor algebra in Section \ref{sec:notation_preliminary}, we state the fast higher order orthogonal iteration and the non-convex maximum likelihood estimation for tensor SVD in Section \ref{sec:mle-and-spectral}. The statistical limits in the context of minimax optimality are provided for strong, weak, and moderate SNR cases respectively in Section \ref{sec:theory}.  Then we further discuss the computational barriers in the moderate SNR case in Section \ref{sec:moderate}. Simulation studies are provided in Section \ref{sec:simulation} to justify the theoretical results of this paper. We briefly discuss the extension of the results to fourth or higher order tensors and i.i.d. sub-Gaussian noise cases in Section \ref{sec:discuss}. The proofs of all technical results are given in Section \ref{sec:proof} and the supplementary materials.

\section{Tensor SVD: Methodology}\label{sec:method}

\subsection{Notation, Preliminaries, and Tensor Algebra}\label{sec:notation_preliminary}

In this section, we start with basic notation, preliminaries, and tensor algebra to be used throughout the paper. For $a, b\in \mathbb{R}$, let $a\wedge b= \min\{a, b\}$, $a\vee b = \max\{a, b\}$. For two sequences $\{a_i\}, \{b_i\}$, if there are two constants $C, c>0$ such that $ca_i \leq b_i \leq Ca_i$ for all $i\geq 1$, we denote $a\asymp b$. We use $C, c, C_0, c_0, \ldots$ to denote generic constants, whose actual values may vary from time to time. Particularly, the uppercase and lowercase letters represent large and small constants, respectively. The matrices are denoted as capital letters, $U_1, V_1, A$, etc. Especially, $\mathbb{O}_{p,r}:=\{U\in\mathbb{R}^{p\times r}: U^\top U=I_r\}$ is the set of all $p$-by-$r$ matrices with orthonormal columns. For any matrix $A\in\mathbb{R}^{p_1\times p_2}$, let $\sigma_1(A)\geq \cdots\geq \sigma_{p_1\wedge p_2}(A) \geq 0$ be the singular values in non-increasing order. We are particularly interested in the smallest singular value of $A$: $\sigma_{\min}(A)=\sigma_{p_1\wedge p_2}(A)$. In addition, the class of matrix Schatten $q$-norms will be used: $\|A\|_q=\Big(\sum_{j=1}^{p_1\wedge p_2}\sigma_j^q(A)\Big)^{1/q}$. Specific instances of Schatten $q$-norms include the Frobenius norm (i.e., Schatten 2-norm), $\|A\|_{\rm F} = \sqrt{\sum_{i=1}^{p_1} \sum_{j=1}^{p_2} A_{ij}^2} = \sqrt{\sum_{j=1}^{p_1\wedge p_2}\sigma_j^2(A)}$, and spectral norm (i.e., Schatten $\infty$-norm), $\|A\| = \sigma_1(A) = \max_{v\in \mathbb{R}^{p_2}} \frac{\|A v\|_2}{\|v\|_2}$. We also use $\SVD_r(A)$ to denote the leading $r$ left singular vectors of $A$, so that $\SVD_r(A)\in \mathbb{O}_{p_1, r}$. Define the projection operator $P_A = A(A^\top A)^\dagger A^\top$. Here $(\cdot)^\dagger$ represents the psudo-inverse. If $A = U\Sigma V^\top$ is the SVD, $P_A$ can be equivalently written as $P_A = UU^\top$. For any two matrices, say $U\in \mathbb{R}^{p_1\times r_1}, V\in \mathbb{R}^{p_2\times r_2}$, we also let $U\otimes V\in \mathbb{R}^{(p_1p_2)\times (r_1r_2)}$ be their outer product matrix, such that $\left(U \otimes V\right)_{[(i-1)p_2+j, (k-1)r_2+l]} = U_{ik}\cdot V_{jl}$, for $i=1,\ldots, p_1, j=1,\ldots, p_2, k=1,\ldots, r_1$, and $l=1,\ldots, r_2$. We adopt the R convention to denote submatrices: $A_{[a:b, c:d]}$ represents the submatrix formed by $a$-to-$b$-th rows and $c$-to-$d$-th columns of the matrix $A$; we also use $A_{[a:b, :]}$ and $A_{[:, c:d]}$ to represent $a$-to-$b$-th full rows of $A$ and $c$-to-$d$-th full columns of $A$, respectively. 

We use $\sin\Theta$ distances to measure the difference between singular subspaces. To be specific, for any two $p\times r$ matrices with orthonormal columns, say $U$ and $\hat U$, we define the principal angles between $U$ and $\hat{U}$ as $\Theta(U,\hat U)={\rm diag}\big(\arccos(\sigma_1),\ldots,\arccos(\sigma_r)\big) \in \mathbb{R}^{r\times r}$, where $\sigma_1\geq \ldots\geq \sigma_r\geq 0$ are the singular values of $U^\top \hat U$. 
The Schatten $q$-$\sin\Theta$-norm is then defined as 
$$\|\sin\Theta(U,\hat U)\|_q = \left(\sum_{i=1}^r \sin^q\left(\arccos(\sigma_i)\right)\right)^{1/q} = \left(\sum_{i=1}^r \left(1-\sigma_i^2\right)^{q/2}\right)^{1/q},\quad 1\leq q \leq +\infty.$$ 
The readers are referred to Lemma \ref{lemma:norms} in the supplementary materials and Lemma 1 in \cite{cai2016rate} for more discussions on basic properties of $\sin\Theta$ distances.

Throughout this paper, we use the boldface capital letters, e.g. $\X, \Y, \Z$, to note tensors. To simplify the presentation, the main context of this paper is focused on third order tensor. The extension to 4-th or higher tensors is briefly discussed in Section \ref{sec:discuss}. The readers are also referred to \cite{kolda2009tensor} for a more detailed tutorial of tensor algebra. For any tensor $\X\in \RR^{p_1\times p_2\times p_3}$, define its mode-1 matricization as a $p_1$-by-$(p_2p_3)$ matrix $\mathcal{M}_1(\X)$ such that
$$
\left[\mathcal{M}_1(\X)\right]_{i, (j-1)p_3 + k}=X_{ijk},\qquad \forall 1\leq i \leq p_1, 1\leq j \leq p_2, 1\leq k \leq p_3.
$$
In other words, $\mathcal{M}_1(\X)$ is composed of all mode-1 fibers, $\{(\X_{[:,i_2,i_3]})\in\mathbb{R}^{p_1}: 1\leq i_2\leq p_2, 1\leq i_3 \leq p_3\}$, of $\X$. The mode-2 and mode-3 matricizations, i.e., $\mathcal{M}_2(\X)\in \mathbb{R}^{p_2\times(p_3p_1)}$ and $\mathcal{M}_3\in \mathbb{R}^{p_3\times(p_1p_2)}$, are defined in the same fashion. We also define the marginal multiplication $\times_1:\RR^{p_1\times p_2\times p_3}\times \RR^{r_1\times p_1}\to \RR^{r_1\times p_2\times p_3}$ as
$$
\X\times_1 Y=\left(\sum_{i'=1}^{p_1} X_{i'jk}Y_{i, i'}\right)_{1\leq i \leq r_1, 1\leq j \leq p_2, 1\leq k\leq p_3}.
$$
Marginal multiplications $\times_2$ and $\times_3$ can be defined similarly. 

Different from matrices, there is no universal definition for tensor ranks. We particularly introduce the following Tucker ranks (also called multilinear ranks) of $\X$ as
\begin{equation*}
	\begin{split}
	r_1&= \rank_1(\X) = {\rm rank}(\mathcal{M}_1(\X))\\
	&={\rm dim}({\rm span}\{\X_{[:,i_2,i_3]}\in\mathbb{R}^{p_1}: 1\leq i_2\leq p_2, 1\leq i_3\leq p_3\}).\\
	\end{split}
\end{equation*}
$r_2 = \rank_2(\X)$ and $r_3 = \rank_3(\X)$ can be similarly defined. Note that, in general, $r_1, r_2, r_3$ satisfy $r_1 \leq r_2r_3, r_2 \leq r_3 r_1, r_3 \leq r_1 r_2$, but are not necessarily equal. We further denote $\rank(\X)$ as the triplet: $(r_1, r_2, r_3)$. The Tucker rank $(r_1, r_2, r_3)$ is also closely associated with the following Tucker decomposition. Let $U_1\in \mathbb{O}_{p_1, r_1}, U_2\in \mathbb{O}_{p_2, r_2}, U_3\in \mathbb{O}_{p_3, r_3}$ be the left singular vectors of $\mathcal{M}_1(\X)$, $\mathcal{M}_2(\X)$ and $\mathcal{M}_3(\X)$ respectively, then there exists a core tensor $\S\in \mathbb{R}^{r_1\times r_2\times r_3}$ such that
\begin{equation}
\label{eq:tucker}
\X = \S \times_1 U_1 \times_2 U_2 \times_3 U_3, \quad \text{or}\quad X_{ijk} = \sum_{i'=1}^{r_1} \sum_{j'=1}^{r_2} \sum_{k'=1}^{r_3} S_{i'j'k'}(U_{1})_{i, i'}(U_{2})_{j, j'}(U_{3})_{k, k'}.
\end{equation}
Expression \eqref{eq:tucker} is widely referred to as the Tucker decomposition of $\X$. 
Finally, to measure the tensor estimation error, we introduce the following tensor Frobenius norm,
$$
\|\X\|_{\rm F}=\Big(\sum_{i=1}^{p_1}\sum_{j=1}^{p_2} \sum_{k=1}^{p_3} X^2_{ijk}\Big)^{1/2}.$$

\subsection{Maximum Likelihood Estimator and Higher Order Orthogonal Iteration}\label{sec:mle-and-spectral}

In this section, we discuss the methodology for tensor SVD. Given the knowledge of Tucker decomposition, the original tensor SVD model \eqref{eq:tensor_pca_model} can be cast as follows,
\begin{equation}\label{eq:tensor_pca_mdel_2}
\Y = \X + \Z = \S \times_1 U_1\times_2 U_2\times_3 U_3 + \Z, \quad \Z \overset{i.i.d.}{\sim}N(0, \sigma^2),
\end{equation}
where $U_1 \in \mathbb{O}_{p_1, r_1}$, $U_2 \in \mathbb{O}_{p_2, r_2}$, $U_3 \in \mathbb{O}_{p_3, r_3}$, and $\S\in \mathbb{R}^{r_1\times r_2\times r_3}$. Our goal is to estimate $U_1, U_2, U_3$, and $\X$ from $\Y$. Clearly, the log-likelihood of Model \eqref{eq:tensor_pca_model} can be written (ignoring the constants) as $\mathcal{L}\left(\Y | \X\right) = -\frac{1}{\sigma^2} \left\|\Y - \X \right\|_{\rm F}^2,$ then it is straightforward to apply the maximum likelihood estimator for estimation, 
$$\hat{\X}^{\rm mle} = \argmin_{\rank(\X) \leq r_1, r_2, r_3}\left\|\Y - \X\right\|_{\rm F},$$ 
$$\hat{U}_k^{\rm mle} = \SVD_{r_k}(\hat{\X}^{\rm mle}), \quad k=1,2,3.$$
Intuitively speaking, the MLE seeks the best rank-($r_1, r_2, r_3$) approximation for $\Y$ in Frobenius norm. By Theorems 4.1 and 4.2 in \cite{de2000best}, MLE can be equivalently written as 
\begin{equation}\label{eq:mle}
\begin{split}
& \hat{U}_1^{\rm mle}, \hat{U}_2^{\rm mle}, \hat{U}_3^{\rm mle} = \argmax_{V_k\in \mathbb{O}_{p_k ,r_k}} \left\|\Y\times_1 V_1^\top \times_2 V_2^\top \times_3 V_3^\top\right\|_{\rm F}^2,\\
& \hat\X^{\rm mle}=\Y\times_1 \Proj_{\hat U_1^{\rm mle}}\times_2 \Proj_{\hat{U}_2^{\rm mle}}\times_3 \Proj_{\hat{U}_3^{\rm mle}}.
\end{split}
\end{equation}
As we will illustrate later in Section \ref{sec:theory}, such estimators achieve optimal rate of convergence in estimation errors. On the other hand, \eqref{eq:mle} is non-convex and computationally NP-hard even when $r=1$ (see, e.g., \citep{hillar2013most}). Then MLE may not be applicable in practice.

To overcome the computational difficulties of MLE, we consider a version of higher order orthogonal iteration (HOOI) \cite{de2000best}. The procedure includes three steps: spectral initialization, power iteration, and tensor projection. The first two steps produce optimal estimations of loadings $U_1, U_2, U_3$. The final step yields an optimal estimator of the underlying low-rank tensor $\X$. It is helpful to present the procedure of HOOI in detail here.
\begin{enumerate}
	\item[Step 1] (Spectral initialization) Since $U_1, U_2$, and $U_3$ respectively represent the singular subspaces of $\mathcal{M}_1(\X)$, $\mathcal{M}_2(\X)$, and $\mathcal{M}_3(\X)$, it is natural to perform singular value decomposition (SVD) on $\mathcal{M}_k(\Y)$ to obtain initial estimators for $U_k$: 
	\begin{equation*}
	\hat{U}_k^{(0)} = \SVD_{r_k}(\mathcal{M}_k(\Y)) = \text{the first $r_k$ left singular vectors of }\mathcal{M}_k(\Y).
	\end{equation*}
	In fact, $\hat{U}^{(0)}_k$ is exactly the higher order SVD (HOSVD) estimator introduced by De Lathauwer, De Moor, and Vandewalle \cite{de2000multilinear}. As we will show later, $\hat{U}_k^{(0)}$ serves as a good starting point but not as an optimal estimator for $U_k$.
	\item[Step 2] (Power Iteration) Then one applies power iterations to update the estimations. Given $\hat{U}_2^{(t-1)}, \hat{U}_3^{(t-1)}$, $\Y$ can be denoised via mode-2 and 3 projections: $\Y\times_2 (\hat{U}_2^{(t-1)})^\top \times_3 (\hat{U}_3^{(t-1)})^\top$. As we will illustrate via theoretical analysis, the mode-1 singular subspace of $\X$ is preserved while the amplitude of the noise is highly reduced after such the projection. Thus, for $t = 1, 2, \ldots$, we calculate
	\begin{equation}
	\begin{split}
	\hat{U}_1^{(t)} = & \text{ first $r_1$ left singular vectors of } \mathcal{M}_1(\Y\times_2 (\hat{U}_2^{(t-1)})^\top \times_3 (\hat{U}_3^{(t-1)})^\top),\\
	\hat{U}_2^{(t)} = & \text{ first $r_2$ left singular vectors of } \mathcal{M}_2(\Y\times_1 (\hat{U}_1^{(t)})^\top \times_3 (\hat{U}_3^{(t-1)})^\top),\\
	\hat{U}_3^{(t)} = & \text{ first $r_3$ left singular vectors of } \mathcal{M}_3(\Y\times_1 (\hat{U}_1^{(t)})^\top \times_2 (\hat{U}_2^{(t)})^\top).
	\end{split}
	\end{equation}
	The iteration is stopped when either the increment is no more than the tolerance $\varepsilon$, i.e.,
	\begin{equation}
	\begin{split}
	& \left\|\Y\times_1 (\hat{U}_1^{(t)})^{\top} \times_2 (\hat{U}_2^{(t)})^{\top} \times_3 (\hat{U}_3^{(t)})^{\top}\right\|_{\rm F}\\
	& - \left\|\Y\times_1 (\hat{U}_1^{(t-1)})^{\top} \times_2 (\hat{U}_2^{(t-1)})^{\top} \times_3 (\hat{U}_3^{(t-1)})^{\top}\right\|_{\rm F} \leq \varepsilon, 
	\end{split}
	\end{equation}
	or the maximum number of iterations is reached.
	\item[Step 3] (Projection) With the final estimates $\hat{U}_1, \hat{U}_2, \hat{U}_3$, we estimate $\S$ and $\X$ as
	\begin{equation*}
	\hat{\S} = \Y \times_1 \hat{U}_1^\top \times_2 \hat{U}_2^\top \times_3 \hat{U}_3^\top,\quad \hat{\X} = \hat{\S}\times_1 \hat{U}_1 \times_2 \hat{U}_2 \times_3 \hat{U}_3 = \Y \times_1 P_{\hat{U}_1} \times_2 P_{\hat{U}_2} \times_3 P_{\hat{U}_3}.
	\end{equation*}
\end{enumerate}
The procedure of HOOI is summarized in Algorithm~\ref{alg:strong}. The further generalization to order-4 or higher tensors SVD will be discussed in Section \ref{sec:discuss}.
\begin{algorithm}
	\caption{Higher Order Orthogonal Iteration (HOOI) \cite{de2000best}}
	\label{alg:strong}
	\begin{algorithmic}[1]
		\State Input: $\Y\in \mathbb{R}^{p_1\times p_2\times p_3}$, $(r_1, r_2, r_3)$, increment tolerance $\varepsilon>0$, and maximum number of iterations $t_{\max}$.
		\State Let $t = 0$, initiate via matricization SVDs
		$$\hat{U}_1^{(0)} = \text{SVD}_{r_1}(\mathcal{M}_1(\Y)),\quad \hat{U}_2^{(0)} = \text{SVD}_{r_2}(\mathcal{M}_2(\Y)),\quad \hat{U}_3^{(0)} = \text{SVD}_{r_3}(\mathcal{M}_3(\Y)).$$
		\Repeat
		\State Let $t = t+1$, calculate
			\begin{equation*}
			\hat{U}_1^{(t)} = \SVD_{r_1}\left(\mathcal{M}_1(\Y\times_2 (\hat{U}_2^{(t-1)})^\top \times_3 (\hat{U}_3^{(t-1)})^\top)\right),
			\end{equation*}
			\begin{equation*}
			\hat{U}_2^{(t)} = \SVD_{r_2}\left(\mathcal{M}_2(\Y\times_1 (\hat{U}_1^{(t)})^\top \times_3 (\hat{U}_3^{(t-1)})^\top)\right),
			\end{equation*}
			\begin{equation*}
			\hat{U}_3^{(t)} = \SVD_{r_3}\left(\mathcal{M}_3(\Y\times_1 (\hat{U}_1^{(t)})^\top \times_2 (\hat{U}_2^{(t)})^\top)\right).
			\end{equation*}
		\Until $t = t_{\max}$ or
		\begin{equation*}
		\begin{split}
		& \left\|\Y\times_1 (\hat{U}_1^{(t)})^{\top} \times_2 (\hat{U}_2^{(t)})^{\top} \times_3 (\hat{U}_3^{(t)})^{\top}\right\|_{\rm F}\\ 
		& - \left\|\Y\times_1 (\hat{U}_1^{(t-1)})^{\top} \times_2 (\hat{U}_2^{(t-1)})^{\top} \times_3 (\hat{U}_3^{(t-1)})^{\top}\right\|_{\rm F} \leq \varepsilon.
		\end{split}
		\end{equation*}
		\State Estimate and output:
		$$\hat{U}_1 = \hat{U}_1^{(t)},\quad \hat{U}_2 = \hat{U}_2^{(t)}, \quad \hat{U}_3 = \hat{U}_3^{(t)};$$
		$$\hat{\X} = \Y\times_1 P_{\hat{U}_1} \times_2 P_{\hat{U}_2} \times_3 P_{\hat{U}_3}. $$
	\end{algorithmic}\label{al:strong-signal}
\end{algorithm}

\section{Statistical Limits: Minimax Upper and Lower Bounds}\label{sec:theory}

In this section, we develop the statistical limits for tensor SVD. Particularly, we analyze the estimation error upper bounds of HOOI and MLE, then develop the corresponding lower bounds. For any $\X\in \mathbb{R}^{p_1\times p_2\times p_3}$, denote $\lambda = \min_{k=1,2,3} \sigma_{r_k}(\mathcal{M}_k(\X))$ as the minimal singular values of each matricization, which essentially measures the signal level in tensor SVD model. Suppose the signal-to-noise ratio (SNR) is $\lambda/\sigma = p^{\alpha}$, where $p = \min\{p_1, p_2, p_3\}$. Then the problem of tensor SVD exhibits three distinct phases: $\alpha\geq 3/4$ (strong SNR), $\alpha < 1/2$ (weak SNR), and $1/2\leq \alpha < 3/4$ (moderate SNR).

We first analyze the statistical performance of HOOI, i.e., Algorithm \ref{al:strong-signal}, under the strong SNR setting that $\lambda/\sigma \geq Cp^{3/4}$. 
\begin{Theorem}[Upper Bound for HOOI]\label{th:upper_bound_strong}
	Suppose there exist constants $C_0, c_0>0$ such that $p_k \leq C_0p$, $\|\X\|_{\rm F} \leq  C_0\sigma\exp(c_0p)$, $r_k \leq C_0p^{1/2}$ for $p = \min\{p_1, p_2, p_3\}$, and $k=1,2,3$. Then there exist absolute constants $C_{gap}, C>0$, which do not depend on $p_k, r_k, \lambda, \sigma, q$, such that whenever 
	$$\lambda/\sigma\geq C_{gap} p^{3/4}, \quad \text{(i.e., in the strong SNR case)},$$ 
	after at most $t_{\max} = C\left(\log\left( \frac{p}{\lambda}\right)\vee 1\right)$ iterations in Algorithm \ref{al:strong-signal}, 	the following upper bounds hold,
	\begin{equation}\label{ineq:upper_bound_U}
	\mathbb{E} r_k^{-1/q}\left\|\sin\Theta\left(\hat{U}_k, U_k\right)\right\|_q \leq C \frac{\sqrt{p_k}}{\lambda/\sigma}, \quad k=1,2,3,\quad 1\leq q \leq \infty,
	\end{equation}
	\begin{equation}\label{ineq:upper_bound_X} 
	\begin{split}
	& \mathbb{E}\left\|\hat{\X} - \X\right\|_{\rm F}^2 \leq C\sigma^2\left(p_1r_1 + p_2r_2 + p_3r_3\right), \quad  \mathbb{E}\frac{\|\hat{\X} - \X\|_{\rm F}^2}{\|\X\|_{\rm F}^2} \leq C\left(\frac{\left(p_1 + p_2 + p_3\right)}{\lambda^2/\sigma^2}\bigwedge 1\right).
	\end{split}
	\end{equation}
\end{Theorem}
\begin{Remark}
	In contrast to the error bound for final estimators $\hat{U}_k$ in \eqref{ineq:upper_bound_U}, an intermediate step in the proof for Theorem \ref{th:lower_bound} yields the following upper bound for initializations $\hat{U}_k^{(0)}$, i.e., the output from Algorithm \ref{al:strong-signal} Step 1,
	\begin{equation}\label{ineq:upper_bound_initialization}
	\begin{split}
	& \mathbb{E}r_k^{-1/q}\left\|\sin\Theta\left(\hat{U}_k^{(0)}, U_k^{(0)}\right)\right\|_q \leq C\frac{\sqrt{p_k}}{\lambda/\sigma} + \frac{ Cp^{3/2}}{\lambda^2/\sigma^2},\quad  k=1,2, 3.
	\end{split}
	\end{equation}
	Compared to Theorem \ref{th:upper_bound_strong}, the bound in \eqref{ineq:upper_bound_initialization} is suboptimal as long as $\lambda/\sigma=p^{\alpha}$ when $3/4 \leq \alpha< 1$. Thus, the higher order SVD (HOSVD) $\hat{U}_k^{(0)}$  \cite{de2000multilinear} may yield sub-optimal result. We will further illustrate this phenomenon by numerical analysis in Section \ref{sec:simulation}. 
\end{Remark}
\begin{Remark}
	Especially when $r=1$, Theorem \ref{th:upper_bound_strong} confirms the conjecture in Richard and Montanari \cite{richard2014statistical} that the tensor unfolding method achieves reliable estimates for order-3 spiked tensors if $\lambda/\sigma > \Omega(p^{3/4})$. Moreover, Theorem \ref{th:upper_bound_strong} further shows the power iterations are necessary in order to refine the reliable estimates to minimax-optimal estimates.

	Our result in Theorem \ref{th:upper_bound_strong} outperforms the ones by Sum-of-Squares (SOS) scheme (see, e.g., \cite{hopkins2015tensor,anandkumar2016homotopy}), where an additional logarithm factor on the assumption of $\lambda$ is required. In addition, the method we analyze here, i.e., HOOI, is efficient, easy to implement, and achieves the optimal rate of convergence for estimation error.
\end{Remark}
\begin{Remark}
	The strong SNR assumption ($\lambda/\sigma \geq Cp^{3/4}$) is crucial to guarantee the performance of Algorithm \ref{al:strong-signal}. Actually, to ensure that Step 1 in Algorithm \ref{al:strong-signal} provides meaningful initializations, $\lambda$ should be at least of order $p^{3/4}$ according to our theoretical analysis.
\end{Remark}

Moreover, the estimators with high likelihood, such as MLE, achieve the following upper bounds under the weaker assumption that $\lambda/\sigma \geq Cp^{1/2}$. 
\begin{Theorem}[Upper Bound for Estimators with Large Likilihood and MLE]\label{th:upper_bound_mle}
	Suppose there exist constants $C_0, c_0>0$ such that $p_k \leq C_0p$, $r_k \leq C_0p^{1/2}$ for $p=\min\{p_1,p_2, p_3\}$, $\|\X\|_{\rm F} \leq  C_0\sigma \left(\exp(c_0p)\right)$, $\max\{r_1, r_2, r_3\}\leq C_0\min\{r_1, r_2, r_3\}$ for $k=1,2,3$. Suppose $\hat{U}_k^{\bullet}\in \mathbb{O}_{p_k, r_k}$ are estimators satisfying 
	\begin{equation}\label{eq:high-likelihood}
	\min_{\hat{\S}^{\bullet}}\|\hat{\Y} - \hat{\S}^{\bullet}\times_1 \hat{U}_1^{\bullet} \times\hat{U}_2^{\bullet} \times_3\hat{U}_3^{\bullet}\|_F^2 \leq \min_{\hat{\S}}\|\hat{\Y} - \hat{\S} \times_1 U_1 \times U_2 \times_3 U_3 \|_F^2,
	\end{equation}
	i.e., the likelihood value of $\hat{U}_k^{\bullet}$ is no less than $U_k$. Then there exists a uniform constant $C_{gap}>0$ (which does not depend on $p_k, r_k, \lambda, \sigma, q$) such that whenever 
	$$\lambda/\sigma \geq C_{gap}p^{1/2}, \quad \text{(i.e., in moderate or strong SNR cases)},$$ 
	$\hat{U}_1^{\bullet}, \hat{U}_2^{\bullet}, \hat{U}_3^{\bullet}$, and $\hat{\X}^{\bullet} = \hat{\S}^{\bullet}\times_1 \hat{U}_1^{\bullet} \times_2 \hat{U}_2^{\bullet} \times_3 \hat{U}_3^\bullet$ satisfy 
	\begin{equation}\label{ineq:MLE-upper-bound}
	\begin{split}
	& \mathbb{E} r_k^{1/q}\left\|\sin\Theta(\hat{U}_k^\bullet, U_k)\right\|_q \leq \frac{C\sqrt{p_k}}{\lambda/\sigma}, \quad k = 1, 2, 3, \quad 1\leq q \leq 2,\\
	& \mathbb{E}\left\|\hat{\X}^\bullet - \X\right\|_{\rm F}^2 \leq C\sigma^2\left(p_1r_1 + p_2r_2 + p_3r_3\right),\\
	& \mathbb{E}\frac{\|\hat{\X}^\bullet - \X\|_{\rm F}^2}{\|\X\|_{\rm F}^2} \leq C\left(\frac{p_1 + p_2 + p_3}{\lambda^2/\sigma^2} \bigwedge 1\right).
	\end{split}
	\end{equation}
	Especially, the upper bounds of \eqref{ineq:MLE-upper-bound} hold for maximum likelihood estimators  \eqref{eq:mle}.
\end{Theorem}

Then we establish the lower bound for tensor SVD. We especially consider the following class of general low-rank tensors,
\begin{equation}\label{eq:lower_rank_class}
\mathcal{F}_{\bp, \br}(\lambda) = \left\{\X\in \mathbb{R}^{p_1\times p_2\times p_3}: \rank_k(\X) \leq r_k, \sigma_{r_k}\left(\mathcal{M}_k(\X)\right)\geq \lambda, k=1,2,3\right\}.
\end{equation}
Here $\bp = (p_1, p_2, p_3)$, $\br = (r_1, r_2, r_3)$ represent the dimension and rank triplets, $\lambda$ is the smallest non-zero singular value for each matricization of $\X$, which essentially measures the signal strength of the problem. The following lower bound holds over $\mathcal{F}_{\bp, \br}(\lambda)$.
\begin{Theorem}[Lower Bound]\label{th:lower_bound}
	Suppose $p=\min\{p_1, p_2, p_3\}$, $\max\{p_1,p_2,p_3\}\leq C_0p$, $\max\{r_1,r_2,r_3\}\leq C_0\min\{r_1,r_2,r_3\}$, $4r_1 \leq r_2r_3, 4r_2 \leq r_3r_1, 4r_3 \leq r_1r_2$, $1\leq r_k \leq p_k/3$, and $\lambda > 0$, 
	then there exists a universal constant $c>0$ such that for $1\leq q\leq \infty$,
	\begin{equation}\label{eq:stronglowerbound1}
	\begin{split}
	& \inf_{\tilde{U}_k}\sup_{\X \in \mathcal{F}_{\bp, \br}(\lambda)} \mathbb{E}r_k^{-1/q}\left\|\sin\Theta\left(\tilde{U}_k, U_k\right)\right\|_q \geq c\left(\frac{\sqrt{p_k}}{\lambda/\sigma} \bigwedge 1\right), \quad k=1,2,3,
	\end{split}
	\end{equation}
	\begin{equation}\label{eq:stronglowerbound2}
	\begin{split}
	& \inf_{\hat{\X}} \sup_{\X \in \mathcal{F}_{\bp, \br}(\lambda)} \mathbb{E} \left\|\hat{\X} - \X\right\|^2_{\rm F} \geq c\sigma^2\left(p_1r_1+p_2r_2+p_3r_3\right),\\
	& \inf_{\hat{\X}} \sup_{\X \in \mathcal{F}_{\bp, \br}(\lambda)} \mathbb{E} \frac{\|\hat{\X} - \X\|^2_{\rm F}}{\|\X\|_{\rm F}^2} \geq c\left(\frac{p_1+p_2+p_3}{\lambda^2/\sigma^2} \bigwedge 1\right).
	\end{split}
	\end{equation}
\end{Theorem}
\begin{Remark}
	Theorem \ref{th:lower_bound} implies when $\lambda/\sigma \leq cp^{1/2}$ for some small constant $c>0$, i.e., under the weak SNR setting, the constant term dominates in \eqref{eq:stronglowerbound1} and there are no consistent estimates for $U_1, U_2, U_3$. On the other hand, when $\lambda/\sigma \geq Cp^{1/2}$, i.e., under the strong and moderate SNR settings, $\frac{\sqrt{p_k}}{\lambda/\sigma}$ dominates in \eqref{eq:stronglowerbound1} and provides non-trivial minimax lower bounds for the estimation errors.
\end{Remark}

We further define $\tau^2 = \mathbb{E}\|\Z\|_F^2$ as the expected squared Frobenius norm of the whole noisy tensor. In summary, Theorems \ref{th:upper_bound_strong}, \ref{th:upper_bound_mle}, and \ref{th:lower_bound} together yield the following statistical limits for tensor SVD.
\begin{enumerate}
	\item Under the strong SNR case that $\lambda/\sigma \geq Cp^{3/4}$ (or $\lambda/\tau\geq Cp^{-3/4}$), the higher order orthogonal iteration, i.e., Algorithm \ref{al:strong-signal}, provides minimax rate-optimal estimators for $U_1, U_2, U_3$, and $\X$.
	\begin{equation}\label{eq:optimal-rate}
	\begin{split}
	& \inf_{\hat{U}_k} \sup_{\X \in \mathcal{F}_{\bp, \br}(\lambda)}\mathbb{E} r_k^{-1/q} \left\|\sin\Theta(\hat{U}_k, U_k)\right\|_q \asymp \frac{\sqrt{p_k}}{\lambda/\sigma} \asymp \frac{\tau\sqrt{p_k}}{\lambda \sqrt{p_1p_2p_3}},\quad k=1,2,3,1\leq q\leq +\infty,\\
	& \inf_{\hat{\X}}\sup_{\X \in \mathcal{F}_{\bp, \br}(\lambda)}\mathbb{E}\left\|\hat{\X} - \X\right\|_{\rm F}^2 \asymp \sigma^2\left(p_1r_1 + p_2r_2 + p_3r_3\right) \asymp \frac{\tau^2(p_1r_1+p_2r_2+p_3r_3)}{p_1p_2p_3}, \\ 
	& \inf_{\hat{\X}}\sup_{\X \in \mathcal{F}_{\bp, \br}(\lambda)}\mathbb{E}\frac{\|\hat{\X} - \X\|_{\rm F}^2}{\|\X\|_{\rm F}^2} \asymp \left(\frac{p_1 + p_2 + p_3}{\lambda^2/\sigma^2} \wedge 1\right) \asymp \left(\frac{\tau^2(p_1+p_2+p_3)}{\lambda^2p_1p_2p_3}\wedge 1\right).
	\end{split}
	\end{equation}	
	\item Under the moderate SNR case that $Cp^{1/2} \leq \lambda/\sigma \leq cp^{3/4}$ (or $Cp^{-1}\leq \lambda/\tau \leq cp^{-3/4}$), the estimators with high likelihood \eqref{eq:high-likelihood}, including the MLE \eqref{eq:mle}, are minimax rate-optimal. The rate here is exactly the same as \eqref{eq:optimal-rate}.
	\item Under the weak SNR case that $\lambda/\sigma \leq c p^{1/2}$ (or $\lambda/\tau \leq cp^{-1}$), there are no consistent estimators for $U_1, U_2, U_3$, or $\X$.
\end{enumerate}

However, as we have discussed in Section \ref{sec:mle-and-spectral}, MLE is not applicable even with the moderate dimension. It is still crucial to know whether there is any fast and efficient algorithm for tensor SVD under the moderate SNR setting.

\section{Computational Limits in Moderate SNR Case}\label{sec:moderate}

In this section, we focus on the computational aspect of tensor SVD under the moderate SNR setting. If $\lambda/\sigma = p^{\alpha}$ with $p=\min\{p_1, p_2, p_3\}$ and $\alpha<3/4$, we develop the computational lower bound to show that every polynomial-time algorithm is statistically inconsistent in estimating $U_1$, $U_2$, $U_3$, and $\X$ based on the computational hardness assumption.

In recent literature, we have seen achievements in obtaining computational lower bounds via computational hardness assumptions for many problems, such as sparse PCA \citep{berthet2013computational,berthet2013optimal,wang2014statistical,gao2014sparse}, submatrix localization \citep{ma2015computational,cai2015computational,chen2014statistical}, tensor completion \citep{barak2016noisy}, sparse CCA \citep{gao2014sparse}, and community detection \citep{hajek2015computational}. The computational hardness assumptions, such as planted clique detection and Boolean satisfiability, has been widely studied and conjectured that no polynomial-time algorithm exists under certain settings. For tensor SVD, our computational lower bound is established upon the hardness hypothesis of hypergraphic planted clique detection, which is discussed in detail in the next section.

\subsection{Planted clique detection in hypergraphs}
Let $G=(V,E)$ be a graph, where $V=\{1,2,\ldots, N\}$ and $E$ are the vertex and edge sets, respectively. For a standard graph, the edge $e=(i,j)\in E$ indicates certain relation exists between vertices $i$ and $j$ in $V$. A 3-hypergraph (or simply noted as a hypergraph, without causing any confusion) is a natural extension, where each hyper-edge is represented by an unordered group of three different vertices, say $e = (i,j,k)\in E$. Given a hypergraph $G=(V,E)$ with $|V|=N$, its adjacency tensor ${\bf A}\in\{0, 1\}^{N\times N\times N}$ is defined as 
$$
A_{ijk}=
\begin{cases}
1,&\textrm{if } (i, j, k)\in E;\\
0,&\textrm{otherwise}.
\end{cases}
$$
We denote the Erd\H{o}s-R\'enyi hypergraph of $N$ vertices as $\mathcal{G}_3(N,1/2)$, if for each $1\leq i<j<k\leq N$, $(i,j,k)$ is included into the hyper-edge set independently with probability $1/2$. For $V_1\subset V$ and certain integer $1\leq \kappa_N\leq |V_1|$, we use $\mathcal{G}_3(N,1/2,\kappa_N,V_1)$ to denote a random hypergraph where a clique of size $\kappa_N$ is planted inside $V_1$. More precisely, we first sample a random graph from $\mathcal{G}_3(N,1/2)$, then pick $\kappa_N$ vertices uniformly at random from $V_1$, denote them as $C$, and connecting all hyper-edges $(i, j, k)$ for all distinct triplets $i,j,k \in C$.
Conventionally, the planted clique detection is referred to as the problem for distinguishing whether there is any planted clique hidden in the Erd\H{o}s-R\'enyi graph. To simplify our analysis in tensor SVD later, we propose a slightly different version of hypergraphic planted clique detection problem as follows.
\begin{Definition}
Let $G$ be drawn from either $\mathcal{G}_3(N,1/2,\kappa_N,V_1)$ or $\mathcal{G}_3(N,1/2,\kappa_N,V_2)$, 
where $V_1=\{1,2,\ldots,\floor{N/2}\}$ and $V_2=\big\{\floor{N/2}+1,\floor{N/2}+2,\ldots,N\big\}$.
The {\it hypergraphic planted clique detection problem}, noted as $ {\bf PC}_3(N,\kappa_N)$, refers to the hypothesis testing problem
	\begin{equation}\label{eq:hypergraph}
	H_0: G\sim \mathcal{G}_3(N,1/2,\kappa_N,V_1)\quad \textrm{vs}.\quad H_1: G\sim \mathcal{G}_3(N,1/2,\kappa_N,V_2).
	\end{equation}
Given a hypergraph $G$ sampled from either $H_0$ or $H_1$ with adjacency tensor ${\bf A}\in\{0,1\}^{N\times N\times N}$, let $\psi(\cdot): \{0, 1\}^{N\times N\times N}\mapsto \{0,1\}$ be a binary-valued function on ${\bf A}$ such that $\psi({\bf A})=1$ indicates rejection of $H_0$. Then the risk of test $\psi$ is defined as the sum of Type-I and II errors,
$$
\mathcal{R}_{N,\kappa_N}(\psi)=\mathbb{P}_{H_0}\big\{\psi({\bf A})=1\big\}+\mathbb{P}_{H_1}\big\{\psi({\bf A})=0\big\}.
$$
\end{Definition}
Putting it differently, given a random hypergraph $G\sim H_0$ or $H_1$, our goal is to identify whether the clique is planted in the first or second half of vertices. 

When we replace the hyper-edges (involving three vertices each) of $\mathcal{G}_3(N, 1/2, \kappa_N, V_1)$ by the regular edges (involving two vertices each), the above hypergraphic planted clique detection becomes the traditional planted clique detection problem. 
To provide circumstantial evidence to the hardness of ${\bf PC}_3(N,\kappa_N)$, it is helpful for us to review some well-known results of the traditional planted clique detection here. First, the difficulty of traditional planted clique detection depends crucially on the planted clique size: $\kappa_N$.  \cite{bollobas1976cliques} and \cite{feldman2013statistical} showed that if $\kappa_N =o\big(\log N\big)$, it is statistically impossible to determine whether a planted clique exists since a random graph $G\sim \mathcal{G}_2(N,1/2)$ contains a clique of size $2\log N$ with high probability. When $\kappa_N\geq C\sqrt{N}$, it has been shown that the planted clique can be located by performing polynomial-time operations by spectral methods \citep{alon1998finding,ames2011nuclear}. If the size clique further increases, say $\kappa_N\geq C\sqrt{N\log N}$, \cite{kuvcera1995expected} developed an algorithm to find exactly the planted clique with high probability in polynomial time. However, when $\log N \ll \kappa_N \ll \sqrt{N}$, there is still no known polynomial-time algorithm for planted clique detection, and it is currently widely conjectured by the theoretical computer science and graph theory community that such polynomial-time algorithm may not exist (see \cite{feldman2013statistical,jerrum1992large,feige2003probable} and the references therein).

When moving to hypergraphs, the hardness of ${\bf PC}_3(N,\kappa_N)$, to the best of our knowledge, remains unclear. In an extreme case of exhaustive search, it needs an exponential number of operations, i.e., $\binom{N}{\kappa_N}$, to verify a solution. In addition, the performance of the simple matricization-spectral method (which shares similar idea as the proposed Algorithm \ref{al:strong-signal})
highly depends on the size of the clique $\kappa_N$. We particularly have the following Proposition \ref{pr:hard}.

\begin{Proposition}\label{pr:hard}
	 Suppose $G\sim \mathcal{G}_3(N, 1/2, \kappa_N, V_1)$, so there exists $C\subseteq V_1$ as a planted clique of size $\kappa_N$ with uniform random position. Let $\A$ be the corresponding adjacency tensor, and
	 $1_C \in \mathbb{R}^{|V_1|}$ be the indicator for the hidden clique that $(1_C)_i = 1_{\{i\in C\}}$. We further partition $V_1 = \{1,\ldots, \lfloor N/2\rfloor \}$ into three equal subsets: $D_k = \{\lfloor kN/6\rfloor + 1,\ldots, \lfloor (k+1)N/6\rfloor\}$ for $k=1,2,3$.
	 Then we can calculate $\hat{u}_k\in \mathbb{R}^{|V_1|}$ as the leading left singular vector of $\mathcal{M}_k(2\cdot\A_{[D_1, D_2, D_3]} - 1_{|D_1|\times |D_2|\times |D_3|})$, where $1_{|D_1|\times |D_2|\times |D_3|}$ is a $|D_1|$-by-$|D_2|$-by-$|D_3|$ tensor with all entries 1. If the sequence $\{\kappa_N\}$ satisfies $\liminf_{N\to \infty} \frac{\kappa_N}{N^{1/2}} = \infty$, then 
	 $$\sin\Theta\left(\hat{u}_k, (1_{C})_{D_k}\right) \overset{d}{\to} 0,\quad  \text{as } N\to \infty, \quad k=1,2,3.$$ 
	 In another word, the angle between $\hat{u}_k$ and $(1_{C})_{D_k}$ tends to 0 in probability.
\end{Proposition}
\begin{Remark}
	For technical convenience, we partition $V_1$ into three parts and perform SVD on $\mathcal{M}_k(2\A_{[D_1, D_2, D_3]}-1_{|D_1|\times |D_2|\times |D_3|})$ to ensure that most of the entries of $\mathcal{M}_k(\A)$ are i.i.d. Rademacher distributed.
\end{Remark}

Proposition \ref{pr:hard} suggests that $\hat{u}_k$ can be used to locate $C$ when $\kappa_N \gg N^{1/2}$. However, the theoretical analysis in Proposition \ref{pr:hard} fails when $\kappa_N = N^{(1-\tau)/2}$ for $\tau>0$, and we conjecture that such computational barrier is essential. Particularly, we propose the following computational hardness assumption on hypergraphic planted clique detection.
\begin{Hypothesis}{\bf H($\tau$).}
	For any sequence $\{\kappa_N\}$ such that $\underset{N\to\infty}{\lim}\sup\frac{\log \kappa_N}{\log \sqrt{N}}\leq (1-\tau)$ and any sequence of polynomial-time tests $\{\psi_N\}$,
	$$
	\underset{N\to\infty}{\lim\inf}\ \mathcal{R}_{N,\kappa_N}(\psi_N)\geq \frac{1}{2}.
	$$

\end{Hypothesis}
\subsection{The computational lower bound of tensor SVD}

Now we are ready to develop the computational lower bound for tensor SVD based on Hypothesis {\bf H($\tau$)}. Recall
$$
{\bf Y}={\bf X}+{\bf Z}\in\mathbb{R}^{p_1\times p_2\times p_3},\quad \X = \S\times_1 U_1\times_2 U_2\times_3 U_3,\quad  {\bf Z}\overset{iid}{\sim} N(0,\sigma^2).
$$
To better present the asymptotic argument, we add a superscript of dimension, $p = \min\{p_1, p_2, p_3\}$, to the estimators, i.e., $\hat{U}_k^{(p)}$, $\hat{\X}^{(p)}$. The computational lower bound is then presented as below. 
\begin{Theorem}[Computational Lower Bound]\label{th:moderatelowerbound}
	Suppose the hypergraphic planted clique assumption ${\bf H(\tau)}$ holds for some $\tau \in (0, 1)$. Then there exist absolute constants $c_0,c_1>0$ such that if $\lambda/\sigma\leq c_0\Big(\frac{p^{3(1-\tau)/4}}{\sqrt{\log 3p}}\Big)$, for any integers $r_1, r_2, r_3 \geq1$ and any polynomial time estimators $\hat{U}_k^{(p)}$, $\hat\X^{(p)}$, the following inequalities hold
	\begin{equation}\label{eq:moderate1}
	\underset{p\to\infty}{\lim\inf} \underset{{\bf X}\in\mathcal{F}_{\bp, \br}(\lambda)}{\sup}\ \mathbb{E}\Big\|\sin\Theta\big(\hat{U}_k^{(p)}, U_k\big)\Big\|^2\geq c_1,\quad k=1,2,3,
	\end{equation}
	\begin{equation}\label{eq:moderate3}
	\underset{p\to\infty}{\lim\inf} \underset{{\bf X}\in\mathcal{F}_{\bp,\br}(\lambda)}{\sup}\ \frac{\mathbb{E}\|\hat{\X}^{(p)}-\X\|_{\rm F}^2}{\|\X\|_{\rm F}^2}\geq c_1.
	\end{equation}
	\end{Theorem}

\begin{Remark}
For technical reasons, there is an additional logarithmic factor in the condition $\lambda/\sigma\leq c_0\Big(\frac{p^{3(1-\tau)/4}}{\sqrt{\log 3p}}\Big)$, compared with the statistical lower bound in Theorem \ref{th:lower_bound}. Since $\tau$ is a strictly positive number, the effect of logarithmic factor is dominated by $p^{c}$ for any $c>0$ asymptotically. 
\end{Remark}
Theorem \ref{th:moderatelowerbound} illustrates the computational hardness for tensor SVD under moderate scenario that $\lambda/\sigma=p^\alpha, 1/2\leq \alpha<3/4$, if the hypergraphic planted clique assumption $\mathbf{H}(\tau)$ holds for any $\tau>0$.

\section{Simulations}\label{sec:simulation}
In this section, we further illustrate the statistical and computational limits for tensor SVD via numerical studies. 

We first consider the average Schatten $q$-$\sin\Theta$-norm losses for initial estimators $\hat{U}_k^{(0)}$ (HOSVD) and final estimators $\hat{U}_k$ (HOOI) under the following simulation setting. For any given triplet $(p, r, \lambda)$, we let $p = p_1= p_2= p_3, r = r_1= r_2= r_3$, generate $\tilde{U}_k\in \mathbb{R}^{p_k\times r_k}$ as i.i.d. standard Gaussian matrices, and apply QR decomposition on $\tilde{U}_k$ and assign the Q part to $U_k$. In other words, the singular subspaces $U_1, U_2, U_3$ are drawn randomly from Haar measure. Then we construct $\tilde{\S}\in \mathbb{R}^{r_1\times r_2\times r_3}$ as an i.i.d. Gaussian tensor, and rescale it as $\S = \tilde{\S}\cdot \frac{\lambda}{\min_{k=1,2,3} \sigma_{r_k}(\mathcal{M}_k(\tilde{\S}))}$ to ensure $\min_{k=1,2,3}\sigma_{r_k}\left(\mathcal{M}_k(\X)\right) \geq \lambda$. Next, we construct $\Y = \X + \Z$, where the signal tensor $\X = \S\times_1 U_1\times_2 U_2 \times_3 U_3$, the noisy tensor $\Z$ are drawn from i.i.d. standard Gaussian distribution. We apply Algorithm \ref{al:strong-signal} to $\Y$ and record the average numerical performance for different values of $(p, r, \lambda)$. The results based on 100 replications are shown in Table \ref{table:setting1}. We can clearly see that the power iterations (Step 2 in Algorithm \ref{al:strong-signal}, i.e., HOOI) significantly improve upon spectral initializations (Step 1 in Algorithm \ref{al:strong-signal}, i.e., HOSVD) in different Schatten-$q$ $\sin\Theta$ losses under different settings.
\begin{table}
	\resizebox{\textwidth}{!}{
\begin{tabular}{r|cc|cc|cc|cc}
\hline\hline
$(p, r,\lambda)$&$l_1(\hat{U})$ & $l_1(\hat{U}^{(0)})$ & $l_2(\hat{U})$ &$l_2(\hat{U}^{(0)})$& $l_5(\hat{U})$ & $l_5(\hat{U}^{(0)})$ & $l_{\infty}(\hat{U})$ & $l_{\infty}(\hat{U}^{(0)})$\\
\hline
  (50, 5, 20) &  1.1094  &  2.1192  &  0.5194  &  1.0535 &   0.3572 &   0.7991 &   0.3286  &  0.7699  \\
  (50, 5, 50) &  0.4297 &   0.5243 &   0.2016  &  0.2519 &   0.1392  &  0.1815&    0.1283 &   0.1713 \\
  (50, 10, 20) & 2.4529 &   4.5208 &   0.8179  &  1.5674  &  0.4629  &  0.9611 &   0.3955 &   0.8762 \\
  (50, 10, 50) &   0.9111  &  1.1210 &   0.3030  &  0.3771 &   0.1707  &  0.2175 &   0.1452 &   0.1890  \\
  (100, 5, 40) & 0.7952  &  1.5649  &  0.3695 &   0.7707  &  0.2509  &  0.5778 &   0.2294 &   0.5543 \\
  (100,    5,   60) &   0.5301 &   0.8132 &   0.2463 &   0.3938 &   0.1673 &    0.2878 &   0.1530 &   0.2731\\
  (100, 10, 40) &  1.7448 &   3.5371  &  0.5688 &   1.1943 &   0.3087 &   0.7015  &  0.2554  &  0.6246 \\
 (100, 10, 60) &     1.1466  &   1.8055 &    0.3735 &   0.6015  &  0.2021 &   0.3427  &  0.1660 &   0.2975\\
\hline\hline
\end{tabular}}
\caption{The average Schatten $q$-$\sin\Theta$ loss of the final estimations $\hat{U}_k$ and the spectral initializations $\hat{U}_k^{(0)}$ based on $100$ repetitions. Here, $p_1=p_2=p_3=p$, $r_1=r_2=r_3=r$, $l_q(\hat{U}) = \frac{1}{3}\sum_{k=1}^3 \|\sin\Theta(\hat{U}_k, U_k)\|_q$.}
\label{table:setting1}
\end{table}

Then we consider another setting that $\X$ has three different dimensions. Specifically, we generate $\Y = \X + \Z$ by the same scheme as the previous setting with varying $(p_1, p_2, p_3)$ and fixed $r_1 =r_2 = r_3 = 5$. We repeat the experiment for 100 times, then record the average estimation errors in Table \ref{table:setting2}. Again, we can see HOOI performs well under various values of dimensions.
\begin{table}
	\resizebox{\textwidth}{!}{%
	\begin{tabular}{r|cc|cc|cc|cc}
		\hline\hline
		$(p_1, p_2, p_3, \lambda)$ & $l_\infty(\hat{U}_1)$ & $l_2(\hat{U}_1)$ & $l_\infty(\hat{U}_2)$ & $l_2(\hat{U}_2)$ & $l_\infty(\hat{U}_3)$ & $l_2(\hat{U}_3)$ & $\|\hat{\X} - \X\|_{\rm F}$ & $\frac{\|\hat{\X} - \X\|_{\rm F}}{\|\X\|_{\rm F}}$\\
		\hline
		$(20, 30, 50, 20)$ &  0.2082  &  0.3032 & 0.2530 & 0.3858 & 0.3109  &0.4975 & 24.7037 &  0.3276\\
		$(20, 30, 50, 100)$ & 0.0409  &  0.0596  &  0.0498  &  0.0761  &  0.0641 &   0.1017 &  23.5708  &  0.0631\\
		$(30, 50, 100, 20)$ & 0.2674  &  0.4036  &  0.3354  &  0.5247 &   0.4456 &   0.7252 &  33.6219 &   0.4479\\
		$(30, 50, 100, 100)$ & 0.0490 &   0.0753 &   0.0640 &   0.1012 &   0.0911 &  0.1469  & 30.9540  &  0.0822\\
		$(100, 200, 300, 50)$ & 0.1840 &   0.2982 &   0.2551 &   0.4301 &   0.3161  &  0.5155 &   57.8482 &   0.3090\\		
		$(100, 200, 300, 100)$ & 0.0940  &  0.1506 &   0.1259  &   0.2117 & 0.1638 &   0.2627 &  55.9009  &  0.1505\\
		$(200, 300, 400, 50)$ & 0.2579  &  0.4335 &   0.3331 &   0.5523  &  0.3420  &  0.6017 &  72.2912  &  0.4026\\
		$(200, 300, 400, 150)$ & 0.0825  &  0.1389  &  0.1076&    0.1739 &   0.1277  &  0.2024 &  68.0305 &    0.1199\\
		\hline\hline
	\end{tabular}}
\caption{The average spectral and Frobenius $\sin\Theta$ loss for $\hat{U}_1$, $\hat{U}_2$, $\hat{U}_3$ and average Frobenius loss for $\hat{\X}$ under various settings. Here $l_\infty(\hat{U}_k) = \|\sin\Theta(\hat{U}_k, U_k)\|$, $l_2(\hat{U}) = \|\sin\Theta(\hat{U}_k, U_k)\|_{\rm F}$.}
\label{table:setting2}
\end{table}

Next, we illustrate the phase transition phenomenon of tensor SVD. Let $\X=\S\times_1 U_1\times_2 U_2\times_3 U_3$ be a $p$-by-$p$-by-$p$ tensor, where  $U_1,U_2,U_3$ are randomly generated $p$-by-$r$ orthogonal matrices from Haar measure and $\S\in\mathbb{R}^{r\times r\times r}$ is a fixed diagonal tensor such that $S_{i,j,k}= p^{\alpha}\cdot 1_{\{i=j=k\}}, 1\leq i, j, k \leq r$ for $\alpha \in [0.4, 0.9]$. Then $p^{\alpha}$ is the signal strength in our context. The entries of $\Z$ are generated as either i.i.d. $N(0, 1)$ or $\text{Unif}[-\sqrt{3}, \sqrt{3}]$, which are sub-Gaussian, mean 0, and variance 1. To demonstrate the phase transitions at both $p^{3/4}$ and $p^{1/2}$, ideally one wishes to implement both MLE and HOOI. Since MLE, i.e., the best low-rank approximation estimator \eqref{eq:mle}, is computationally intractable, we instead consider the following oracle warm-start HOOI to obtain an approximation for MLE: suppose an oracle provides a warm start as
$$
\hat{U}^{(0){\rm warm}}_{k} =\frac{1}{\sqrt{2}}U_k+\frac{1}{\sqrt{2}}U_{k}',\quad k=1,2,3,
$$
where $U_k$ is the true underlying loading and $U_{k}'$ is a $p$-by-$r$ random orthonormal matrix in the complementary space of $U_k$. $\{U_{k}'\}_{k=1}^3$ here are generated based on the following scheme: first calculate $U_{k\perp}\in\mathbb{O}_{p, p-r}$ as the orthogonal complement of $U_k$, then construct $U_k' = U_{k\perp}O$ for some random orthogonal matrix $O\in \mathbb{O}_{p-r, r}$. Based on the oracle warm-start, we apply Steps 2 and 3 of Algorithm \ref{al:strong-signal} to obtain the warm-start HOOI estimator $\hat{U}_k^{\rm warm}$ as an approximation for MLE.

\begin{figure}
\centering
 \includegraphics[width=0.8\textwidth]{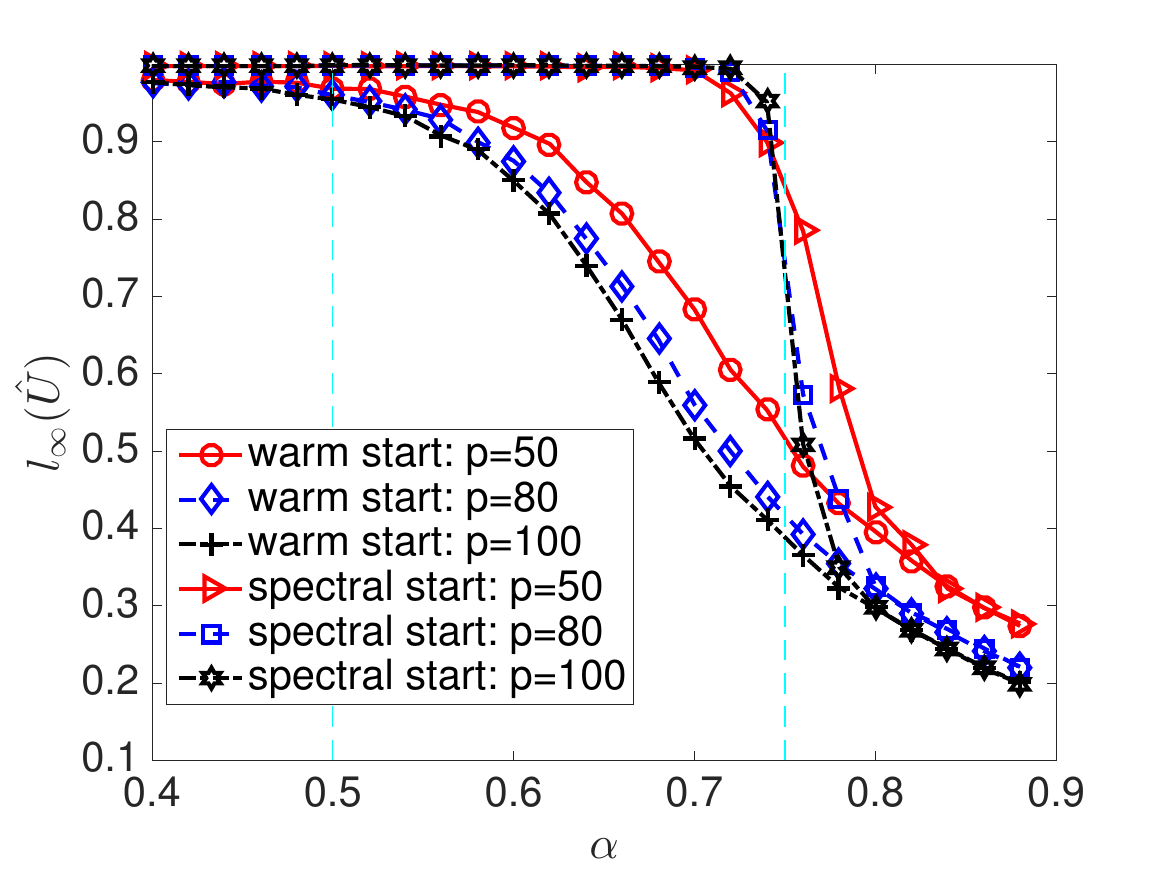}\\
 \includegraphics[width=0.8\textwidth]{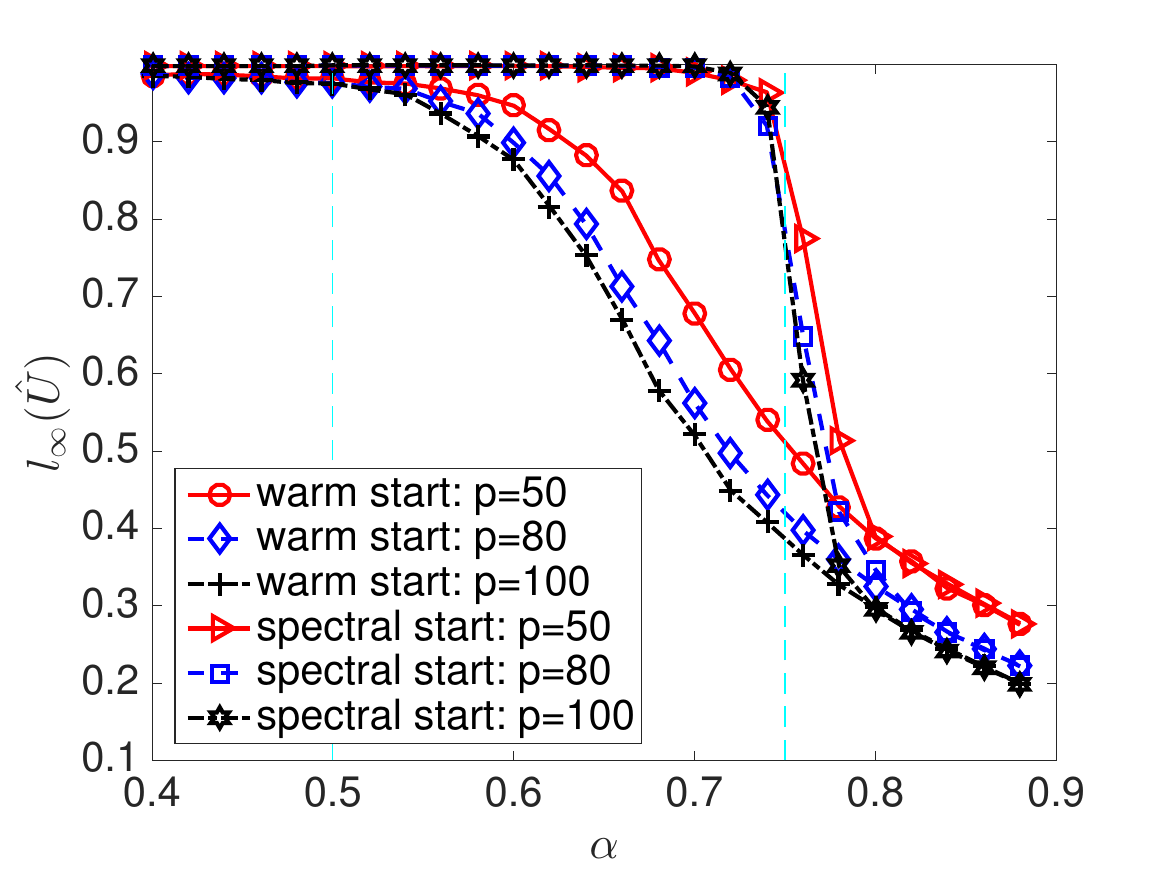}
\caption{Phase transitions in tensor SVD at SNR $ = p^{.5}$ and $ p^{.75}$. Upper panel: Gaussian noise $N(0, 1)$. Lower panel: uniform noise ${\rm Unif}[-\sqrt{3}, \sqrt{3}]$.}\label{fig:simu_1}
\end{figure}

We let $p$ vary from 50 to 100, $r = 5$, and apply both the spectral-start HOOI (i.e., the original HOOI and Algorithm \ref{al:strong-signal}) and the oracle warm-start HOOI. The average spectral $\sin\Theta$ loss, i.e., $l_{\infty}(\hat{U}) = \frac{1}{3}\sum_{k=1}^3\|\sin\Theta(\hat{U}_k, U_k)\|$, from 100 repetitions are presented in Figure \ref{fig:simu_1}, where the upper panel and lower panel correspond to the i.i.d. Gaussian noise and i.i.d. uniform noise cases, respectively. Both panels of Figure \ref{fig:simu_1} clearly demonstrate the phase transition effects: the estimation error significantly decreases around SNR $ = p^{3/4}$ and around SNR $=p^{1/2}$ for spectral-start HOOI and oracle warm-start HOOI, respectively. This exactly matches our theoretical findings in Section \ref{sec:theory}. In addition, there is little difference between two plots in the upper and lower panels, which implies that the statistical estimation error for tensor SVD mainly relies on the SNR and is less influenced by the particular sub-Gaussian noise type.

\section{Discussions: Further Generalizations}\label{sec:discuss}

In this article, we propose a general framework for tensor singular value decomposition (tensor SVD), which focuses on extracting the underlying Tucker low-rank structure from the noisy tensor observations. We provide a comprehensive analysis for tensor SVD in aspects of both statistics and computation. The problem exhibits three distinct phases according to the signal-to-noise ratio (SNR): with strong SNR, the higher order orthogonal iteration (HOOI) performs efficiently and achieves statistical optimal results; with weak SNR, no method performs consistently; with moderate SNR, the estimators with high likelihood, such as the computational intractable MLE, perform optimally in statistical convergence rate, and no polynomial algorithm can do so unless we have a polynomial-time algorithm for the hypergraphic planted clique problem.

The results of this paper are mainly presented under the i.i.d. Gaussian noise setting. When the noise is more generally i.i.d. sub-Gaussian distributed, say
$$\Z \overset{iid}{\sim} Z, \quad \text{where}\quad \|Z\|_{\psi_2} = \sup_{q\geq 1}q^{-1/2}(\mathbb{E} |Z|^q)^{1/q} \leq \sigma,$$
we can derive the upper bounds similarly to Theorems \ref{th:upper_bound_strong} and \ref{th:upper_bound_mle}, as the proofs of main technical tools, including Lemmas \ref{lm:random_tensor_max_projection} and \ref{lm:gaussian-vector-projection}, still hold for the i.i.d. sub-Gaussian noise case.

We have also focused our presentations mainly on order-3 tensors throughout this article. The results can be additionally generalized to order-$d$ tensor SVD for any $d\geq 2$. Suppose one observes an order-$d$ tensor as follows,
\begin{equation}
\Y = \X + \Z, \quad \X = \S\times_1 U_1 \cdots \times_d U_d,
\end{equation}
where $\Y, \X, \Z \in \mathbb{R}^{p_1\times \cdots \times p_d}$, $U_k \in \mathbb{O}_{p_k, r_k}$ for $k=1,\ldots, d$, and $\S\in \mathbb{R}^{r_1\times \cdots \times r_d}$. The higher order orthogonal iteration can be written as follows \cite{de2000best}, 
\begin{enumerate}
	\item[Step 1] Initialize by singular value decomposition of each matricizations,
$$\hat{U}_k^{(0)} = \SVD_{r_k}\left(\mathcal{M}_k(\Y)\right),\quad k=1,\ldots, d. $$
	\item[Step 2] Let $t = 0, 1, \ldots$, update the estimates for $U_k$ sequentially for $k=1,\ldots, d$,
	\begin{equation*}
	\begin{split}
	\hat{U}_k^{(t+1)} = \SVD_{r_k}\bigg(\Y& \times_1 (\hat{U}_1^{(t+1)})^\top\times \cdots \times_{(k-1)} (\hat{U}_{k-1}^{(t+1)})^\top \\
	& \times_{(k+1)} (\hat{U}_{k+1}^{(t)})^\top \times \cdots \times_{d} (\hat{U}_{d}^{(t)})^\top \bigg).
	\end{split}
	\end{equation*}
	The iteration is continued until convergence or maximum number of iteration is reached.
	\item[Step 3] With the final estimators $\{\hat{U}_k\}_{k=1}^d$ from Step 2, one estimates $\X$ as
	$$\hat{\X} = \Y \times_1 P_{\hat{U}_1} \times \cdots \times_d P_{\hat{U}_d}.$$
\end{enumerate}
Meanwhile, the non-convex maximum likelihood estimates can be written as
\begin{equation}
\begin{split}
(\hat{U}_1^{\rm mle},\ldots, \hat{U}_d^{\rm mle}) = & \argmax_{\substack{V_k\in \mathbb{O}_{p_k, r_k}\\k=1,\ldots, d}} \left\|\Y\times_1 V_1^\top \times \cdots \times_d V_d^\top\right\|_{\rm F}, \\
\hat{\X}^{\rm mle} = & \Y\times_1 P_{\hat{U}_1^{\rm mle}}\times \cdots \times_d P_{\hat{U}_d^{\rm mle}}.
\end{split}
\end{equation}
Again, let $\lambda = \min_{1\leq k \leq d} \sigma_{r_k}(\mathcal{M}_k(\X))$ measure the signal strength. For fixed $d$, when $p = \min\{p_1,\ldots, p_d\}$, $\max\{p_1,\ldots, p_d\}\leq Cp$, $r_k \leq Cp^{1/(d-1)}$, similarly to the proofs for Theorems \ref{th:upper_bound_strong}, \ref{th:upper_bound_mle}, and \ref{th:lower_bound}, it is possible to show under strong SNR case, where $\lambda/\sigma = p^{\alpha}$ for $\alpha \geq d/4$, HOOI achieves optimal rate of convergence over the following class of low-rank tensors
\begin{equation*}
\mathcal{F}_{\bp, \br}(\lambda)=\left\{\X \in \mathbb{R}^{p_1\times \cdots \times p_d}: \sigma_{r_k}(\mathcal{M}_k(\X))\geq \lambda, k=1,\ldots, d\right\};
\end{equation*}
under the weak SNR where $\lambda/\sigma = p^{\alpha}$ for $\alpha < 1/2$, it is impossible to have consistent estimators for $U_1,\ldots, U_d$, or $\X$ in general; under the moderate SNR, where $\lambda/\sigma = p^{\alpha}$ for $\frac{1}{2}\leq \alpha < \frac{d}{4}$, the estimators with high likelihood, such as MLE, achieve optimal statistical performance; while one can develop a computational lower bound with the computational hardness assumption of a higher order hypergraphic planted clique detection problem similarly to Theorem \ref{th:moderatelowerbound}. We can also see that the gap between statistical and computational limits vanishes if $d=2$. This coincides with the previous results in matrix denoising literature (see, e.g. \cite{donoho2014minimax,gavish2014optimal,cai2016rate}), where the standard singular value decomposition achieves both statistical optimality and computational efficiency. 

Additionally, if $d$ grows rather than stays as a fixed constant, the asymptotics of tensor SVD in both statistical and computational aspects will be an interesting future project.

\section{Proofs}\label{sec:proof}

We collect the proofs for the main results in this paper in this section. Specifically, the proof for Theorems \ref{th:upper_bound_strong}, \ref{th:lower_bound}, and \ref{th:moderatelowerbound} are presented in Sections \ref{sec:proofstrong}, \ref{sec:proof_lower_bound}, and \ref{sec:proof_moderatelowerbound}, respectively. Proofs for Theorem \ref{th:upper_bound_mle}, Proposition \ref{pr:hard}, and additional technical lemmas are postponed to the supplementary materials.

\subsection{Proof of Theorem \ref{th:upper_bound_strong}}\label{sec:proofstrong}

We first consider the proof for Theorem \ref{th:upper_bound_strong}. Throughout the proof, we assume the noise level $\sigma^2=1$ without loss of generality. For convenience, we denote 
$$X_1 = \mathcal{M}_1(\X), \quad X_2 = \mathcal{M}_2(\X), \quad X_3 = \mathcal{M}_3(\X)$$
as the matricizations of $\X$. We also denote $Y_1, Y_2, Y_3, Z_1, Z_2, Z_3$ in the similar fashion. We also let $r = \max\{r_1, r_2, r_3\}$.

We divide the proof into steps. 
\begin{enumerate}[leftmargin=*]
	\item In this first step, we consider the performance of initialization step, we particularly prove that for any small constant $c_0>0$, there exists large constant $C_{gap} > 0$ such that whenever $\lambda \geq C_{gap} p^{3/4}$, we have
	\begin{equation}\label{ineq:initial_inequality}
	\left\|\sin\Theta(\hat{U}_k^{(0)}, U_k)\right\| \leq c_0\Big(\frac{\sqrt{p_k}\lambda + (p_1p_2p_3)^{1/2}}{\lambda^2}\Big)
	\end{equation}
	with probability at least $1 - C\exp(-cp)$. The proof of this step is closely related to the proof for Theorem 3 in \cite{cai2016rate}. Note that
	\begin{equation*}
	\hat{U}^{(0)}_1 = \SVD_{r_1}\left(Y_1\right), \quad Y_1 = X_1 + Z_1,
	\end{equation*}
	where $X_1$ is a fixed matrix satisfying $\rank(X_1) = r_1$; $Z_1\in \mathbb{R}^{p_1\times (p_2p_3)}, \\ \{(Z_1)_{ij}\}_{i,j=1}^{p_1, p_2p_3} \overset{iid}{\sim} N(0, 1)$. This shares the same setting as the one in Theorem 3 in \cite{cai2016rate}, if one sets $p_1, p_2$ in the statement of Theorem 3 in \cite{cai2016rate} respectively as $p_2p_3, p_1$ in our context. Thus we can essentially follow their proof. Let $U_{1\perp}$ be the orthogonal complement of $U_1$. Then the Appendix Equations (1.15), (1.16) in \cite{cai2016rate} yields
	\begin{equation}\label{ineq:derived_from_perturbation1}
	\begin{split}
	& \mathbb{P}\left(\left\|\sin\Theta(\hat{U}_1^{(0)}, U_1)\right\|^2 \leq \frac{C(\lambda^2 + p_2p_3) \|U_{1\perp}^\top Y_1P_{Y_1^\top U_1}\|^2}{\lambda^4}\right)\\
	\geq & 1 - C\exp\left\{-c\frac{\lambda^4}{\lambda^2 + p_2p_3}\right\}
	\end{split}
	\end{equation}
	and
	\begin{equation}\label{ineq:derived_from_perturbation2}
	\begin{split}
	& \mathbb{P}\left(\|U_{1\perp}^\top Y_1P_{Y_1^\top U_1}\| \geq x\right)\\
	\leq & C\exp\left\{Cp_1 - c\min\left(x^2, x\sqrt{\lambda^2 + p_2p_3}\right)\right\} + C\exp\left\{-c(\lambda^2 + p_2p_3)\right\}
	\end{split}
	\end{equation}
	for some uniform constant $C, c>0$. Since 
	$$\lambda \geq C_{gap} p^{3/4} \geq C_{gap}c\left(p_1^{1/2} + (p_1p_2p_3)^{1/4}\right),$$ 
	if we set $x = C\sqrt{p_1}$, \eqref{ineq:derived_from_perturbation2} further leads to 
	\begin{equation}\label{ineq:derived_from_perturabation3}
	\mathbb{P}\left(\|U_{1\perp}^\top Y_1P_{Y_1^\top U_1}\| \geq C\sqrt{p_1}\right) \leq Ce^{-cp} + C\exp(-c(\lambda^2+p^2))\leq C\exp(-cp).
	\end{equation}
	Combining \eqref{ineq:derived_from_perturabation3} and \eqref{ineq:derived_from_perturbation1}, we have proved \eqref{ineq:initial_inequality} for $k=1$. The proof for \eqref{ineq:initial_inequality} for $k=2,3$ can be similarly written down.

	\item After spectral initialization, we assume the algorithm evolves from $t=0$ to $t=t_{\max}$, where $t_{\max} \geq C\left(\log(\frac{p}{\lambda})\vee 1\right)$. In this step, we derive the perturbation bounds for $\hat{U}_1^{(t_{\max})}, \hat{U}_2^{(t_{\max})}, \hat{U}_3^{(t_{\max})}$ under the assumptions that $\lambda \geq C_{gap}^{3/4}$ for large constant $C_{gap} >0$ and the following inequalities all holds,
	\begin{equation}\label{ineq:initial_performance}
	\max\left\{\left\|\sin\Theta(\hat{U}_1^{(0)}, U_1)\right\|, \left\|\sin\Theta(\hat{U}_2^{(0)}, U_2)\right\|, \left\|\sin\Theta(\hat{U}_3^{(0)}, U_3)\right\|\right\} \leq \frac{1}{2}, 
	\end{equation}
	\begin{equation}\label{ineq:th1-check-1}
	\begin{split}
	& \max_{\substack{V_2\in \mathbb{R}^{p_2\times r_2}\\V_3\in \mathbb{R}^{p_3\times r_3}}}\frac{\left\|Z_1 \cdot \left(V_2\otimes V_3\right) \right\|}{\|V_2\|\|V_3\|} \leq C_1\sqrt{pr}, \quad  \max_{\substack{V_3\in \mathbb{R}^{p_3\times r_3}\\V_1\in \mathbb{R}^{p_1\times r_1}}}\frac{\left\|Z_2 \cdot \left(V_3\otimes V_1\right) \right\|}{ \|V_3\|\|V_1\|} \leq C_1\sqrt{pr},\\
	& \max_{\substack{V_1\in \mathbb{R}^{p_1\times r_1}\\V_2\in \mathbb{R}^{p_2\times r_2}}}\frac{\left\|Z_3 \cdot \left(V_1\otimes V_2\right) \right\|}{ \|V_1\|\|V_2\|} \leq C_1\sqrt{pr},
	\end{split}
	\end{equation}
	\begin{equation}\label{ineq:th1-check-2}
	\begin{split}
	& \left\|Z_1 \left(U_2  \otimes U_3\right) \right\| \leq C_2\sqrt{p_1},  \left\|Z_2 \left(U_3 \otimes U_1\right) \right\| \leq C_2\sqrt{p_2}, \left\|Z_3 \left(U_1 \otimes U_2\right) \right\| \leq C_2\sqrt{p_3}.
	\end{split}
	\end{equation}
	Recall here that $U_1, U_2, U_3$ are left singular subspaces for $X_1, X_2, X_3$, respectively.	We let $L_t$ be the spectral $\sin\Theta$ norm error for $\hat{U}_k^{(t)}$,
	\begin{equation}\label{eq:def_L_t}
	L_t = \max_{k = 1, 2, 3} \left\|\sin\Theta(\hat{U}_k^{(t)}, U_k)\right\|, \quad t = 0, 1, 2, \ldots 
	\end{equation} 
	Given \eqref{ineq:initial_performance}, $L_0 \leq \frac{1}{2}$.
	Next we aim to prove that for $t = 0, 1, \ldots$,
	\begin{equation}\label{ineq:L_t+1<=L_t}
	L_{t+1} = \max_{k=1,2,3} \left\|\sin\Theta(\hat{U}_k^{(t+1)}, U_k)\right\| \leq \frac{C_1\sqrt{pr}}{\lambda} L_{t} + \frac{C_2\sqrt{p}}{\lambda} \leq \frac{1}{2}.
	\end{equation}
	To show \eqref{ineq:L_t+1<=L_t}, we first focus on the upper bound of $\|\sin\Theta(\hat{U}_1^{(t+1)}, U_1)\|$ when $t=0$. Define the following key components in our analysis as follows,
	\begin{equation*}
	\begin{split}
	& Y_1^{(t)} = \mathcal{M}_1\Big(\Y \times_2 \big(\hat{U}_2^{(t)}\big)^{\top} \times_3 \big(\hat{U}_3^{(t)}\big)^{\top}\Big) \overset{\text{Lemma \ref{lm:tensor-algebra}}}{=} Y_1\cdot \left(\hat{U}_2^{(t)}\otimes \hat{U}_3^{(t)}\right) \in \mathbb{R}^{p_1\times r_2r_3}, \\
	& X_1^{(t)} = \mathcal{M}_1\Big(\X \times_2 \big(\hat{U}_2^{(t)}\big)^{\top} \times_3 \big(\hat{U}_3^{(t)}\big)^\top\Big) \overset{\text{Lemma \ref{lm:tensor-algebra}}}{=} X_1\cdot \left(\hat{U}_2^{(t)}\otimes \hat{U}_3^{(t)}\right) \in \mathbb{R}^{p_1\times r_2r_3},\\
	& Z_1^{(t)} = \mathcal{M}_1\Big(\Z \times_2 \big(\hat{U}_2^{(t)}\big)^{\top} \times_3 \big(\hat{U}_3^{(t)}\big)^{\top}\Big) \overset{\text{Lemma \ref{lm:tensor-algebra}}}{=} Z_1\cdot \left(\hat{U}_2^{(t)}\otimes \hat{U}_3^{(t)}\right) \in \mathbb{R}^{p_1\times r_2r_3}. 
	\end{split}
	\end{equation*}
	By definition, the left and right singular subspaces of $X_1$ are $U_1 \in \mathbb{O}_{p_1, r_1}$ and $U_2\otimes U_3\in \mathbb{O}_{p_2p_3, r_2r_3}$.
	Then,
	\begin{equation}\label{ineq:sigma_r_X_1}
	\begin{split}
	& \sigma_{r_1}\left(X_1^{(t)}\right) = \sigma_{r_1}\left(X_1 \cdot \left(\hat{U}_2^{(t)}\otimes \hat{U}_3^{(t)}\right)\right) =  \sigma_{r_1}\left(X_1 \cdot \Proj_{U_2\otimes U_3} \cdot \left(\hat{U}_2^{(t)}\otimes \hat{U}_3^{(t)}\right)\right) \\
	= & \sigma_{r_1}\left(X_1 \cdot \left(U_2\otimes U_3\right) \cdot \left(U_2\otimes U_3\right)^\top \cdot \left(\hat{U}_2^{(t)}\otimes \hat{U}_3^{(t)}\right)\right) \\
	\geq & \sigma_{r_1}\left(X_1 \cdot (U_2\otimes U_3)\right)\cdot \sigma_{\min}\left(\left(U_2\otimes U_3\right)^\top \cdot \left(\hat{U}_2^{(t)}\otimes \hat{U}_3^{(t)}\right)\right)\\
	= & \sigma_{r_1}(X_1) \cdot \sigma_{\min}\left(\left(U_2^\top \hat{U}_2^{(t)}\right)\otimes \left(U_3^\top \hat{U}_3^{(t)}\right)\right)\\
	\geq & \sigma_{r_1}(X_1)\cdot \sigma_{\min}\left(U_2^\top \hat{U}_2^{(t)}\right)\cdot \sigma_{\min}\left(U_3^\top \hat{U}_3^{(t)}\right)\\
	\geq & \lambda \cdot \left(1 - L_t^2\right)\quad \text{(by \eqref{eq:def_L_t} and Lemma 1 in \cite{cai2016rate})}.
	\end{split}
	\end{equation}
	Meanwhile,
	\begin{equation}\label{ineq:sigma_r_Z_1}
	\begin{split}
	& \left\|Z_1^{(t)}\right\| = \left\|Z_1 \left(\hat{U}_2^{(t)}\otimes \hat{U}_3^{(t)}\right)\right\|\\
	= & \left\|Z_1 \left(P_{U_2\otimes U_3} + P_{U_{2\perp}\otimes U_3} + P_{I_p\otimes U_{3\perp}} \right) \left(\hat{U}_2^{(t)}\otimes \hat{U}_3^{(t)} \right) \right\| \quad \text{(by \text{Lemma \ref{lm:tensor-algebra}})}\\
	\leq & \left\|Z_1\left(\Proj_{U_2\otimes U_3}\right)\left(\hat{U}_2^{(t)}\otimes \hat{U}_3^{(t)}\right)\right\| + \left\|Z_1\left(\Proj_{U_{2\perp}\otimes U_3}\right)\left(\hat{U}_2^{(t)}\otimes \hat{U}_3^{(t)}\right)\right\|\\
	& + \left\|Z_1\left(\Proj_{U_2\otimes U_{3\perp}}\right)\left(\hat{U}_2^{(t)}\otimes \hat{U}_3^{(t)}\right)\right\| + \left\|Z_1\left(\Proj_{U_{2\perp}\otimes U_{3\perp}}\right)\left(\hat{U}_2^{(t)}\otimes \hat{U}_3^{(t)}\right)\right\|\\
	= & \left\|Z_1 (U_2\otimes U_3) \left(U_2\otimes U_3\right)^\top \left(\hat{U}_2^{(t)}\otimes \hat{U}_3^{(t)}\right) \right\|\\ 
	& +\left\|Z_1\left((\Proj_{U_{2\perp}}\hat{U}_2^{(t)})\otimes(\Proj_{U_3}\hat{U}_3^{(t)})\right)\right\|\\
	& + \left\|Z_1\left((\Proj_{U_{2}}\hat{U}_2^{(t)})\otimes(\Proj_{U_{3\perp}}\hat{U}_3^{(t)})\right)\right\| + \left\|Z_1\left((\Proj_{U_{2\perp}}\hat{U}_2^{(t)})\otimes(\Proj_{U_{3\perp}}\hat{U}_3^{(t)})\right)\right\|\\
	\overset{\eqref{ineq:th1-check-1}}{\leq} & \|Z_1(U_2\otimes U_3)\| + C_1\sqrt{pr} \left\|P_{U_{2\perp}} \hat{U}_2^{(t)}\right\| \left\|P_{U_{3}} \hat{U}_3^{(t)}\right\| \\
	& + C_1\sqrt{pr} \left\|P_{U_{2}} \hat{U}_2^{(t)}\right\| \left\|P_{U_{3\perp}} \hat{U}_3^{(t)}\right\| + C_1\sqrt{pr} \left\|P_{U_{2\perp}} \hat{U}_2^{(t)}\right\| \left\|P_{U_{3\perp}} \hat{U}_3^{(t)}\right\|\\
	\overset{\eqref{ineq:th1-check-2}}{\leq} & C_2\sqrt{p}_1 + C_1\sqrt{pr} L_t + C_1\sqrt{pr} L_t + C_1\sqrt{pr} L_t^2 
	\\ & \quad\quad  \text{(by \eqref{ineq:th1-check-2} and Lemma 1 in \cite{cai2016rate})}\\
	\leq & C_2\sqrt{p}_1 + 3C_1\sqrt{pr}L_t \quad \text{(since the spectral $\sin\Theta$ norm is at most 1)}.
	\end{split}
	\end{equation}
	Since $U_1$ and $\hat{U}_1^{(t+1)}$ are respectively the leading $r$ singular vectors of $X_1^{(t)}$ and $Y_1^{(t)}$, by Wedin's $\sin\Theta$ theorem \citep{wedin1972perturbation},
	\begin{equation}
	\begin{split}
	& \left\|\sin\Theta\left(\hat{U}^{(t+1)}_1, U_1\right)\right\| \leq \frac{\|Z_1^{(t)}\|}{\sigma_r\big(X_1^{(t)}\big)} \leq \frac{C_2 \sqrt{p_1} + C_1\sqrt{pr}L_t}{\lambda(1-L_t^2)}\\ 
	\overset{\eqref{ineq:initial_performance}}{\leq} & \frac{2C_2\sqrt{p_1}}{\lambda} + \frac{4C_1\sqrt{pr}}{\lambda}L_t.
	\end{split}
	\end{equation}
	We can similarly prove that 
	\begin{equation*}
	\begin{split}
	&\left\|\sin\Theta\left(\hat{U}^{(t+1)}_2, U_2\right)\right\| \leq \frac{2C_2\sqrt{p_2}}{\lambda} + \frac{4C_1\sqrt{pr}}{\lambda}L_t,\\ &\left\|\sin\Theta\left(\hat{U}^{(t+1)}_3, U_3\right)\right\|\leq \frac{2C_2\sqrt{p_3}}{\lambda} + \frac{4C_1\sqrt{pr}}{\lambda}L_t.
	\end{split}
	\end{equation*}
	Finally, since $\max\{p_1, p_2, p_3\}\leq C_0p$ and $\max\{r_1, r_2, r_3\} \leq C_0p^{1/2}$, there exists a large constant $C_{gap} > 0$ such that when $\lambda\geq C_{gap}p^{3/4}$,
	\begin{equation}\label{ineq:1/2-useful}
	\frac{2C_2\sqrt{p_1}}{\lambda} + \frac{4C_1\sqrt{pr}}{\lambda}L_t \leq \frac{1}{2}\quad  \text{and}\quad \frac{4C_1\sqrt{pr}}{\lambda} \leq \frac{1}{2}
	\end{equation}
	Then we have finished the proof for \eqref{ineq:L_t+1<=L_t} for $t = 0$. By induction, we can sequentially prove that \eqref{ineq:L_t+1<=L_t} for all $t \geq 0$.

	At this point, \eqref{ineq:L_t+1<=L_t} yields
	\begin{equation*}
	\begin{split}
	& L_{t+1} \leq \frac{C\sqrt{p}}{\lambda} + \frac{4C_1\sqrt{pr}}{\lambda}L_t, \quad \quad t = 1, 2, \ldots, t_{\max}-1 \\ 
	\Rightarrow \quad  & L_{t+1} - \frac{2C\sqrt{p}}{\lambda} \leq \frac{4C_1\sqrt{pr}}{\lambda}\left(L_t - \frac{2C\sqrt{p}}{\lambda}\right),\quad \text{(since \eqref{ineq:1/2-useful})},\\
	\Rightarrow \quad & L_{t_{\max}} - \frac{2C\sqrt{p}}{\lambda} \leq \left(\frac{4C_1\sqrt{pr}}{\lambda}\right)^{t_{\max}} \cdot \left(L_0 - \frac{2C\sqrt{p}}{\lambda} \right)\\
	\Rightarrow \quad & L_{t_{\max}} \leq \frac{2C\sqrt{p}}{\lambda} +  \frac{L_0}{2^{t_{\max}}}  = \frac{2C\sqrt{p}}{\lambda} +  \frac{1}{2^{t_{\max}}} C\left(\frac{\sqrt{p}\lambda + p^{3/2}}{\lambda^2}\right) \leq \frac{3C\sqrt{p}}{\lambda}
	\end{split}
	\end{equation*}
	when $t_{\max}\geq C\left(\log\left( \frac{p}{\lambda}\right)\vee 1\right)$. Therefore, we have the following upper bound for spectral $\sin\Theta$ norm loss for $\hat{U}_k^{t_{\max}} = \hat{U}_k$,
	\begin{equation}
	\|\sin\Theta(\hat{U}_k, U_k)\|\leq L_{t_{\max}} \leq \frac{C\sqrt{p_k}}{\lambda}.
	\end{equation}	
	when \eqref{ineq:initial_performance}, \eqref{ineq:th1-check-1}, \eqref{ineq:th1-check-2} holds.
	
	By the same calculation, we can also prove $L_{t_{\max}-1}$ satisfies $L_{t_{\max}-1} \leq C\sqrt{p}/\lambda$. We prepare this inequality for the use in the next step.
	
	\item In this step, we develop the upper bound for $\left\|\hat{\X} - \X\right\|_{\rm F}$ under the assumptions of \eqref{ineq:initial_performance}, \eqref{ineq:th1-check-1}, \eqref{ineq:th1-check-2}, and 
	\begin{equation}\label{ineq:th1-check-3}
	\left\|\Z\times_1 \hat{U}_1^\top \times_2 \hat{U}_2^\top \times_3 \hat{U}_3^\top \right\|_{\rm F} \leq C\left(\sqrt{p_1r_1} + \sqrt{p_2r_2} + \sqrt{p_3r_3}\right).
	\end{equation}
	Instead of working on $\|\X\|_{\rm F}$ and $\hat{U}_k^{(t_{\max})}$ directly, we take one step back and work on the evolution of $\hat{U}_k^{(t_{\max}-1)}$ to $\hat{U}_k^{(t_{\max})}$. 
	
	Recall that $\hat{U}_1=\hat{U}_1^{(t_{\max})}, \hat{U}_2=\hat{U}_2^{(t_{\max})}, \hat{U}_3=\hat{U}_3^{(t_{\max})}$; $\hat{U}_{1\perp}$, $\hat{U}_{2\perp}$, $\hat{U}_{3\perp}$ are the orthogonal complements of $\hat{U}_{1}$, $\hat{U}_{2}$, $\hat{U}_{3}$, respectively; $\hat{\X} = \Y \times_1 \Proj_{\hat{U}_1} \times_2 \Proj_{\hat{U}_2}\times_3 \Proj_{\hat{U}_3}$. In the previous step we have also proved that	
	$$\left\|\sin\Theta(\hat{U}_k, U_k)\right\| \leq C\frac{\sqrt{p_k}}{\lambda}, \quad k = 1, 2, 3.$$
	Then we have the following decomposition for the estimation error
	\begin{equation}\label{ineq:hat_X-X-decomposition}
	\begin{split}
	& \left\|\hat{\X} - \X\right\|_{\rm F} \\
	\leq & \left\|\X - \X \times_1 \Proj_{\hat{U}_1} \times_2 \Proj_{\hat{U}_2}\times_3 \Proj_{\hat{U}_3} \right\|_{\rm F} + \left\|\Z \times_1 \Proj_{\hat{U}_1} \times_2 \Proj_{\hat{U}_2}\times_3 \Proj_{\hat{U}_3} \right\|_{\rm F}\\
	= & \left\|\X \times_1 \Proj_{\hat{U}_{1\perp}} + \X\times_1 \Proj_{\hat{U}_1} \times_2 \Proj_{\hat{U}_{2\perp}} + \X\times_1 \Proj_{\hat{U}_1} \times_2 \Proj_{\hat{U_2}} \times_3 \Proj_{\hat{U}_{3\perp}} \right\|_{\rm F}\\
	& + \left\|\Z\times_1 \hat{U}_1^\top \times_2 \hat{U}_2^\top \times_3 \hat{U}_3^\top \right\|_{\rm F}\\
	\leq & \left\|\X \times_1 \hat{U}_{1\perp}^\top \right\|_{\rm F} + \left\|\X \times_2 \hat{U}_{2\perp}^\top \right\|_{\rm F} + \left\|\X \times_3 \hat{U}_{3\perp}^\top \right\|_{\rm F}\\
	& + \left\|\Z\times_1 \hat{U}_1^\top \times_2 \hat{U}_2^\top \times_3 \hat{U}_3^\top \right\|_{\rm F}.
	\end{split}
	\end{equation}
	To obtain the upper bound of $\|\hat{\X} - \X\|$, we only need to analyze the four terms in \eqref{ineq:hat_X-X-decomposition} separately. Recall in Step 2, we defined
	\begin{equation*}
	\begin{split}
	& \mathcal{M}_1\left(\Y \times_2 (\hat{U}_2^{(t_{\max}-1)})^\top \times_3 (\hat{U}_3^{(t_{\max}-1)})^\top \right)\\
	 = & Y_1 \cdot \left(\hat{U}_2^{(t_{\max}-1)}\otimes \hat{U}_3^{(t_{\max}-1)}\right) := Y_1^{(t_{\max}-1)},\\
	\end{split}
	\end{equation*}
	$X_1^{(t_{\max}-1)}$, $Z_1^{(t_{\max}-1)}$ are defined similarly. Based on the calculation in \eqref{ineq:sigma_r_X_1} and \eqref{ineq:sigma_r_Z_1}, we have 
		\begin{equation*}
		\begin{split}
		\sigma_{\min}(X_1^{(t_{\max}-1)}) \geq \sigma_r(X_1) \cdot \sigma_{\min}\left(U_2^\top \hat{U}_2^{(t_{\max}-1)}\right)\cdot \sigma_{\min}\left(U_3^\top \hat{U}_3^{(t_{\max}-1)}\right) \geq \frac{3}{4}\lambda,
		\end{split}
		\end{equation*}
	\begin{equation*}
	\begin{split}
	\|Z_1^{(t_{\max}-1)}\| 
	\leq C\sqrt{p_1} + C\sqrt{pr} L_{t_{\max}-1}\leq C\sqrt{p_1} + C\sqrt{pr} \cdot \frac{\sqrt{p_1}}{\lambda} \leq C\sqrt{p_1}.\\
	\end{split}
	\end{equation*}
	Since $\hat{U}_1$ is the leading $r$ left singular vectors of $Y_1^{(t_{\max}-1)}=X_1^{(t_{\max}-1)}+Z_1^{(t_{\max}-1)}$, Lemma \ref{lm:r-pinciple-compoenent} implies
	\begin{equation*}
	\begin{split}
	& \left\|\Proj_{\hat{U}_{1\perp}} \mathcal{M}_1\left(\X \times_2 (\hat{U}_2^{(t_{\max}-1)})^\top \times_3 (\hat{U}_3^{(t_{\max}-1)})^\top \right) \right\|_{\rm F} \\
	= & \left\|\Proj_{\hat{U}_{1\perp}} X_1^{(t_{\max}-1)} \right\|_{\rm F} \leq C \sqrt{p_1r_1}.
	\end{split}
	\end{equation*}
	As a result,
	\begin{equation}\label{ineq:X*P_U1}
	\begin{split}
	& \left\|\X \times_1 \Proj_{\hat{U}_{1\perp}} \right\|_{\rm F} = \left\|\Proj_{\hat{U}_{1\perp}} \cdot X_1 \cdot (\Proj_{U_2}\otimes \Proj_{U_3}) \right\|_{\rm F} = \left\|\Proj_{\hat{U}_{1\perp}} \cdot X_1 \cdot \left(U_2 \otimes U_3 \right) \right\|_{\rm F}\\
	\leq & \left\|\Proj_{\hat{U}_{1\perp}} X_1 \left(\hat{U}_2^{(t_{\max}-1)} \otimes \hat{U}_3^{(t_{\max}-1)}\right) \right\|_{F} \cdot \sigma_{\min}^{-1}(U_2^\top \hat{U}_2^{(t_{\max}-1)}) \cdot \sigma_{\min}^{-1}(U_3^\top \hat{U}_3^{(t_{\max}-1)})\\
	\leq & C\sqrt{p_1r_1} \frac{1}{\sqrt{1 - (1/2)^2}} \frac{1}{\sqrt{1 - (1/2)^2}} \leq C\sqrt{p_1r_1}.
	\end{split}
	\end{equation}
	Similarly, we can show 
	\begin{equation}\label{ineq:X*P_U1'}
	\left\|\hat{\X} \times_2 \Proj_{\hat{U}_{2\perp}}\right\|_{\rm F} \leq C\sqrt{p_2r_2},\quad \left\|\hat{\X} \times_3 \Proj_{\hat{U}_{3\perp}}\right\|_{\rm F} \leq C\sqrt{p_3r_3}.
	\end{equation}
	Now combining \eqref{ineq:hat_X-X-decomposition}, \eqref{ineq:th1-check-3}, \eqref{ineq:X*P_U1}, and \eqref{ineq:X*P_U1'}, we have
	\begin{equation}
	\left\|\hat{\X} - \X\right\|_{\rm F} \leq C\left(\sqrt{p_1r_1} + \sqrt{p_2r_2} + \sqrt{p_3r_3}\right)
	\end{equation}
	for some constant $C>0$.
	
	\item We finalize the proof for Theorem \ref{th:upper_bound_strong} in this step. By Lemma \ref{lm:random_tensor_max_projection}, we know  \eqref{ineq:th1-check-1}, \eqref{ineq:th1-check-2}, and \eqref{ineq:th1-check-3} hold with probability at least $1 - C\exp(-cp)$. By the result in Step 1, we know \eqref{ineq:initial_performance} holds with probability at least $1 - C\exp(-cp)$. Let $Q = \{\text{\eqref{ineq:th1-check-1}, \eqref{ineq:th1-check-2}, \eqref{ineq:th1-check-3}, \eqref{ineq:initial_performance} all hold}\}$, then 
	\begin{equation}\label{ineq:PQ}
	P(Q) \geq 1-C\exp(-cp).
	\end{equation} 
	By Steps 2 and 3, one has $\|\sin\Theta(\hat{U}_k, U_k)\|\leq C\sqrt{p_k}/\lambda, k=1,2,3$, and 
	$$\|\hat{\X}-\X\|_{\rm F}\leq C\left(\sqrt{p_1r_1}+\sqrt{p_2r_2}+\sqrt{p_3r_3}\right) \text{ under $Q$}.$$ 
	It remains to consider situation under $Q^c$. By definition, $\hat{\X}$ is a projection of $\Y$, so
	\begin{equation*}
	\|\hat{\X}\|_{\rm F}\leq \|\Y\|_{\rm F}\leq \|\X\|_{\rm F}+\|\Z\|_{\rm F}.
	\end{equation*}
	Then we have the following rough upper bound for the 4-th moment of recovery error,
	\begin{equation*}
	\begin{split}
	& \mathbb{E}\|\hat{\X} - \X\|_{\rm F}^4 \leq C\left(\mathbb{E}\|\hat{\X}\|_{\rm F}^4+\|\X\|_{\rm F}^4\right)\leq C\|\X\|_{\rm F}^4 + C\mathbb{E}\|\Z\|_{\rm F}^4\\
	\leq & C\exp(c_0p) + C\mathbb{E}\left(\chi_{p_1p_2p_3}^2\right)^2 \leq C\exp(c_0p) + Cp^6.
	\end{split}
	\end{equation*}
	Then the following upper bound holds for the Frobenius norm risk of $\hat{\X}$,
	\begin{equation*}
	\begin{split}
	& \mathbb{E}\|\hat{\X}-\X\|_{\rm F}^2 = \mathbb{E}\|\hat{\X}-\X\|_{\rm F}^21_Q + \mathbb{E}\|\hat{\X}-\X\|_{\rm F}^21_{Q^c}\\
	\leq & C\left(\sqrt{p_1r_1}+\sqrt{p_2r_2}+\sqrt{p_3r_3}\right)^2 + \sqrt{\mathbb{E}\|\hat{\X}-\X\|_{\rm F}^4\cdot \mathbb{E}_{Q^c}}\\
	\overset{\eqref{ineq:PQ}}{\leq} & C\left(\sqrt{p_1r_1}+\sqrt{p_2r_2}+\sqrt{p_3r_3}\right)^2 + C\exp\left((c_0-c)p/2\right) + Cp^3\exp(-cp/2).
	\end{split}
	\end{equation*}
	Thus, one can select $c_0<c$ to ensure that
	$$\mathbb{E}\|\hat{\X}-\X\|_{\rm F}^2 \leq C\left(p_1r_1+p_2r_2+p_3r_3\right). $$
	
	Additionally, since $\sigma_{r_k}(\mathcal{M}_k(\X))\geq \lambda$, we have $\|\X\|_{\rm F^2} = \|\mathcal{M}_k(\X)\|_{\rm F}^2 \geq r_k\lambda^2$ for $k=1,2,3$, which implies $\|\X\|_{\rm F}^2\geq \max\{r_1, r_2, r_3\}\lambda = \lambda r$, then
	$$\mathbb{E}\frac{\|\hat{\X}-\X\|_{\rm F}^2}{\|\X\|_{\rm F}^2} \leq C\frac{p_1+p_2+p_3}{\lambda}. $$

	Moreover, by definition, $\|\sin\Theta(\hat{U}_k, U_k)\|\leq 1$. Thus we have the following upper bound for the spectral $\sin\Theta$ risk for $\hat{U}_k$,
	\begin{equation*}
	\begin{split}
	& \mathbb{E}\|\sin\Theta(\hat{U}_k, U_k)\|\leq \mathbb{E}\|\sin\Theta(\hat{U}_k, U_k)\|1_Q + \mathbb{E}\|\sin\Theta(\hat{U}_k, U_k)\|1_{Q^c}\\
	= & C\frac{\sqrt{p_k}}{\lambda} + \sqrt{\mathbb{E}\|\sin\Theta(\hat{U}_k, U_k)\|^2\cdot \mathbb{E} 1_{Q^c}}\overset{\eqref{ineq:PQ}}{\leq} C\frac{\sqrt{p_k}}{\lambda} + \sqrt{C\exp(-cp)}.
	\end{split}
	\end{equation*}	
	By definition of $\lambda$, we know $\lambda = \sigma_{r_k}(\mathcal{M}_k(\X)) \leq \frac{\|\X\|_{\rm F}}{\sqrt{r_k}}\leq \frac{C\exp(c_0p)}{\sqrt{r_k}}$, so one can select small constant $c_0>0$ to ensure that
	$$\frac{\sqrt{p_k}}{\lambda}\geq \frac{\sqrt{p_kr_k}}{C\exp(c_0p)} \geq c\sqrt{\exp(-cp)}, $$
	which implies $\mathbb{E}\|\sin\Theta(\hat{U}_k, U_k)\| \leq C\frac{p_k}{\lambda}$. Finally, we can derive the general Schatten $q$-$\sin\Theta$-norm risk via H\"older's inequality,
	$$\mathbb{E} r_k^{-1/q} \mathbb{E}\|\sin\Theta(\hat{U}_k, U_k)\|_{q} \leq \mathbb{E}\|\sin\Theta(\hat{U}_k, U_k)\| \leq C\frac{\sqrt{p_k}}{\lambda}.$$
\end{enumerate}

Summarizing from Steps 1-4, we have finished the proof of Theorem \ref{th:upper_bound_strong}. \quad $\square$

\subsection{Proof of Theorem~\ref{th:moderatelowerbound}}\label{sec:proof_moderatelowerbound}
We particularly show that it suffices to consider sparse tensor models and we set $\sigma=1$ for brevity. 
A tensor ${\bf X}\in\mathbb{R}^{p_1\times p_2\times p_3}$ is sparse with respect to parameters $\mathcal{S}(\X)=(s_1,s_2,s_3)$ if there exists $S_k({\bf X})\subset [p_k]:=\{1,2,\ldots,p_k\}, k=1,2,3$ such that
	$$
	{X}_{ijk}=0, \quad\forall (i,j,k)\in [p_1]\times [p_2]\times [p_3]\setminus S_1({\bf X})\times S_2({\bf X})\times S_3({\bf X})
	$$
	with $|S_k({\bf X})|\leq s_k, k=1,2,3$.
It means that the nonzero entries of ${\bf X}$ are constrained in the block $S_1({\bf X})\times S_2({\bf X})\times S_3({\bf X})$. 	Define the subset $\mathcal{M}(\bp,k,\br,\lambda)\subset \mathcal{F}_{\bp,\br}(\lambda)$ for integer $k=\floor{p^{(1-\tau)/2}}$ as follows,
	\begin{align*}
	\mathcal{M}(\bp,k,\br,\lambda):=\Big\{{\bf X}\in\mathcal{F}_{\bp,\br}(\lambda):  \mathcal{S}({\bf X})\leq (20k, 20k, 20k)\Big\},
	\end{align*}
	containing sparse tensors in $\mathcal{F}_{\bp,\br}(\lambda)$. Consider two disjoint subsets of $\mathcal{M}(\bp,k,\br,\lambda)$:
	$$
\mathcal{M}_0(\bp,k,\br,\lambda):=\Big\{{\bf X}\in\mathcal{M}(\bp,k,\br,\lambda), S_1({\bf X})\cup S_2({\bf X})\cup S_3({\bf X})\subset [p/2]\Big\},
$$
and
$$
\mathcal{M}_1(\bp,k,\br,\lambda):=\Big\{{\bf X}\in\mathcal{M}(\bp,k,\br,\lambda), S_1({\bf X})\cup S_2({\bf X})\cup S_3({\bf X})\subset [p]\setminus[p/2]\Big\}
$$
where matrices in $\mathcal{M}_0(\bp,k,\br,\lambda)$ and $\mathcal{M}_1(\bp,k,\br,\lambda)$ are supported on disjoint blocks, so are their singular vectors.
Given the observation:
$$
{\bf Y}={\bf X}+{\bf Z}\in\mathbb{R}^{p_1\times p_2\times p_3},
$$
the following testing problem is studied:
\begin{equation}\label{h01}
H_0: {\bf X}\in\mathcal{M}_0(\bp,k,\br,\lambda)\quad \textrm{V.S.}\quad H_1: {\bf X}\in\mathcal{M}_1(\bp,k,\br,\lambda).
\end{equation}
A test is then defined as $\phi(\cdot): \mathbb{R}^{p_1\times p_2\times p_3}\to \{0,1\}$ whose risk is given as
$$
\mathcal{R}_{\bp,\br,\lambda}(\phi)=\underset{\X\in \mathcal{M}_{0}(\bp,k,\br,\lambda)}{\sup}\mathbb{P}_{\X}\big\{\phi({\bf Y})=1\big\}+\underset{{\bf X}\in\mathcal{M}_{1}(\bp,k,\br,\lambda)}{\sup}\mathbb{P}_{\bf X}\big\{\phi({\bf Y})=0\big\},
$$
the worst case of Type-\rom{1}$+$\rom{2} error.
\begin{Lemma}\label{lemma:detection}
	Suppose Hypothesis {\bf H($\tau$)} for some $\tau\in(0,1)$. 
	Let $\{\phi_p\}$ be any sequence of polynomial-time tests of (\ref{h01}). There exists an absolute constant $c_0>0$ such that if $\lambda\leq c_0\Big(\frac{p^{3(1-\tau)/4}}{\sqrt{\log p}}\Big)$, then as long as $\min\{r_1,r_2,r_3\}\geq 1$,
	$$
	\underset{p\to\infty}{\lim\inf}\ \mathcal{R}_{\bp,\br,\lambda}(\phi_p)\geq \frac{1}{2}.
	$$
\end{Lemma}
Now we move back to the proof for Theorem \ref{th:moderatelowerbound}. Suppose that, on the contradiction, for any $k=1,2,3$, there exists a sub-sequence $(\hat U_k^{(p)})$ such that 
$$
\lim_{p\to\infty}\underset{\X\in \mathcal{F}_{\bp,\br}(\lambda)}{\sup} \mathbb{E}\big\|\sin\Theta\big(\hat U_k^{(p)},U_k(\X)\big)\big\|=0,
$$
which implies that
\begin{equation}\label{eq:UU-hatUhatU}
\lim_{p\to\infty}\underset{\X\in \mathcal{F}_{\bp,\br}(\lambda)}{\sup} \mathbb{P}\Big(\big\|U_kU_k^{\top}-\hat U_{k}^{(p)}\big(\hat U_{k}^{(p)}\big)^{\top}\big\|\leq\frac{1}{3}\Big)=1.
\end{equation}
Define a sequence of tests $\phi_{p}:\mathbb{R}^{p_1\times p_2\times p_3}\mapsto \{0,1\}$ as
$$
\phi_{p}(\Y):=
\begin{cases}
0,& {\rm if}\ \Big\|\big(\hat U_{k}^{(p)}\big)_{[1:p/2,:]}\big(\hat U_{k}^{(p)}\big)_{[1:p/2,:]}^{\top}\Big\|\geq \frac{2}{3},\\
1,&{\rm otherwise},
\end{cases}
$$
where $\big(\hat U_{k}^{(p)}\big)_{[1:p/2,:]}$ denote the first $p/2$ rows of $\hat U_{k}^{(p)}$. Clearly,
\begin{eqnarray*}
\mathcal{R}_{\bp,\br,\lambda}(\phi_{p})\leq \underset{\X\in\mathcal{M}_0(\bp,k,\br,\lambda)}{\sup}\mathbb{P}_{\X}\Big(\big\|U_kU_k^{\top}-\hat U_{k}^{(p)}\big(\hat U_{k}^{(p)}\big)^{\top}\big\|>\frac{1}{3}\Big)\\
+\underset{\X\in\mathcal{M}_1(\bp,k,\br,\lambda)}{\sup}\mathbb{P}_{\X}\Big(\big\|U_kU_k^{\top}-\hat U_{k}^{(p)}\big(\hat U_{k}^{(p)}\big)^{\top}\big\|\geq \frac{2}{3}\Big),
\end{eqnarray*}
which implies $\lim_{p\to\infty}\mathcal{R}_{\bp,\br,\lambda}\big(\phi_{p}\big)=0$, contradicting Lemma~\ref{lemma:detection}. Now, we prove claim (\ref{eq:moderate3}). Suppose that, on the contradiction, there exists a sub-sequence $(\hat\X^{(p)})$ such that
$$
\lim_{p\to\infty} \underset{\X\in\mathcal{F}_{\bp,\br}(\lambda)}{\sup}\mathbb{E}\frac{\|\hat\X^{(p)}-\X\|_{\rm F}^2}{\|\X\|_{\rm F}^2}=0,
$$
which implies
\begin{equation}\label{eq:hatX-X}
\lim_{p\to\infty}\underset{\X\in\mathcal{F}_{\bp,\br}(\lambda)}{\sup}\mathbb{P}\Big(\|\hat\X^{(p)}-\X\|_{\rm F}\leq \frac{1}{3}\|\X\|_{\rm F}\Big)=1.
\end{equation}
Define a sequence of test $\phi_p:\mathbb{R}^{p_1\times p_2\times p_3}\mapsto \{0,1\}$ as
$$
\phi_p(\Y):=
\begin{cases}
0,& {\rm if}\ \big\|\big(\hat\X^{(p)}\big)_{[V_1,V_1,V_1]}\big\|_{\rm F}\geq \big\|\big(\hat\X^{(p)}\big)_{[V_2,V_2,V_2]}\big\|_{\rm F},\\
1,& {\rm otherwise},
\end{cases}
$$
where $V_1=[p/2], V_2=[p]\setminus V_1$ and $\X_{[V_1,V_1,V_1]}$ denotes the sub-tensor on the block $V_1\times V_1\times V_1$. Under $H_0$, if $\phi_p(\Y)=1$, then
\begin{eqnarray*}
\|\hat\X^{(p)}-\X\|_{\rm F}^2\geq\big\|\big(\hat\X^{(p)}-\X\big)_{[V_1,V_1,V_1]}\big\|_{\rm F}^2+\big\|\big(\hat\X^{(p)}\big)_{[V_2,V_2,V_2]}\big\|_{\rm F}^2\\
\geq\big\|\big(\hat\X^{(p)}-\X\big)_{[V_1,V_1,V_1]}\big\|_{\rm F}^2+\big\|\big(\hat\X^{(p)}\big)_{[V_1,V_1,V_1]}\big\|_{\rm F}^2\geq \frac{1}{2}\|\X\|_{\rm F}^2.
\end{eqnarray*}
Clearly,
\begin{eqnarray*}
\mathcal{R}_{\bp,\br,\lambda}(\phi_p)=\underset{\X\in\mathcal{M}_0(\bp,k,\br,\lambda)}{\sup} \mathbb{P}\Big(\|\hat\X^{(p)}-\X\|_{\rm F}\geq \frac{\sqrt{2}}{2}\|\X\|_{\rm F}\Big)\\
+\underset{\X\in\mathcal{M}_1(\bp,k,\br,\lambda)}{\sup}\mathbb{P}\Big(\|\hat\X^{(p)}-\X\|_{\rm F}\geq \frac{\sqrt{2}}{2}\|\X\|_{\rm F}\Big),
\end{eqnarray*}
which implies that $\mathcal{R}_{\bp,\br,\lambda}(\phi_p)\to 0$ as $p\to\infty$ based on (\ref{eq:hatX-X}), which contradicts Lemma~\ref{lemma:detection}. \quad $\square$

\subsection{Proof of Theorems \ref{th:lower_bound}}\label{sec:proof_lower_bound} 
Without loss of generality we can assume $\sigma=1$ throughout the proof. First, we construct the core tensor $\tilde{\S} \in \mathbb{R}^{r_1\times r_2\times r_3}$ with i.i.d. standard Gaussian entries, then according to random matrix theory (c.f. Corollary 5.35 in \cite{vershynin2010introduction}), with probability at least $1 - 6e^{-x}$, we have
$$\sqrt{r_{k+1}r_{k+2}} - \sqrt{r_k} - x \leq \sigma_{\min}(\mathcal{M}_k(\tilde{\S})) \leq \sigma_{\max}(\mathcal{M}_k(\tilde{\S})) \leq \sqrt{r_{k+1}r_{k+2}} + \sqrt{r_k} + x,$$
for $k = 1,2,3$. Plug in $x = 1.8$, by simple calculation, we can see there is a positive probability that
\begin{equation}
\sqrt{r_{k+1}r_{k+2}} - \sqrt{r_k} - 1.8 \leq \sigma_{\min}(\mathcal{M}_k(\tilde{\S})) \leq \sigma_{\max}(\mathcal{M}_k(\tilde{\S})) \leq \sqrt{r_{k+1}r_{k+2}} + \sqrt{r_k} + 1.8.
\end{equation}
Note that
\begin{equation*}
\begin{split}
r_1r_2 \geq 4 r_3, r_2r_3\geq 4 r_1, r_3r_1 \geq 4 r_2, \quad \Rightarrow & \quad r_1r_2r_3 \geq 4\max_{1\leq k \leq3}\{r_k\}^2\\
\Rightarrow & \quad r_k\frac{r_{k+1}}{\max\{r_k\}}\frac{r_{k+2}}{\max\{r_k\}} \geq 4, \Rightarrow r_k \geq 4\\
\Rightarrow & \quad r_{k+1}r_{k+2} \geq \frac{4\max\{r_k\}^2}{r_k} \geq 4 r_k \geq 16,
\end{split}
\end{equation*}
we know 
\begin{equation*}
\begin{split}
& \sqrt{r_{k+1}r_{k+2}}-\sqrt{r_k} - 1.8 \\
\geq & \sqrt{r_{k+1}r_{k+2}} - \sqrt{\frac{r_{k+1}r_{k_2}}{4}} - \frac{1.8}{4}\sqrt{r_{k+1}r_{k+2}} = 0.05\sqrt{r_{k+1}r_{k+2}};
\end{split}
\end{equation*}
\begin{equation*}
\begin{split}
& \sqrt{r_{k+1}r_{k+2}}-\sqrt{r_k} + 1.8 \\
\leq & \sqrt{r_{k+1}r_{k+2}} + \sqrt{\frac{r_{k+1}r_{k_2}}{4}} + \frac{1.8}{4}\sqrt{r_{k+1}r_{k+2}} = 1.95\sqrt{r_{k+1}r_{k+2}}.
\end{split}
\end{equation*}
By previous arguments, there exists $\tilde{\S}\in \mathbb{R}^{r_1\times r_2\times r_3}$ such that $c\sqrt{r_{k+1}r_{k+2}}\leq \sigma_{\min}(\mathcal{M}_k(\tilde{\S}))\leq C\sqrt{r_{k+1}r_{k+2}}$ for $k=1,2,3$. Now, we construct the scaled core tensor $\S = \tilde{\S} \frac{\lambda}{\min_{k=1,2,3}\sigma_{\min}(\mathcal{M}_k(\tilde{\S}))}$. Given $r\leq r_1, r_2, r_3 \leq C_0 r$, we know $\S \in \mathbb{R}^{r_1\times r_2 \times r_3}$ satisfies the following property
\begin{equation}
\lambda \leq \sigma_{\min}\left(\mathcal{M}_k(\S)\right) \leq \sigma_{\max}\left(\mathcal{M}_k(\S)\right)\leq C\lambda, \quad k=1,2,3.
\end{equation}

\paragraph{Proof of the first claim.}
It suffices to consider $k=1$.
We construct a large subset of $\mathbb{O}_{p_1,r_1}$ whose elements are well separated in Schatten $q$-norms for all $1\leq q\leq +\infty$. To this end, we need some preliminary facts about the packing number in Grassmann manifold $\mathcal{G}_{p,r}$, which is the set of all $r$-dimensional subspaces of $\mathbb{R}^{p}$. Given such a subspace $L\subset \mathbb{R}^p$ with ${\rm dim}(L)=r$, let $U_L\in\mathbb{O}(p,r)$ denote the orthonormal basis of $L$. Denote $\mathcal{B}_{p,r}:=\{U_L, L\in\mathcal{G}_{p,r}\}$ which is actually a subset of $\mathbb{O}_{p,r}$ and will be equipped with Schatten $q$-norm distances for all $q\in[1,+\infty]$: $d_q(U_{L_1}, U_{L_2}):=\|U_{L_1}U_{L_1}^{\top}-U_{L_2}U_{L_2}^{\top}\|_q$. Recall that the $\varepsilon$-packing number of a metric space $(T,d)$ is defined as
$$
D(T,d,\varepsilon):=\max\Big\{n: {\rm there\ are}\ t_1,\ldots,t_n\in T,\ \rm{such\ that }\ \min_{i\neq j}\ d(t_i,t_j)>\varepsilon\Big\}.
$$ 
The following lemma can be found in Lemma 5 in~\cite{koltchinskii2015optimal} which controls the packing numbers of $\mathcal{B}_{p,r}$ with respect to Schatten distances $d_q$.
\begin{Lemma}\label{lemma:Grassmann}
	For all integers $1\leq r\leq p$ such that $r\leq p-r$, and all $1\leq q\leq +\infty$, the following bound holds
	$$
	\Big(\frac{c}{\varepsilon}\Big)^{r(p-r)}\leq D\big(\mathcal{B}_{p,r},d_q,\varepsilon r^{1/q}\big)\leq \Big(\frac{C}{\varepsilon}\Big)^{r(p-r)}
	$$
	with absolute constants $c,C>0$.
\end{Lemma}
We are in position to construct a well-separated subset of $\mathbb{O}_{p_1,r_1}$.
According to Lemma~\ref{lemma:Grassmann} by choosing $\varepsilon=\frac{c}{2}$, there exists a subset $\mathcal{V}_{p_1-r_1,r_1}\subset \mathbb{O}_{p_1-r_1,r_1}$ with ${\rm Card}\big(\mathcal{V}_{p_1-r_1,r_1}\big)\geq 2^{r_1(p_1-2r_1)}$ such that for each $V_1\neq V_2\in \mathcal{V}_{p_1-r_1,r_1}$,
$$
\|V_1V_1^{\top}-V_2V_2^{\top}\|_q\geq \frac{c}{2}r_1^{1/q}.
$$
Now, fix a $\delta>0$ whose value is to be determined later. For every $V\in \mathcal{V}_{p_1-r_1,r_1}$, define $\tilde{V}\in \mathbb{O}_{p_1,r_1}$ as follows
$$
\tilde{V}=\left(
\begin{array}{c}
\sqrt{1-\delta}I_{r_1}\\
\sqrt{\delta} V
\end{array}
\right).
$$
It is easy to check that $\tilde{V}\in \mathbb{O}_{p_1,r_1}$ as long as $V\in\mathbb{O}_{p_1-r_1,r_1}$. We conclude with a subset $\mathcal{V}_{p_1,r_1}\subset \mathbb{O}_{p_1,r_1}$ with ${\rm Card}(\mathcal{V}_{p_1,r_1})={\rm Card}(\mathcal{V}_{p_1-r_1,r_1})\geq 2^{r_1(p_1-2r_1)}$. Moreover, for $\tilde{V}_1\neq \tilde{V}_2\in \mathcal{V}_{p_1,r_1}$,
$$
\|\tilde{V}_1\tilde{V}_1^{\top}-\tilde{V}_2\tilde{V}_2^{\top}\|_q\geq \sqrt{\delta(1-\delta)}\|V_1-V_2\|_q\geq \frac{c}{2}\sqrt{\delta(1-\delta)}r_1^{1/q}.
$$
Meanwhile,
$$
\|\tilde{V}_1-\tilde{V}_2\|_{\rm F}\leq \sqrt{\delta}\|V_1-V_2\|_{\rm F}\leq \sqrt{2\delta r_1}.
$$
Then we construct a series of fixed signal tensors: $\X_i = \S\times_1 \tilde{V}_i \times_2 U_2\times_3 U_3, i=1,\ldots, m$, where $U_2\in \mathbb{O}_{p_2,r_2}. U_3\in \mathbb{O}_{p_3,r_3}$ are any fixed orthonormal columns, $\tilde{V}_i\in\mathcal{V}_{p_1,r_1}\subset \mathbb{O}_{p_1,r_1}$ and $m=2^{r_1(p_1-2r_1)}$. By such construction, $\sigma_{\min}^2(\mathcal{M}_k(\X_i)) \geq \lambda$ for $k=1,2,3$, so that $\{\X_i\}_{i=1}^{m} \subseteq \mathcal{F}_{\bp, \br}(\lambda)$.

We further let $\Y_i = \X_i + \Z_i$, where $\Z_i$ are i.i.d. standard normal distributed tensors, which implies $\Y_i \sim N(\X_i, I_{p_1\times p_2\times p_3})$. Then the Kullback-Leibler divergence between the distribution $\Y_i$ and $\Y_j$ is
\begin{equation*}
\begin{split}
& D_{\rm KL}\left(\Y_i|| \Y_j\right) = \frac{1}{2} \left\|\X_i - \X_j\right\|_{\rm F}^2 = \frac{1}{2}\left\|\S\times_1 (\tilde{V}_i - \tilde{V}_j) \times_2 U_2\times_3 U_3\right\|_{\rm F}^2\\
= & \frac{1}{2} \|(\tilde{V}_i-\tilde{V}_j)\cdot \mathcal{M}_1(\S) \cdot (U_2\otimes U_3)^\top \|_{\rm F}^2 \leq C\lambda^2\left\|\tilde{V}_i - \tilde{V}_j\right\|_{\rm F}^2\leq 2C\lambda^2\delta r_1.
\end{split}
\end{equation*}
Then the generalized Fano's lemma yields the following lower bound
\begin{equation}\label{ineq:lower_bound_1}
\inf_{\hat{U}_1} \sup_{U_1 \in \{V_i\}_{i=1}^m} \mathbb{E} \left\|\hat{U}_1\hat{U}_1^{\top}-U_1U_1^{\top}\right\|_{q} \geq \frac{c}{2}\sqrt{\delta(1-\delta)}r_1^{1/q}\left(1 - \frac{C\lambda^2\delta r_1 +\log 2}{r_1(p_1-2r_1)\log 2}\right).
\end{equation}
By setting $\delta=c_1\frac{(p_1-2r_1)}{\lambda^2}$ for a small but absolute constant $c_1>0$, we obtain 
$$
\frac{c}{2}\sqrt{\delta(1-\delta)}r_1^{1/q}\left(1 - \frac{C\lambda^2\delta r_1 +\log 2}{r_1(p_1-2r_1)\log 2}\right)\geq c_0\frac{(p_1-2r_1)^{1/2}}{\lambda}r_1^{1/q}
$$
for an absolute constant $c_0>0$. Then, if $p_1\geq 3r_1$, 
$$
\inf_{\hat{U}_1} \sup_{U_1 \in \{V_i\}_{i=1}^m} \mathbb{E} \left\|\hat{U}_1\hat{U}_1^{\top}-U_1U_1^{\top}\right\|_{q} \geq c_0\Big(\frac{\sqrt{p_1}}{\lambda}r_1^{1/q}\wedge r_1^{1/q}\Big)
$$
where $r_1^{1/q}$ is a trivial term. The first claim in Theorem~\ref{th:lower_bound} it thus obtained by viewing the equivalence between the Schatten $q$-norms and $\sin\Theta$ Schatten $q$-norms, see Lemma~\ref{lemma:norms} in the Appendix.

\paragraph{Proof of second and third claims.} To prove the minimax lower bounds in estimating $\X$, we need a different construction scheme. Specifically, we consider the metric space $\left(\mathbb{O}_{p_1, r_1}, \|\sin\Theta(\cdot, \cdot)\|_2\right)$, fix an $V_0 \in \mathbb{O}_{p_1, r_1}$, and consider the following ball of radius $\varepsilon>0$ and center $V_0$:
\begin{equation*}
B(V_0, \varepsilon) = \left\{V' \in \mathbb{O}_{p_1, r_1}: \|\sin\Theta(V', V)\|_2 \leq \varepsilon\right\}.
\end{equation*}
By Lemma 1 in \cite{cai2013sparse}, for $0 < \alpha<1$ and $0 <\varepsilon \leq 1$, there exists $V'_1, \ldots, V'_m \subseteq B(V_0, \varepsilon)$ such that
\begin{equation*}
m\geq \left(\frac{c_0}{\alpha}\right)^{r_1(p_1-r_1)}, \quad \min_{1\leq i< j\leq m} \left\|\sin\Theta(V'_i, V'_j)\right\|_2 \geq \alpha\varepsilon.
\end{equation*}
By the property of $\sin\Theta$ distance (Lemma 1 in \cite{cai2016rate}), we can find a rotation matrix $O_i\in \mathbb{O}_{r_1}$ such that
$$\|V_0 - V_i'O_i\|_{\rm F} \leq \sqrt{2}\|\sin\Theta(V_0, V_i')\|_2 \leq \sqrt{2}\varepsilon.$$
We denote $V_i = V_i'O_i$, then 
\begin{equation}
\|V_i  - V_0\|_{\rm F}  \leq \sqrt{2}\varepsilon,\quad \|\sin\Theta(V_i, V_j)\|_2 \geq \alpha \varepsilon, \quad 1\leq i<j\leq m.
\end{equation}
Construct $\X_i=\S\times_1 V_i\times _2 U_2\times U_3$ for $i=1,\ldots,m$ in a similar fashion as above.
Then the class of low-rank tensors satisfy the following properties,
\begin{equation}
\begin{split}
& \left\|\X_i - \X_j\right\|_{\rm F}^2 = \left\|\S\times_1 (V_i - V_j) \times_2 U_2\times_3 U_3\right\|_{\rm F}^2\\
= & \frac{1}{2} \|(V_i-V_j)\cdot \mathcal{M}_1(\S) \cdot (U_2\otimes U_3)^\top \|_{\rm F}^2 \geq \frac{1}{2} \sigma_{r_1}\left(\mathcal{M}_1(\S)\right) \left\|V_i - V_j\right\|_{\rm F}^2\\ 
\geq & \frac{\lambda^2}{2} \min_{O\in \mathbb{O}_r} \|V_i - V_jO\|_F^2  \quad \text{(by Lemma 1 in \cite{cai2016rate})}\\
\geq & \frac{\lambda^2}{2}\left\|\sin\Theta\left(V_i, V_j\right)\right\|_2^2 \geq \alpha^2 \varepsilon^2\lambda^2, \quad 1\leq i < j\leq m.
\end{split}
\end{equation}
\begin{equation}\label{ineq:lower_bound_X_upper_bound}
\begin{split}
\min_{1\leq i \leq m} \left\|\X_i\right\|_{\rm F} = \left\|\S\times_1 V_i \times_2 U_2\times_3 U_3\right\|_{\rm F} = \|\S\|_{\rm F} \geq \|\mathcal{M}_1(\S)\|_{\rm F} \geq \lambda\sqrt{r_1}.
\end{split}
\end{equation}
Moreover, under the same model $\Y_i=\X_i+\Z_i$ as above, the KL-divergence between the distributions of $\Y_i$ and $\Y_j$ is
\begin{equation*}
\begin{split}
& D_{\rm KL}\left(\Y_i|| \Y_j\right) = \frac{1}{2} \left\|\X_i - \X_j\right\|_{\rm F}^2 = \frac{1}{2}\left\|\S\times_1 (V_i - V_j) \times_2 U_2\times_3 U_3\right\|_{\rm F}^2\\
= & \frac{1}{2} \|(V_i-V_j)\cdot \mathcal{M}_1(\S) \cdot (U_2\otimes U_3)^\top \|_{\rm F}^2 \leq C\lambda^2\left\|V_i - V_j\right\|_{\rm F}^2\leq C\lambda^2\varepsilon^2.
\end{split}
\end{equation*}
Then the generalized Fano's lemma yields the following lower bound
\begin{equation}\label{ineq:lower_bound_2}
\inf_{\hat{\X}} \sup_{\X \in \{\X_i\}_{i=1}^m} \mathbb{E} \left\|\hat{\X} - \X\right\|_{\rm F}^2 \geq \lambda\alpha\varepsilon\left(1 - \frac{C\lambda^2\varepsilon^2 +\log 2}{r_1(p_1-r_1)\log (c_0/\alpha)}\right).
\end{equation}
By setting $\varepsilon = \sqrt{\frac{r_1(p_1-r_1)}{2C\lambda^2}} \wedge \sqrt{2r_1}$, $\alpha = (c_0 \wedge 1)/8$, we have
\begin{equation*}
\alpha\varepsilon\left(1 - \frac{C\lambda^2\varepsilon^2 +\log 2}{r_1(p_1-r_1)\log (c_0/\alpha)}\right) \geq c_1 \left(\frac{\sqrt{r_1p_1}}{\lambda} \wedge \sqrt{r_1}\right)
\end{equation*}
for some small constant $c_1 > 0$. 
Moreover,
\begin{equation}\label{ineq:lower_bound_3}
\begin{split}
& \inf_{\hat{\X}} \sup_{\X\in\mathcal{F}_{\bp, \br}(\lambda)}  \mathbb{E} \left\|\hat{\X} - \X\right\|_{\rm F}^2 \geq \inf_{\hat{\X}} \sup_{\X\in\mathcal{F}_{\bp, \br}(\lambda\vee \sqrt{p_1})}  \mathbb{E} \left\|\hat{\X} - \X\right\|_{\rm F}^2 \\
\geq & \inf_{\hat{\X}} \sup_{\X\in\{\X_1,\ldots, \X_m\}} \mathbb{E} \left\|\hat{\X} - \X\right\|_{\rm F}^2 \overset{\eqref{ineq:lower_bound_2}}{\geq}c_1(\lambda^2 \vee p_1)\left(\frac{r_1p_1}{\lambda^2 \vee p_1}\wedge r_1\right)\geq c_1p_1r_1.
\end{split}
\end{equation}
\begin{equation}\label{ineq:lower_bound_4}
\begin{split}
& \inf_{\hat{\X}} \sup_{\X\in\mathcal{F}_{\bp, \br}(\lambda)}  \mathbb{E} \frac{\|\hat{\X} - \X\|_{\rm F}^2}{\|\X\|_{\rm F}^2} \geq \inf_{\hat{\X}} \sup_{\X\in\mathcal{F}_{\bp, \br}(\lambda\vee p_1)}  \mathbb{E} \frac{\|\hat{\X} - \X\|_{\rm F}^2}{\|\X\|_{\rm F}^2}\\
\geq & \inf_{\hat{\X}} \sup_{\X\in\{\X_1,\ldots, \X_m\}} \mathbb{E} \frac{\|\hat{\X} - \X\|_{\rm F}^2}{\max_{1\leq i \leq m}\|\X\|_{\rm F}^2} \overset{\eqref{ineq:lower_bound_X_upper_bound}\eqref{ineq:lower_bound_2}}{\geq}\frac{c_1\lambda^2}{\lambda^2r_1}\left(\frac{r_1p_1}{\lambda^2}\wedge r_1\right)\geq c_1\left(\frac{p_1}{\lambda^2}\wedge 1\right).
\end{split}
\end{equation}

Finally, we apply the same argument of \eqref{ineq:lower_bound_2}, \eqref{ineq:lower_bound_3}, and \eqref{ineq:lower_bound_4} to $U_2, U_3$, then we can obtain \eqref{eq:stronglowerbound1} and \eqref{eq:stronglowerbound2}.
\quad $\square$

\section*{Acknowledgment}

The authors would like to thank Ming Yuan and Zongming Ma for helpful discussions. The authors would also like to thank Genevera Allen and Nathaniel Helwig for pointing out several key references. The authors would also like to thank editors and anomalous referees for your helpful comments and suggestions on improving the paper.


\bibliographystyle{ieeetr}
\bibliography{reference}

\begin{thebibliography}{10}

\bibitem{shen2008sparse}
H.~Shen and J.~Z. Huang, ``Sparse principal component analysis via regularized
  low rank matrix approximation,'' {\em Journal of multivariate analysis},
  vol.~99, no.~6, pp.~1015--1034, 2008.

\bibitem{lee2010biclustering}
M.~Lee, H.~Shen, J.~Z. Huang, and J.~Marron, ``Biclustering via sparse singular
  value decomposition,'' {\em Biometrics}, vol.~66, no.~4, pp.~1087--1095,
  2010.

\bibitem{yang2014sparse}
D.~Yang, Z.~Ma, and A.~Buja, ``A sparse singular value decomposition method for
  high-dimensional data,'' {\em Journal of Computational and Graphical
  Statistics}, vol.~23, no.~4, pp.~923--942, 2014.

\bibitem{yang2016rate}
D.~Yang, Z.~Ma, and A.~Buja, ``Rate optimal denoising of simultaneously sparse
  and low rank matrices,'' {\em Journal of Machine Learning Research}, vol.~17,
  no.~92, pp.~1--27, 2016.

\bibitem{candes2013unbiased}
E.~J. Candes, C.~A. Sing-Long, and J.~D. Trzasko, ``Unbiased risk estimates for
  singular value thresholding and spectral estimators,'' {\em IEEE transactions
  on signal processing}, vol.~61, no.~19, pp.~4643--4657, 2013.

\bibitem{shabalin2013reconstruction}
A.~A. Shabalin and A.~B. Nobel, ``Reconstruction of a low-rank matrix in the
  presence of gaussian noise,'' {\em Journal of Multivariate Analysis},
  vol.~118, pp.~67--76, 2013.

\bibitem{donoho2014minimax}
D.~Donoho, M.~Gavish, {\em et~al.}, ``Minimax risk of matrix denoising by
  singular value thresholding,'' {\em The Annals of Statistics}, vol.~42,
  no.~6, pp.~2413--2440, 2014.

\bibitem{gavish2014optimal}
M.~Gavish and D.~L. Donoho, ``The optimal hard threshold for singular values is
  $4/\sqrt{3}$,'' {\em IEEE Transactions on Information Theory}, vol.~60,
  no.~8, pp.~5040--5053, 2014.

\bibitem{zou2006sparse}
H.~Zou, T.~Hastie, and R.~Tibshirani, ``Sparse principal component analysis,''
  {\em Journal of computational and graphical statistics}, vol.~15, no.~2,
  pp.~265--286, 2006.

\bibitem{cai2013sparse}
T.~T. Cai, Z.~Ma, Y.~Wu, {\em et~al.}, ``Sparse pca: Optimal rates and adaptive
  estimation,'' {\em The Annals of Statistics}, vol.~41, no.~6, pp.~3074--3110,
  2013.

\bibitem{birnbaum2013minimax}
A.~Birnbaum, I.~M. Johnstone, B.~Nadler, and D.~Paul, ``Minimax bounds for
  sparse pca with noisy high-dimensional data,'' {\em Annals of statistics},
  vol.~41, no.~3, p.~1055, 2013.

\bibitem{candes2011robust}
E.~J. Cand{\`e}s, X.~Li, Y.~Ma, and J.~Wright, ``Robust principal component
  analysis?,'' {\em Journal of the ACM (JACM)}, vol.~58, no.~3, p.~11, 2011.

\bibitem{zhou2013tensor}
H.~Zhou, L.~Li, and H.~Zhu, ``Tensor regression with applications in
  neuroimaging data analysis,'' {\em Journal of the American Statistical
  Association}, vol.~108, no.~502, pp.~540--552, 2013.

\bibitem{zhang2016cross}
A.~Zhang, ``Cross: Efficient low-rank tensor completion,'' {\em arXiv preprint
  arXiv:1611.01129}, 2016.

\bibitem{karatzoglou2010multiverse}
A.~Karatzoglou, X.~Amatriain, L.~Baltrunas, and N.~Oliver, ``Multiverse
  recommendation: n-dimensional tensor factorization for context-aware
  collaborative filtering,'' in {\em Proceedings of the fourth ACM conference
  on Recommender systems}, pp.~79--86, ACM, 2010.

\bibitem{rendle2010pairwise}
S.~Rendle and L.~Schmidt-Thieme, ``Pairwise interaction tensor factorization
  for personalized tag recommendation,'' in {\em Proceedings of the third ACM
  international conference on Web search and data mining}, pp.~81--90, ACM,
  2010.

\bibitem{lu2008mpca}
H.~Lu, K.~N. Plataniotis, and A.~N. Venetsanopoulos, ``Mpca: Multilinear
  principal component analysis of tensor objects,'' {\em IEEE Transactions on
  Neural Networks}, vol.~19, no.~1, pp.~18--39, 2008.

\bibitem{liu2013tensor}
J.~Liu, P.~Musialski, P.~Wonka, and J.~Ye, ``Tensor completion for estimating
  missing values in visual data,'' {\em IEEE Transactions on Pattern Analysis
  and Machine Intelligence}, vol.~35, no.~1, pp.~208--220, 2013.

\bibitem{anandkumar2014tensor-latent}
A.~Anandkumar, R.~Ge, D.~J. Hsu, S.~M. Kakade, and M.~Telgarsky, ``Tensor
  decompositions for learning latent variable models.,'' {\em Journal of
  Machine Learning Research}, vol.~15, no.~1, pp.~2773--2832, 2014.

\bibitem{anandkumar2014tensor-membership}
A.~Anandkumar, R.~Ge, D.~Hsu, and S.~M. Kakade, ``A tensor approach to learning
  mixed membership community models,'' {\em J. Mach. Learn. Res.}, vol.~15,
  pp.~2239--2312, Jan. 2014.

\bibitem{li2010tensor}
N.~Li and B.~Li, ``Tensor completion for on-board compression of hyperspectral
  images,'' in {\em 2010 IEEE International Conference on Image Processing},
  pp.~517--520, IEEE, 2010.

\bibitem{liu2017characterizing}
T.~Liu, M.~Yuan, and H.~Zhao, ``Characterizing spatiotemporal transcriptome of
  human brain via low rank tensor decomposition,'' {\em arXiv preprint
  arXiv:1702.07449}, 2017.

\bibitem{richard2014statistical}
E.~Richard and A.~Montanari, ``A statistical model for tensor pca,'' in {\em
  Advances in Neural Information Processing Systems}, pp.~2897--2905, 2014.

\bibitem{hopkins2015tensor}
S.~B. Hopkins, J.~Shi, and D.~Steurer, ``Tensor principal component analysis
  via sum-of-square proofs.,'' in {\em COLT}, pp.~956--1006, 2015.

\bibitem{perry2016statistical}
A.~Perry, A.~S. Wein, and A.~S. Bandeira, ``Statistical limits of spiked tensor
  models,'' {\em arXiv preprint arXiv:1612.07728}, 2016.

\bibitem{cai2016rate}
T.~T. Cai and A.~Zhang, ``Rate-optimal perturbation bounds for singular
  subspaces with applications to high-dimensional statistics,'' {\em The Annals
  of Statistics}, vol.~to appear, 2017.

\bibitem{hillar2013most}
C.~J. Hillar and L.-H. Lim, ``Most tensor problems are np-hard,'' {\em Journal
  of the ACM (JACM)}, vol.~60, no.~6, p.~45, 2013.

\bibitem{de2000multilinear}
L.~De~Lathauwer, B.~De~Moor, and J.~Vandewalle, ``A multilinear singular value
  decomposition,'' {\em SIAM journal on Matrix Analysis and Applications},
  vol.~21, no.~4, pp.~1253--1278, 2000.

\bibitem{de2000best}
L.~De~Lathauwer, B.~De~Moor, and J.~Vandewalle, ``On the best rank-1 and
  rank-(r 1, r 2,..., rn) approximation of higher-order tensors,'' {\em SIAM
  Journal on Matrix Analysis and Applications}, vol.~21, no.~4, pp.~1324--1342,
  2000.

\bibitem{vasilescu2003multilinear}
M.~A.~O. Vasilescu and D.~Terzopoulos, ``Multilinear subspace analysis of image
  ensembles,'' in {\em Computer Vision and Pattern Recognition, 2003.
  Proceedings. 2003 IEEE Computer Society Conference on}, vol.~2, pp.~II--93,
  IEEE, 2003.

\bibitem{sheehan2007higher}
B.~N. Sheehan and Y.~Saad, ``Higher order orthogonal iteration of tensors
  (hooi) and its relation to pca and glram,'' in {\em Proceedings of the 2007
  SIAM International Conference on Data Mining}, pp.~355--365, SIAM, 2007.

\bibitem{costantini2008higher}
R.~Costantini, L.~Sbaiz, and S.~Susstrunk, ``Higher order svd analysis for
  dynamic texture synthesis,'' {\em IEEE Transactions on Image Processing},
  vol.~17, no.~1, pp.~42--52, 2008.

\bibitem{haardt2008higher}
M.~Haardt, F.~Roemer, and G.~Del~Galdo, ``Higher-order svd-based subspace
  estimation to improve the parameter estimation accuracy in multidimensional
  harmonic retrieval problems,'' {\em IEEE Transactions on Signal Processing},
  vol.~56, no.~7, pp.~3198--3213, 2008.

\bibitem{liu2014generalized}
Y.~Liu, F.~Shang, W.~Fan, J.~Cheng, and H.~Cheng, ``Generalized higher-order
  orthogonal iteration for tensor decomposition and completion,'' in {\em
  Advances in Neural Information Processing Systems}, pp.~1763--1771, 2014.

\bibitem{zheng2015interpolating}
Q.~Zheng and R.~Tomioka, ``Interpolating convex and non-convex tensor
  decompositions via the subspace norm,'' in {\em Advances in Neural
  Information Processing Systems}, pp.~3106--3113, 2015.

\bibitem{anandkumar2016homotopy}
A.~Anandkumar, Y.~Deng, R.~Ge, and H.~Mobahi, ``Homotopy analysis for tensor
  pca,'' {\em arXiv preprint arXiv:1610.09322}, 2016.

\bibitem{anandkumar2014guaranteed}
A.~Anandkumar, R.~Ge, and M.~Janzamin, ``Guaranteed non-orthogonal tensor
  decomposition via alternating rank-$1 $ updates,'' {\em arXiv preprint
  arXiv:1402.5180}, 2014.

\bibitem{sun2015provable}
W.~W. Sun, J.~Lu, H.~Liu, and G.~Cheng, ``Provable sparse tensor
  decomposition,'' {\em Journal of Royal Statistical Association}, 2015.

\bibitem{allen2012regularized}
G.~I. Allen, ``Regularized tensor factorizations and higher-order principal
  components analysis,'' {\em arXiv preprint arXiv:1202.2476}, 2012.

\bibitem{allen2012sparse}
G.~Allen, ``Sparse higher-order principal components analysis.,'' in {\em
  AISTATS}, vol.~15, 2012.

\bibitem{qi2016uniqueness}
Y.~Qi, P.~Comon, and L.-H. Lim, ``Uniqueness of nonnegative tensor
  approximations,'' {\em IEEE Transactions on Information Theory}, vol.~62,
  no.~4, pp.~2170--2183, 2016.

\bibitem{lu2016tensor}
C.~Lu, J.~Feng, Y.~Chen, W.~Liu, Z.~Lin, and S.~Yan, ``Tensor robust principal
  component analysis: Exact recovery of corrupted low-rank tensors via convex
  optimization,'' in {\em Proceedings of the IEEE Conference on Computer Vision
  and Pattern Recognition}, pp.~5249--5257, 2016.

\bibitem{lesieur2017statistical}
T.~Lesieur, L.~Miolane, M.~Lelarge, F.~Krzakala, and L.~Zdeborov{\'a},
  ``Statistical and computational phase transitions in spiked tensor
  estimation,'' in {\em Information Theory (ISIT), 2017 IEEE International
  Symposium on}, pp.~511--515, IEEE, 2017.

\bibitem{kolda2009tensor}
T.~G. Kolda and B.~W. Bader, ``Tensor decompositions and applications,'' {\em
  SIAM review}, vol.~51, no.~3, pp.~455--500, 2009.

\bibitem{berthet2013computational}
Q.~Berthet and P.~Rigollet, ``Computational lower bounds for sparse pca,'' {\em
  arXiv preprint arXiv:1304.0828}, 2013.

\bibitem{berthet2013optimal}
Q.~Berthet, P.~Rigollet, {\em et~al.}, ``Optimal detection of sparse principal
  components in high dimension,'' {\em The Annals of Statistics}, vol.~41,
  no.~4, pp.~1780--1815, 2013.

\bibitem{wang2014statistical}
T.~Wang, Q.~Berthet, and R.~J. Samworth, ``Statistical and computational
  trade-offs in estimation of sparse principal components,'' {\em arXiv
  preprint arXiv:1408.5369}, 2014.

\bibitem{gao2014sparse}
C.~Gao, Z.~Ma, and H.~H. Zhou, ``Sparse cca: Adaptive estimation and
  computational barriers,'' {\em arXiv preprint arXiv:1409.8565}, 2014.

\bibitem{ma2015computational}
Z.~Ma and Y.~Wu, ``Computational barriers in minimax submatrix detection,''
  {\em The Annals of Statistics}, vol.~43, no.~3, pp.~1089--1116, 2015.

\bibitem{cai2015computational}
T.~T. Cai, T.~Liang, and A.~Rakhlin, ``Computational and statistical boundaries
  for submatrix localization in a large noisy matrix,'' {\em arXiv preprint
  arXiv:1502.01988}, 2015.

\bibitem{chen2014statistical}
Y.~Chen and J.~Xu, ``Statistical-computational tradeoffs in planted problems
  and submatrix localization with a growing number of clusters and
  submatrices,'' {\em arXiv preprint arXiv:1402.1267}, 2014.

\bibitem{barak2016noisy}
B.~Barak and A.~Moitra, ``Noisy tensor completion via the sum-of-squares
  hierarchy,'' in {\em 29th Annual Conference on Learning Theory},
  pp.~417--445, 2016.

\bibitem{hajek2015computational}
B.~E. Hajek, Y.~Wu, and J.~Xu, ``Computational lower bounds for community
  detection on random graphs.,'' in {\em COLT}, pp.~899--928, 2015.

\bibitem{bollobas1976cliques}
B.~Bollob{\'a}s and P.~Erd{\"o}s, ``Cliques in random graphs,'' in {\em
  Mathematical Proceedings of the Cambridge Philosophical Society}, vol.~80,
  pp.~419--427, Cambridge Univ Press, 1976.

\bibitem{feldman2013statistical}
V.~Feldman, E.~Grigorescu, L.~Reyzin, S.~Vempala, and Y.~Xiao, ``Statistical
  algorithms and a lower bound for detecting planted cliques,'' in {\em
  Proceedings of the forty-fifth annual ACM symposium on Theory of computing},
  pp.~655--664, ACM, 2013.

\bibitem{alon1998finding}
N.~Alon, M.~Krivelevich, and B.~Sudakov, ``Finding a large hidden clique in a
  random graph,'' {\em Random Structures and Algorithms}, vol.~13, no.~3-4,
  pp.~457--466, 1998.

\bibitem{ames2011nuclear}
B.~P. Ames and S.~A. Vavasis, ``Nuclear norm minimization for the planted
  clique and biclique problems,'' {\em Mathematical programming}, vol.~129,
  no.~1, pp.~69--89, 2011.

\bibitem{kuvcera1995expected}
L.~Kuvcera, ``Expected complexity of graph partitioning problems,'' {\em
  Discrete Applied Mathematics}, vol.~57, no.~2, pp.~193--212, 1995.

\bibitem{jerrum1992large}
M.~Jerrum, ``Large cliques elude the metropolis process,'' {\em Random
  Structures \& Algorithms}, vol.~3, no.~4, pp.~347--359, 1992.

\bibitem{feige2003probable}
U.~Feige and R.~Krauthgamer, ``The probable value of the lov{\'a}sz--schrijver
  relaxations for maximum independent set,'' {\em SIAM Journal on Computing},
  vol.~32, no.~2, pp.~345--370, 2003.

\bibitem{wedin1972perturbation}
P.-A. Wedin, ``Perturbation bounds in connection with singular value
  decomposition,'' {\em BIT Numerical Mathematics}, vol.~12, no.~1,
  pp.~99--111, 1972.

\bibitem{vershynin2010introduction}
R.~Vershynin, ``Introduction to the non-asymptotic analysis of random
  matrices,'' {\em arXiv preprint arXiv:1011.3027}, 2010.

\bibitem{koltchinskii2015optimal}
V.~Koltchinskii and D.~Xia, ``Optimal estimation of low rank density
  matrices,'' {\em Journal of Machine Learning Research}, vol.~16,
  pp.~1757--1792, 2015.

\bibitem{hush2005concentration}
D.~Hush and C.~Scovel, ``Concentration of the hypergeometric distribution,''
  {\em Statistics \& probability letters}, vol.~75, no.~2, pp.~127--132, 2005.

\bibitem{laurent2000adaptive}
B.~Laurent and P.~Massart, ``Adaptive estimation of a quadratic functional by
  model selection,'' {\em Annals of Statistics}, pp.~1302--1338, 2000.

\end{thebibliography}

\appendix

\newpage

\setcounter{page}{1}
\setcounter{section}{0}
\appendix

\begin{center}
	{\LARGE Supplement to ``Tensor SVD: Statistical and Computational Limits"
		\footnote{Anru Zhang is Assistant Professor, Department of Statistics, University of Wisconsin-Madison, Madison, WI 53706, E-mail: anruzhang@stat.wisc.edu; Dong Xia is Visiting Assistant Professor, Department of Statistics, University of Wisconsin-Madison, Madison, WI 53706, E-mail: dongxia@stat.wisc.edu.}}
	
	\bigskip\medskip
	
	Anru Zhang ~ and ~ Dong Xia
	\ \par University of Wisconsin-Madison
\end{center}

\begin{abstract}
	In this Supplement, we provide additional proofs for the main results and technical lemmas. 
\end{abstract}

\section{Additional Proofs}

\subsection{Proof of Theorem \ref{th:upper_bound_mle}}\label{sec:proof_upper_bound_mle}

We only need to prove upper bounds for $\hat{U}^{\bullet}_k$ and $\hat{\X}^\bullet$ as the ones for $\hat{U}^{\rm mle}_k$ and $\hat{\X}^{\rm mle}$ immediately follow. The proof of this theorem is fairly complicated. For convenience, we assume $\sigma=1$, denote $r = \max\{r_1, r_2, r_3\}$, then $r \leq C_0p^{1/2}$ according to the assumption. For any orthogonal columns, e.g. $U_k \in \mathbb{O}_{p_k, r_k}$, we note $U_{k\perp} \in \mathbb{O}_{p_k, p_k-r_k}$ as the orthogonal complement of $U_k$.  

Let $A\otimes B$ be the Kronecker product between matrices $A$ and $B$, ${\rm vec}(\cdot)$ be vectorization of matrices and tensors. Similarly as the proof of Theorems 1 and 2 in \cite{de2000best}, for any $V_k \in \mathbb{O}_{p_k, r_k}$, $k=1,2,3$, 
\begin{equation*}
\begin{split}
& \min_{\hat{\S}}\|\Y - \hat{\S}\times_1 V_1 \times_2 V_2\times_3 V_3\|_F^2 = \min_{\hat{\S}}\|{\rm vec}(\Y) - V_1\otimes V_2\otimes V_3 {\rm vec}(\hat{\S})\|_2^2\\
= & \min_{\hat{\S}}\|P_{(V_1\otimes V_2\otimes V_3)_\perp}{\rm vec}(\Y) + P_{V_1\otimes V_2\otimes V_3}{\rm vec}(\Y) - V_1\otimes V_2\otimes V_3 {\rm vec}(\hat{\S})\|_2^2\\
= & \|P_{(V_1\otimes V_2\otimes V_3)_\perp}{\rm vec}(\Y)\|_2^2 + \min_{\hat{\S}}\|P_{(V_1\otimes V_2\otimes V_3)}{\rm vec}(\Y) - V_1\otimes V_2\otimes V_3 {\rm vec}(\hat{\S})\|_2^2\\
= & \|P_{(V_1\otimes V_2\otimes V_3)_\perp}{\rm vec}(\Y)\|_2^2 = \|{\rm vec}(\Y)\|_2^2 - \|P_{(V_1\otimes V_2\otimes V_3)}{\rm vec}(\Y)\|_2^2 \\
= & \|\Y\|_{\rm F}^2 - \|\Y \times_1 V_1^\top \times_2 V_2^\top \times_3 V_3^\top\|_{\rm F}^2
\end{split}
\end{equation*}
where the inequality holds if and only if $\hat{\S} = \Y\times_1 V_1^\top\times_2 V_2^\top \times_3V_3^\top$. Therefore, we must have
$$\|\Y\times_1(\hat{U}_1^\bullet)^\top \times_2(\hat{U}_2^\bullet)^\top \times_3(\hat{U}_3^\bullet)^\top\|_{\rm F}^2 \geq \|\Y \times_1 U_1^\top \times_2 U_2^\top \times_3 U_3^\top\|_{\rm F}^2.$$
For convenience, we simply let $\hat{U}_k = \hat{U}_k^\bullet$ for $k=1,2,3$ and $\hat{\X} = \hat{\X}^\bullet$ throughout the proof of this theorem. Without loss of generality, we also assume that
\begin{equation}\label{eq:U_123}
U_k = \begin{bmatrix}
I_{r_k}\\
0_{(p_k-r_k)\times r_k}
\end{bmatrix}, \quad k=1,2,3.
\end{equation}
Such assumption will simplify our notation and make the proof easier to understand. This theorem will be shown by steps.
\begin{enumerate}[leftmargin=*]
	\item In this first step, we establish some basic probability bounds which will be used in the latter steps. We first let 
	\begin{equation}\label{eq:def_tilde_X}
	\tilde{\X} = \Y_{[1:r_1, 1:r_2, 1:r_3]} = \S + \Z_{[1:r_1, 1:r_2, 1:r_3]} \in \mathbb{R}^{r_1\times r_2\times r_3}.
	\end{equation}
	Then we first have
	\begin{equation}\label{ineq:U_F<=hat_U_F}
	\left\|\tilde{\X}\right\|_{\rm F}  = \left\|\Y \times_1 U_1^\top \times_2 U_2^\top  \times_3 U_3^\top\right\|_{\rm F} \leq \left\|\Y\times_1 \hat{U}_1^\top \times_2 \hat{U}_2^\top \times_3 \hat{U}_3^\top\right\|_{\rm F}.
	\end{equation}
	Next, we also note that $\mathcal{M}_1(\tilde{\X}) = \mathcal{M}_1(\S) + \mathcal{M}_1(\Z_{[1:r_1,1:r_2,1:r_3]}) \in \mathbb{R}^{r_1\times (r_2r_3)}$, i.e., the fixed matrix $\mathcal{M}_1(\S)$ plus an i.i.d. Gaussian matrix. Meanwhile, $\sigma_{\min}(\mathcal{M}_1(\S)) = \sigma_{r_1}(\mathcal{M}_1(\X)) \geq \lambda$. Now, by Lemma 4 in \cite{cai2016rate},
	\begin{equation*}
	P\left(\sigma_{r_1}^2(\mathcal{M}_1(\tilde{\X})) \geq (\lambda^2 + r_2r_3)(1-x)\right) \geq 1 - C\exp\left(Cr_1 - c(\lambda^2 + r_2r_3)x^2\wedge x\right), \quad \text{for any } x > 0.
	\end{equation*}
	Similar results also hold for $\sigma_{\min}(\mathcal{M}_2(\tilde{\X}))$ and $\sigma_{\min}(\mathcal{M}_3(\tilde{\X}))$. Let $x = 1/2$, note that $\lambda \geq C_{gap}p^{1/2}$, we have
	\begin{equation}\label{ineq:min_singular_M_k}
	\sigma_{\min}(\mathcal{M}_k(\tilde{\X})) \geq \frac{\lambda}{2^{1/2}}, \quad k= 1, 2, 3
	\end{equation}
	with probability at least $1 - C\exp(-cp)$ for large enough constant $C_{gap}>0$. Additionally, by Lemma \ref{lm:random_tensor_max_projection}, we have
	\begin{equation}\label{ineq:ZtimesV1V2V3}
	\max_{\substack{V_k\in \mathbb{O}_{p_k, r_k}\\k=1,2,3}}\left\|\Z\times_1 V_1^\top \times_2 V_2^\top \times_3 V_3^\top\right\|_{\rm F} \leq C\sqrt{pr},
	\end{equation}
	with probability at least $1 - C\exp(-cp)$.
	\item 
	%
	In the following Steps 2 and 3, we temporarily ignore the randomness of $\Z$ and the definition of $\hat{U}_1, \hat{U}_2, \hat{U}_3$ as the the estimators with high likelihood values. Instead we only assume $\X, \Z \in \mathbb{R}^{p_1\times p_2\times p_3}$, $\X = \S\times_1U_1\times_2U_2\times_3U_3$, and $\hat{U}_k \in \mathbb{O}_{p_k, r_k}$ for $k=1,2,3$ satisfy \eqref{ineq:U_F<=hat_U_F}, \eqref{ineq:min_singular_M_k}, and \eqref{ineq:ZtimesV1V2V3}. By Lemma 1 in \cite{cai2016rate}, $\|\sin\Theta(\hat{U}_k, U_k)\|_{\rm F} = \|U^\top_{k\perp} \hat{U}_k\|_{\rm F}$. Our goal in Steps 2 is to show under such setting, one must has
	\begin{equation}\label{ineq:goal-step23}
	\|\sin\Theta(\hat{U}_k, U_k)\|_{\rm F} = \|U^\top_{k\perp} \hat{U}_k\|_{\rm F} \leq C\frac{\sqrt{pr}}{\lambda}, \quad k = 1, 2, 3.
	\end{equation}
	To simplify our notations, we first perform spectral transformation on $(\hat{U}_k)_{[1:r_k, :]}$ and $(\hat{U}_k)_{[(r_k+1):p_k, :]}$. To be specific, let $(\hat{U}_k)_{[1:r_k, :]} = \bar{U}_k\bar{\Sigma}_k\bar{V}_k^\top$ be the singular value decomposition, and $(\hat{U}_k)_{[(r_k+1):p_k, :]} \bar{V}_k = \bar{Q}_k\bar{R}_k$ be the full QR decomposition (so that $\bar{Q}_k$ is a square orthogonal matrix). Then we transform 
	\begin{equation*}
	\begin{split}
	\begin{bmatrix}
	(\hat{U}_k)_{[1:r_k, :]}\\
	(\hat{U}_k)_{[(r_k+1):p_k, :]}\\
	\end{bmatrix}\quad & \Rightarrow\quad \begin{bmatrix}
	\bar{U}_k^\top & \\
	& \bar{Q}_k^\top \\
	\end{bmatrix}\cdot \begin{bmatrix}
	(\hat{U}_k)_{[1:r_k, :]}\\
	(\hat{U}_k)_{[(r_k+1):p, :]}
	\end{bmatrix} \bar{V}_k = \begin{bmatrix}
	\bar{U}_k^\top (\hat{U}_k)_{[1:r_k, :]}\bar{V}_k\\
	\bar{Q}_k^\top (\hat{U}_k)_{[(r_k+1):p_k, :]}\bar{V}_k
	\end{bmatrix} = \begin{bmatrix}
	\bar{\Sigma}_k\\
	\bar{R}_k
	\end{bmatrix}, \\
	\Z\quad & \Rightarrow \quad \Z \times_1 \begin{bmatrix}
	\bar{U}_1 & \\
	& \bar{Q}_1 \\
	\end{bmatrix} \times_2 \begin{bmatrix}
	\bar{U}_2 & \\
	& \bar{Q}_2 \\
	\end{bmatrix} \times_3 \begin{bmatrix}
	\bar{U}_3 & \\
	& \bar{Q}_3 \\
	\end{bmatrix}, \\
	\X \quad & \Rightarrow \quad \X \times_1 \begin{bmatrix}
	\bar{U}_1 & \\
	& \bar{Q}_1 \\
	\end{bmatrix} \times_2 \begin{bmatrix}
	\bar{U}_2 & \\
	& \bar{Q}_2 \\
	\end{bmatrix} \times_3 \begin{bmatrix}
	\bar{U}_3 & \\
	& \bar{Q}_3 \\
	\end{bmatrix}. 
	\end{split}
	\end{equation*}
	We can check that \eqref{ineq:U_F<=hat_U_F}, \eqref{ineq:min_singular_M_k}, \eqref{ineq:ZtimesV1V2V3} still hold after this transformation.
	Suppose $\diag(R_1) = (a_1,\ldots, a_r)$. Since $\bar{\Sigma}_1$ is diagonal and $\bar{R}_1$ is upper diagonal, $\begin{bmatrix}
	\bar{\Sigma}_1\\
	\bar{R}_1
	\end{bmatrix}$ is orthogonal, we must have all off-diagonal entries of $\bar{R}_k$ are zero, and $\bar{\Sigma}_1 = \diag(\sqrt{1-a_1^2}, \ldots, \sqrt{1-a_r^2})$. For convenient we also denote 
	\begin{equation}\label{eq:def_a_i^0^1}
	a_i^{(0)} = \sqrt{1-a_i^2}, a_i^{(1)} = a_i, b_j^{(0)} = \sqrt{1-b_j^2}, b_j^{(1)} = b_j, c_k^{(0)} = \sqrt{1-c_k^2}, c_k^{(1)} = c_k.
	\end{equation} 
	Therefore, without loss of generality we can assume there exist real numbers $0\leq a_i, b_j, c_k \leq 1$ such that
	\begin{equation}\label{eq:U_123_simple_form}
	\begin{split}
	\hat{U}_1 = \begin{bmatrix}
	a_1^{(0)} & & \\
	& \ddots & \\
	& & a_{r_1}^{(0)}\\
	a_1^{(1)} & & \\
	& \ddots & \\
	& & a_{r_1}^{(1)}\\
	& 0_{p_1-2r_1, r_1} &
	\end{bmatrix}, \hat{U}_2 = \begin{bmatrix}
	b_1^{(0)} & & \\
	& \ddots & \\
	& & b_{r_2}^{(0)}\\
	b_1^{(1)} & & \\
	& \ddots & \\
	& & b_{r_2}^{(1)}\\
	& 0_{p_2-2r_2, r_2} &
	\end{bmatrix},
	\hat{U}_3 = \begin{bmatrix}
	c_1^{(0)} & & \\
	& \ddots & \\
	& & c_{r_3}^{(0)}\\
	c_1^{(1)} & & \\
	& \ddots & \\
	& & c_{r_3}^{(1)}\\
	& 0_{p_3-2r_3, r_3} &
	\end{bmatrix}.
	\end{split}
	\end{equation}
	where $0_{p_k-2r_k, r_k}$ represents the zero matrix with dimension $(p_k-2r_k)$-by-$r_k$. By the form of $U_1, U_2, U_3$ in \eqref{eq:U_123}, we must have
	\begin{equation}\label{ineq:sin-theta-a-b-c}
	\max_{k=1,2,3} \|\sin\Theta(\hat{U}_k, U_k)\|_{\rm F} = \max\left\{(\sum_{i=1}^{r_1} a_i^2)^{1/2}, (\sum_{i=1}^{r_2} b_i^2)^{1/2}, (\sum_{i=1}^{r_3} c_i^2)^{1/2}\right\}.
	\end{equation}
	In order to show \eqref{ineq:goal-step23} we only need to prove that $\max\left\{(\sum_{i=1}^{r_1} a_i^2)^{1/2}, (\sum_{i=1}^{r_2} b_i^2)^{1/2}, (\sum_{i=1}^{r_3} c_i^2)^{1/2}\right\} \leq C\sqrt{pr}.$
	
	Next, we decompose the noise tensor $\Z$ to the following eight pieces,
	\begin{equation}
	\Z^{(t_1t_2t_3)} = \Z_{[(t_1r_1+1):(t_1r_1 + r_1), (t_2r_2+1):(t_2r_2 + r_2), (t_3r_3+1):(t_3r_3 + r_3)]}, \quad t_1, t_2, t_3 \in\{0, 1\}.
	\end{equation}
	By \eqref{ineq:ZtimesV1V2V3}, we know $\|\Z^{(t_1t_2t_3)}\|_{\rm F} \leq C\sqrt{pr}$, $t_1, t_2, t_3 \in\{0, 1\}$. Based on the form of $\hat{U}_1, \hat{U}_2, \hat{U}_3$ \eqref{eq:U_123_simple_form}, we have
	\begin{equation*}
	\begin{split}
	& \left(\Y\times_1 \hat{U}_1 \times_2 \hat{U}_2 \times_3 \hat{U}_3\right)_{ijk}\\
	= & \tilde{X}_{ijk} a_i^{(0)}b_j^{(0)}c_k^{(0)} + Z^{(100)}_{ijk}a_i^{(1)}b_j^{(0)}c_k^{(0)} + Z^{(010)}_{ijk}a_i^{(0)}b_j^{(1)}c_k^{(0)}  + Z^{(001)}_{ijk}a_i^{(0)}b_j^{(0)}c_k^{(1)}\\ 
	& + Z^{(110)}_{ijk}a_i^{(1)}b_j^{(1)}c_k^{(0)} + Z^{(101)}_{ijk}a_i^{(1)}b_j^{(0)}c_k^{(1)} + Z^{(011)}_{ijk}a_i^{(0)}b_j^{(1)}c_k^{(1)} + Z^{(111)}_{ijk}a_i^{(1)}b_j^{(1)}c_k^{(1)}\\
	= & \tilde{X}_{ijk}a_i^{(0)}b_j^{(0)}c_k^{(0)} + \sum_{\substack{t_1, t_2, t_3 \in \{0, 1\}\\t_1, t_2, t_3 \text{ are not all 0}}} Z_{ijk}^{(t_1t_2t_3)}a_i^{(t_1)}b_j^{(t_2)}c_k^{(t_3)}.
	\end{split}
	\end{equation*}
	Therefore,
	\begin{equation*}
	\begin{split}
	0 \overset{\eqref{ineq:U_F<=hat_U_F}}{\leq} & \left\|\Y\times_1 \hat{U}_1 \times_2 \hat{U}_2 \times_3 \hat{U}_3\right\|_{\rm F}^2 - \left\|\Y\times_1 U_1 \times_2 U_2 \times_3 U_3 \right\|_{\rm F}^2\\ 
	= & \sum_{i,j,k=1}^{r_1, r_2, r_3}\left[\left(\tilde{X}_{ijk}a_i^{(0)}b_j^{(0)}c_k^{(0)} + \sum_{\substack{t_1, t_2, t_3 \in \{0, 1\}\\t_1, t_2, t_3 \text{ are not all 0}}} Z_{ijk}^{(t_1t_2t_3)}a_i^{(t_1)}b_j^{(t_2)}c_k^{(t_3)}\right)^2 - \tilde{X}_{ijk}^2\right]\\
	\leq & \sum_{i,j,k=1}^{r_1, r_2, r_3}\tilde{X}_{ijk}^2\left[\left(a_i^{(0)}b_j^{(0)}c_k^{(0)}\right)^2 - 1\right]\\ 
	& + 2\sum_{\substack{t_1, t_2, t_3 \in \{0, 1\}\\t_1, t_2, t_3 \text{ are not all 0}}}\sum_{i,j,k=1}^{r_1, r_2, r_3} \tilde{X}_{ijk} a_i^{(0)}b_j^{(0)}c_k^{(0)}  Z_{ijk}^{(t_1t_2t_3)}a_i^{(t_1)}b_j^{(t_2)}c_k^{(t_3)}\\
	& + 7\sum_{\substack{t_1, t_2, t_3 \in \{0, 1\}\\t_1, t_2, t_3 \text{ are not all 0}}} \sum_{i,j,k=1}^{r_1, r_2, r_3} \left( Z_{ijk}^{(t_1t_2t_3)}a_i^{(t_1)}b_j^{(t_2)}c_k^{(t_3)}\right)^2\\
	\leq & \sum_{i,j,k=1}^{r_1, r_2, r_3}\tilde{X}_{ijk}^2\left[(1-a_i^2)(1-b_j^2)(1-c_k^2) - 1\right]\\ 
	& + 63\max_{\substack{t_1, t_2, t_3 \in \{0, 1\}\\t_1, t_2, t_3 \text{ are not all 0}}} \sum_{i,j,k=1}^{r_1,r_2,r_3}\Bigg\{\tilde{X}_{ijk}\sqrt{1-a_i^2}\sqrt{1-b_j^2}\sqrt{1-c_k^2}  Z_{ijk}^{(t_1t_2t_3)}a_i^{(t_1)}b_j^{(t_2)}c_k^{(t_3)}, \\
	& \quad \sum_{i,j,k=1}^{r_1, r_2, r_3} \left( Z_{ijk}^{(t_1t_2t_3)}a_i^{(t_1)}b_j^{(t_2)}c_k^{(t_3)}\right)^2\Bigg\}.
	\end{split}
	\end{equation*}
	By the inequality above, one of the following inequalities must hold for some $t_1, t_2, t_3\in \{0, 1\}$ and $t_1, t_2, t_3$ are not all 0:
	\begin{equation}\label{ineq:must_hold1}
	\begin{split}
	0\leq & \sum_{i,j,k=1}^{r_1, r_2, r_3}\tilde{X}_{ijk}^2\left[(1-a_i^2)(1-b_j^2)(1-c_k^2) - 1\right]\\
	& + 63\tilde{X}_{ijk}\sqrt{1-a_i^2}\sqrt{1-b_j^2}\sqrt{1-c_k^2}  Z_{ijk}^{(t_1t_2t_3)}a_i^{(t_1)}b_j^{(t_2)}c_k^{(t_3)},
	\end{split}
	\end{equation}
	\begin{equation}\label{ineq:must_hold2}
	0 \leq \sum_{i,j,k=1}^{r_1, r_2, r_3}\tilde{X}_{ijk}^2\left[(1-a_i^2)(1-b_j^2)(1-c_k^2) - 1\right] + 63 \sum_{i,j,k=1}^{r_1, r_2, r_3} \left( Z_{ijk}^{(t_1t_2t_3)}a_i^{(t_1)}b_j^{(t_2)}c_k^{(t_3)}\right)^2.
	\end{equation}
	
	\item Next we discuss in two different situations: \eqref{ineq:must_hold1} or \eqref{ineq:must_hold2} hold.
	\begin{enumerate}[leftmargin=*]
		\item When \eqref{ineq:must_hold1} holds, we first assume $t_1=1, t_2 = t_3 = 0$ as the other situations follow similarly. Then
		\begin{equation*}
		\begin{split}
		0\leq & \sum_{i,j,k=1}^{r_1, r_2, r_3}\tilde{X}_{ijk}^2\left[(1-a_i^2) - 1\right] + 63\sum_{i, j, k=1}^{r_1, r_2, r_3}|\tilde{X}_{ijk}||Z_{ijk}^{(100)}|a_i \quad (\text{since } 1\leq a_i, b_j, c_k\leq 1)\\
		\leq & -\sum_{i,j,k=1}^{r_1, r_2, r_3} \tilde{X}_{ijk}^2 a_i^2 + 63 \left(\sum_{i,j,k=1}^{r_1, r_2, r_3}\tilde{X}_{ijk}^2a_i^2\right)^{1/2}\left(\sum_{i,j,k=1}^{r_1, r_2, r_3}(Z_{ijk}^{(100)})^2\right)^{1/2}, \quad \text{(Cauchy-Schwarz inequality)}\\
		\end{split}
		\end{equation*}
		Thus, 
		\begin{equation}\label{ineq:useful_later}
		\sum_{i,j,k=1}^{r_1, r_2, r_3} \tilde{X}_{ijk}^2 a_i^2 \leq C\sum_{i,j,k=1}^{r_1, r_2, r_3}(Z_{ijk}^{(100)})^2. 
		\end{equation}
		Additionally, since 
		$$\sum_{j,k=1}^{r_2, r_3} \tilde{X}_{ijk}^2 = \left\|\left(\mathcal{M}_1(\tilde{\X})\right)_{[i, :]}\right\|_2\geq  \sigma_{\min}^2(\mathcal{M}_1(\tilde{\X}))$$ 
		and \eqref{ineq:min_singular_M_k}, we have
		\begin{equation*}
		\sum_{i=1}^{r_1} a_i^2 \frac{\lambda^2}{2} \leq \sum_{i,j,k=1}^{r_1, r_2, r_3} \tilde{X}_{ijk}^2 a_i^2 \leq 63\sum_{i,j,k=1}^{r_1, r_2, r_3} \left(Z_{ijk}^{(100)}\right)^2 = 63\|\Z^{(100)}\|_{\rm F}^2\overset{\eqref{ineq:ZtimesV1V2V3}}{\leq} Cpr,
		\end{equation*}
		which means $\|\sin\Theta(\hat{U}_1, U_1)\|_{\rm F} \overset{\eqref{ineq:sin-theta-a-b-c}}{=} \sqrt{\sum_{i=1}^{r_1} a_i^2} \leq C\sqrt{pr}/\lambda$. On the other hand, by \eqref{ineq:must_hold1},
		\begin{equation}\label{ineq:useful_later2}
		\begin{split}
		0\leq & \sum_{i,j,k=1}^{r_1, r_2, r_3} \tilde{X}_{ijk}^2[(1-a_i^2)(1-b_j^2)-1] + 63 \sum_{i,j,k=1}^{r_1, r_2, r_3} |\tilde{X}_{ijk}||Z_{ijk}^{(100)}|a_i \sqrt{1-a_i^2}(1-b_j^2)\\
		& \quad  \quad \text{(by Algorithmic-geometric inequality)}\\
		\leq & \sum_{i,j,k=1}^{r_1, r_2, r_3} \tilde{X}_{ijk}^2[-b_j^2-a_i^2 + a_i^2b_j^2] + \sum_{i,j,k=1}^{r_1, r_2, r_3} \left(\tilde{X}_{ijk}^2a_i^2(1-b_j^2) + \frac{63^2}{4} (Z_{ijk}^{(100)})^2(1-a_i^2)(1-b_j^2)\right)\\
		\leq & -\sum_{i,j,k=1}^{r_1, r_2, r_3} \tilde{X}_{ijk}^2 b_j^2 + \frac{63^2}{4}\sum_{i,j,k}^{r}(Z_{ijk}^{(100)})^2 \overset{\eqref{ineq:min_singular_M_k}}{\leq} -\frac{\lambda^2}{2} \sum_{j=1}^{r_2} b_j^2 + \frac{63^2}{4} \|Z^{(100)}\|_{\rm F}^2.
		\end{split}
		\end{equation}
		Therefore, 
		$$\|\sin\Theta(\hat{U}_2, U_2)\|_{\rm F} \overset{\eqref{ineq:sin-theta-a-b-c}}{=} \sqrt{\sum_{j=1}^{r_2} b_j^2} \leq C\|\Z^{(100)}\|_{\rm F}/\lambda\overset{\eqref{ineq:ZtimesV1V2V3}}{\leq}  C\sqrt{pr}/\lambda.$$
		By symmetry, one can also show that $\|\sin\Theta(\hat{U}_3, U_3)\|_{\rm F} \leq C\sqrt{pr}/\lambda$. In summary, we must have
		$$\max_{k=1, 2, 3} \|\sin\Theta(\hat{U}_k, U_k)\|_{\rm F} \leq C\frac{\sqrt{pr}}{\lambda}$$
		for some constant $C>0$ when \eqref{ineq:must_hold1} holds.
		\item When \eqref{ineq:must_hold2} holds for some $t_1, t_2, t_3\in \{0, 1\}$, 
		\begin{equation*}
		\begin{split}
		0\leq & \sum_{i,j,k=1}^{r_1, r_2, r_3} \tilde{X}_{ijk}^2[(1-a_i^2)(1-b_j^2)(1-c_k^2) - 1] + 63\sum_{i,j,k=1}^{r_1, r_2, r_3} (Z_{ijk}^{(t_1t_2t_3)})^2(a_i^{(t_1)})^2(b_j^{(t_2)})^2(c_k^{(t_3)})^2\\
		\leq & \sum_{i,j,k=1}^{r_1, r_2, r_3} \tilde{X}_{ijk}^2[(1-a_i^2) - 1] + 63\sum_{i,j,k=1}^{r_1, r_2, r_3} (Z_{ijk}^{(t_1t_2t_3)})^2\\
		\leq & -\sum_{i,j,k}^{r_1, r_2, r_3} a_i^2 \tilde{X}_{ijk}^2 + 63\|\Z^{(t_1t_2t_3)}\|_{\rm F}^2 \leq \frac{\lambda^2}{2}\sum_{i,j,k=1}^{r_1, r_2, r_3} a_i^2 + Cpr,
		\end{split}
		\end{equation*}
		which means $\|\sin\Theta(\hat{U}_1, U_1)\|_{\rm F} \overset{\eqref{ineq:sin-theta-a-b-c}}{=} \sqrt{\sum_{i=1}^{r_1} a_i^2} \leq C\sqrt{pr}/\lambda$. One can similarly prove the parallel results for $\|\sin\Theta(\hat{U}_2, U_2)\|_{\rm F}$ and $\|\sin\Theta(\hat{U}_3, U_3)\|_{\rm F}$. In summary, we also have
		$$\max_{k=1, 2, 3} \|\sin\Theta(\hat{U}_k, U_k)\|_{\rm F} \leq C\frac{\sqrt{pr}}{\lambda}$$
		for some constant $C>0$ when \eqref{ineq:must_hold2} holds.	
	\end{enumerate}
	To sum up, we have the derived perturbation bound: under \eqref{ineq:U_F<=hat_U_F}, \eqref{ineq:min_singular_M_k}, \eqref{ineq:ZtimesV1V2V3}, one must have $\|\sin\Theta(\hat{U}_k, U_k)\|_{\rm F}\leq C\sqrt{pr}/\lambda$.
	\item Next we consider the recovery loss for $\hat{\X}$. Similarly as Steps 2-3, we temporarily ignore the randomness of $\Z$, and the definition of $\hat{U}_1, \hat{U}_2, \hat{U}_3$ as the estimators with high likelihood values in this step. We aim to prove
	$$\left\|\hat{\X} - \X\right\|_{\rm F} \leq C\sqrt{pr},$$
	under the assumptions of \eqref{ineq:U_F<=hat_U_F}, \eqref{ineq:min_singular_M_k}, and \eqref{ineq:ZtimesV1V2V3}. First, without loss of generality we can assume $\hat{U}_1, \hat{U}_2, \hat{U}_3$ have the simple form \eqref{eq:U_123_simple_form}. Based on the structure of $\hat{U}_1, \hat{U}_2, \hat{U}_3$, we know
	$$P_{\hat{U}_1} = \begin{bmatrix}
	(a_1^{(0)})^2 & & & a_1^{(1)}a_1^{(0)} & & &\\
	& \ddots & & & \ddots & & 0_{r_1, p_1-2r_1} \\
	& & (a_{r_1}^{(0)})^2 & & & a_{r_1}^{(1)}a_{r_1}^{(0)} & \\
	a_1^{(1)}a_1^{(0)} & & & (a_1^{(1)})^2 & & & \\
	& \ddots & & & \ddots & & 0_{r_1, p_1-2r_1}\\
	& & a_r^{(1)}a_{r_1}^{(0)} & & & (a_{r_1}^{(1)})^2 & \\	 
	& 0_{p_1-2r_1, r_1} & & & 0_{p_1-2r_1, r_1} & & 0_{p_1-2r_1, p_1-2r_1} &
	\end{bmatrix}, $$
	while $P_{\hat{U}_2}$ and $P_{\hat{U}_3}$ can be written in similar forms.
	We have the following decomposition for $\|\hat{\X} - \X\|_{\rm F}$,
	\begin{equation}\label{eq:hat_X-X_derivation}
	\begin{split}
	& \left\|\hat{\X} - \X\right\|_{\rm F} = \left\|\Y\times_1 \Proj_{\hat{U}_1} \times_2 \Proj_{\hat{U}_2} \times_3 \Proj_{\hat{U}_3} - \S\times_1 U_1\times_2 U_2\times_3 U_3\right\|_{\rm F}\\
	\leq & \left\|\Z \times_1 P_{\hat{U}_1} \times_2 P_{\hat{U}_2} \times_3 P_{\hat{U}_3}\right\|_{\rm F} + \left\|\X \times_1 P_{\hat{U}_1} \times_2 P_{\hat{U}_2} \times_3 P_{\hat{U}_3} - \S\times_1 U_1\times_2 U_2 \times_3 U_3\right\|_{\rm F}\\
	\leq & C\sqrt{pr} + \left\|\S\times_1 (P_{\hat{U}_1} U_1^\top) \times_2 (P_{\hat{U}_2} U_2^\top) \times_3 (P_{\hat{U}_3} U_3^\top) - \S \times_1U_1 \times_2 U_2 \times_3 U_3 \right\|_{\rm F}\\
	\leq & C\sqrt{pr} + \left\|\tilde{\X}\times_1 (P_{\hat{U}_1} U_1^\top) \times_2 (P_{\hat{U}_2} U_2^\top) \times_3 (P_{\hat{U}_3} U_3^\top) - \tilde{\X} \times_1U_1 \times_2 U_2 \times_3 U_3\right\|_{\rm F}\\
	& + \left\|\Z_{[1:r,1:r,1:r]}\times_1 (P_{\hat{U}_1} U_1^\top) \times_2 (P_{\hat{U}_2} U_2^\top) \times_3 (P_{\hat{U}_3} U_3^\top)\right\| + \left\|\Z_{[1:r, 1:r, 1:r]}\times_1U_1 \times_2 U_2 \times_3 U_3\right\|_{\rm F}\\
	\leq & C\sqrt{pr} + \left\|\tilde{\X}\times_1 (P_{\hat{U}_1} U_1^\top) \times_2 (P_{\hat{U}_2} U_2^\top) \times_3 (P_{\hat{U}_3} U_3^\top) - \tilde{\X} \times_1U_1 \times_2 U_2 \times_3 U_3\right\|_{\rm F}.
	\end{split}
	\end{equation}
	Based on the form of $U_1, U_2, U_3, \hat{U}_1, \hat{U}_2, \hat{U}_3$, we have
	\begin{equation*}
	\begin{split}
	& \left\|\tilde{\X}\times_1 (P_{\hat{U}_1} U_1^\top) \times_2 (P_{\hat{U}_2} U_2^\top) \times_3 (P_{\hat{U}_3} U_3^\top) - \tilde{\X} \times_1U_1 \times_2 U_2 \times_3 U_3\right\|_{\rm F}^2\\
	= & \sum_{i,j,k=1}^{r_1, r_2, r_3} \tilde{X}_{ijk}^2 \left(\left(a_i^{(0)}b_j^{(0)}c_k^{(0)}\right)^2 - 1\right)^2 + \sum_{\substack{t_1, t_2, t_3 \in \{0, 1\}\\t_1, t_2, t_3 \text{ are not all 0}}} \sum_{i,j,k=1}^{r_1, r_2, r_3} \tilde{X}_{ijk}^2 \left(a_i^{(0)}b_j^{(0)}c_k^{(0)} a_i^{(t_1)} b_j^{(t_2)}c_k^{(t_0)}\right)^2\\
	= & \sum_{i,j,k=1}^{r_1, r_2, r_3}\tilde{X}_{ijk}^2 \left(a_i^{(0)}b_j^{(0)}c_k^{(0)}\right)^2\left((a_i^{(0)})^2 + (a_i^{(1)})^2\right) \left((b_j^{(0)})^2 + (b_j^{(1)})^2\right) \left((c_k^{(0)})^2 + (c_k^{(1)})^2\right)\\
	& + \sum_{i, j, k=1}^{r_1, r_2, r_3} \tilde{X}_{ijk}^2\left(-2\left(a_i^{(0)}b_j^{(0)}c_k^{(0)}\right)^2 + 1\right).
	\end{split}
	\end{equation*}
	Recall the actual values of $a_i^{(0)}, a_i^{(1)}, b_j^{(0)}, b_j^{(1)}, c_k^{(0)}, c_k^{(1)}$ in \eqref{eq:def_a_i^0^1}, we further have
	\begin{equation}\label{ineq:mle_step3_intermediate_process_1}
	\begin{split}
	& \left\|\tilde{\X}\times_1 (P_{\hat{U}_1} U_1^\top) \times_2 (P_{\hat{U}_2} U_2^\top) \times_3 (P_{\hat{U}_3} U_3^\top) - \tilde{\X} \times_1U_1 \times_2 U_2 \times_3 U_3\right\|_{\rm F}^2\\
	& \sum_{i,j,k=1}^{r_1, r_2, r_3} \tilde{X}_{ijk}^2 \left(a_i^{(0)}b_j^{(0)}c_k^{(0)}\right)^2 + \sum_{i, j, k=1}^{r_1, r_2, r_3} \tilde{X}_{ijk}^2\left(-2\left(a_i^{(0)}b_j^{(0)}c_k^{(0)}\right)^2 + 1\right)\\
	= & \sum_{i, j, k=1}^{r_1, r_2, r_3} \tilde{X}_{ijk}^2\left(1-(1-a_i^2)(1-b_j^2)(1-c_k^2)\right).
	\end{split}
	\end{equation}
	By the analysis in Step 2, we know under \eqref{ineq:U_F<=hat_U_F}, \eqref{ineq:min_singular_M_k}, \eqref{ineq:ZtimesV1V2V3}, at least one of \eqref{ineq:must_hold1} and \eqref{ineq:must_hold2} must hold for some $(t_1, t_2, t_3)\in \{0, 1\}^3\backslash\{(0,0,0)\}$. Again, we discuss in two different situations to show no matter which of \eqref{ineq:must_hold1} or \eqref{ineq:must_hold1} happen, we must have
	\begin{equation}\label{ineq:goal_recovery_X_middle}
	\sum_{i,j,k=1}^{r_1, r_2, r_3} \tilde{X}_{ijk}^2\left[1 - (1-a_i^2)(1-b_j^2)(1-c_k^2)\right] \leq Cpr.
	\end{equation}
	\begin{enumerate}[leftmargin=*]
		\item When \eqref{ineq:must_hold1} holds, we again assume $t_1 = 1, t_2 = t_3=0$ as the other situations follow similarly. Particularly, we have shown in Step 2 (a), \eqref{ineq:useful_later} and \eqref{ineq:useful_later2}, 
		$$\sum_{i,j,k=1}^{r_1, r_2, r_3}\tilde{X}_{ijk}^2a_i^2 \leq C \sum_{i,j,k=1}^{r_1, r_2, r_3} (Z_{ijk}^{(100)})^2, \quad \sum_{i,j,k=1}^{r_1, r_2, r_3}\tilde{X}_{ijk}^2b_j^2 \leq C \sum_{i,j,k=1}^{r_1, r_2, r_3} (Z_{ijk}^{(100)})^2 $$
		Clearly, $\sum_{i,j,k=1}^{r_1, r_2, r_3} \tilde{X}_{ijk}^2c_k^2 \leq C \sum_{i,j,k=1}^{r_1, r_2, r_3} (Z_{ijk}^{(100)})^2$ can be derived by symmetry. Then,
		\begin{equation*}
		\begin{split}
		& \sum_{i,j,k=1}^{r_1, r_2, r_3}\tilde{X}_{ijk}^2 \left(1 - (1-a_i^2)(1-b_j^2)(1-c_k^2)\right)\\ 
		= & \sum_{i,j,k=1}^{r_1, r_2, r_3}\tilde{X}_{ijk}^2 \left(a_i^2+b_j^2+c_k^2 - a_i^2b_j^2-a_i^2c_k^2-b_j^2c_k^2 + a_i^2b_j^2c_k^2\right) \\
		\leq & \sum_{i,j,k=1}^{r_1, r_2, r_3} \tilde{X}_{ijk}^2(a_i^2+b_j^2+c_k^2) \quad \text{(since $0\leq a_i, b_j, c_k \leq 1$)}\\
		\leq & C\sum_{i,j,k=1}^{r_1,r_2,r_3}(Z_{ijk}^{(100)})^2 \leq Cpr.
		\end{split}
		\end{equation*}
		\item When \eqref{ineq:must_hold2} holds, one has
		\begin{equation*}
		\begin{split}
		& \sum_{i,j,k=1}^{r_1, r_2, r_3} \tilde{X}_{ijk}^2\left(1 - (1-a_i^2)(1-b_j^2)(1-c_k^2)\right) \leq 63\sum_{ijk=1}^{r_1, r_2, r_3} \left(Z_{ijk}^{(t_1t_2t_3)}a_i^{(t_1)}b_j^{(t_2)}c_k^{(t_3)}\right)^2\\
		\leq & C\sum_{i, j, k=1}^{r_1, r_2, r_3} \left(Z_{ijk}^{(t_1t_2t_3)}\right)^2 = C\|\Z^{(t_1t_2t_3)}\|_{\rm F} \leq Cpr.
		\end{split}
		\end{equation*}
	\end{enumerate}
	In summary of Cases (a)(b), we must have \eqref{ineq:goal_recovery_X_middle}. 
	Combining \eqref{ineq:mle_step3_intermediate_process_1}, \eqref{eq:hat_X-X_derivation}, and \eqref{ineq:goal_recovery_X_middle}, we have
	shown	
	\begin{equation}
	\|\hat{\X} - \X\|_{\rm F} \leq C\sqrt{pr},\quad \text{under \eqref{ineq:U_F<=hat_U_F}, \eqref{ineq:min_singular_M_k}, \eqref{ineq:ZtimesV1V2V3}.}
	\end{equation}

	\item We finalize the proof for Theorem \ref{th:upper_bound_mle} in this step. We let $Q = \{\text{\eqref{ineq:U_F<=hat_U_F}, \eqref{ineq:min_singular_M_k}, \eqref{ineq:ZtimesV1V2V3} all hold}\}$. By Step 1, $P(Q) \geq 1-C\exp(-cp)$; by Steps 2-4, one has $\|\sin\Theta(\hat{U}_k, U_k)\|_{\rm F} \leq C\sqrt{pr}/\lambda, k=1,2,3$, and $\|\hat{\X} - \X\|_{\rm F}\leq C\sqrt{pr}$ under $Q$. The rest of the proof is essentially the same as the Step 4 in the proof of Theorem \ref{th:upper_bound_strong}.
	
	Since $\hat{\X}$ is a projection of $\Y$ by definition, so $\|\hat{\X}\|_{\rm F}\leq \|\Y\|_{\rm F} \leq \|\X\|_{\rm F} + \|\Z\|_{\rm F}.$ Then we have the following upper bound for 4-th moment of recovery error,
	\begin{equation*}
	\begin{split}
	& \mathbb{E} \|\hat{\X} - \X\|_{\rm F}^4\leq C\left(\mathbb{E}\|\hat{\X}\|_{\rm F}^4 + \|\X\|_{\rm F}^4\right) \leq C\|\X\|_{\rm F}^4 + C\mathbb{E} \|\Z\|_{\rm F}^4 \\
	\leq & C\exp(c_0p) + C\mathbb{E} \left(\chi^2_{p^3}\right)^2 = C\exp(c_0p) + Cp^6.
	\end{split}
	\end{equation*}
	The we have the following upper bound for the risk of $\hat{\X}$,
	\begin{equation*}
	\begin{split}
	& \mathbb{E}\|\hat{\X} - \X\|_{\rm F}^2 = \mathbb{E}\|\hat{\X} - \X\|_{\rm F}^21_{Q} + \mathbb{E}\|\hat{\X} - \X\|_{\rm F}^21_{Q^c} = Cpr +\sqrt{\mathbb{E} \|\hat{\X} - \X\|_{\rm F}^4 \mathbb{E}1_{Q^c}} \\
	\leq & Cpr + C\exp\left((c_0-c)p\right) + Cp^6\exp(-cp).
	\end{split}
	\end{equation*}
	Thus, one can select $c_0 < c$ to ensure that 
	$$\mathbb{E}\|\hat{\X} - \X\|_{\rm F}^2\leq Cpr \leq p_1r_1+p_2r_2+p_3r_3.$$
	Additionally, when $\sigma_{\min}(\mathcal{M}_k(\X)) \geq \lambda$, we have $\|\X\|_{\rm F}^2 = \|\mathcal{M}_k(\X)\|_{\rm F}^2 \geq r_k\lambda^2$ for $k=1,2,3$, which implies $\|\X\|_{\rm F}^2 \geq Cr\lambda^2$. Thus we also have
	$$\E\frac{\|\hat{\X} - \X\|_{\rm F}^2}{\|\X\|_{\rm F}^2} \leq \frac{p_1+p_2+p_3}{\lambda^2}.$$
	
	Now we consider the Frobenius $\sin\theta$ norm risk for $\hat{U}_k$. Since $\sin\Theta(\hat{U}_k, U_k)$ is a $r_k$-by-$r_k$ matrix with spectral norm no more than 1, definition $\|\sin\Theta(\hat{U}_k, U_k)\|_{\rm F}^2 \leq r_k\leq r$. Therefore, one has
	\begin{equation*}
	\begin{split}
	& \mathbb{E}\|\sin\Theta(\hat{U}_k, U_k)\|_{\rm F} = \mathbb{E}\|\sin\Theta(\hat{U}_k, U_k)\|_{\rm F} 1_{Q} + \mathbb{E}\|\sin\Theta(\hat{U}_k, U_k)\|_{\rm F} 1_{Q^c}\\
	= & C\frac{\sqrt{pr}}{\lambda} + \sqrt{\mathbb{E} \|\sin\Theta(\hat{U}_k, U_k)\|_{\rm F}^2 \cdot \mathbb{E} 1_{Q^c}} \leq C\frac{\sqrt{pr}}{\lambda} + \sqrt{r\cdot C\exp(-cp)}
	\end{split}
	\end{equation*}
	By the definition of $\lambda$, we know $\lambda =\sigma_{r_k}(\mathcal{M}_k(\X)) \leq \frac{\|\X\|_{\rm F}}{\sqrt{r_k}} \leq C\frac{\exp(c_0p)}{\sqrt{r_k}}$, so we can select $c_0>0$ small enough to ensure that
	$$\frac{\sqrt{pr}}{\lambda}\geq  \frac{\sqrt{pr^2}}{C\exp(c_0p)} \geq c\sqrt{r\cdot C\exp(-cp)}.$$
	This means $\mathbb{E}\|\sin\Theta(\hat{U}_k, U_k)\|_{\rm F} \leq C\frac{\sqrt{pr}}{\lambda}.$
	Finally, for any $1\leq q\leq 2$, we have
	\begin{equation*}
	\mathbb{E}r_k^{-1/q}\left\|\sin\Theta(\hat{U}_k, U_k)\right\|_q \leq \mathbb{E}r_k^{-1/2}\left\|\sin\Theta(\hat{U}_k, U_k)\right\|_{\rm F} \leq C\frac{\sqrt{p}}{\lambda} \leq C\frac{\sqrt{p_k}}{\lambda}.
	\end{equation*}
\end{enumerate}
To sum up, we have finished the proof for this theorem.\quad $\square$

\subsection{Proof of Proposition \ref{pr:hard}}

For convenient, we introduce the following notations: $m_k = |C\cap D_k|$, $u^{(k)} = (1_{C})_{D_k}$, $U^{(k)} = u^{(k)}/\|u^{(k)}\|_2$ is the normalized vector of $u^{(k)}$, $U^{(k)}_{\perp}\in \mathbb{R}^{\frac{N}{6}\times \left(\frac{N^2}{36}-1\right)}$ is the orthogonal complement of $U^{(k)}$, $\tilde{A}_k = \mathcal{M}_k(2\cdot \A_{[D_1, D_2, D_3]} - 1_{|D_1|\times|D_2|\times|D_3|})$, for $k=1,2,3$. Without loss of generality and for convenience of the presentation, we assume $N$ is a multiple of 6. Based on the statement, 
	\begin{equation}\label{eq:tilde_A}
	(\tilde{A}_1)_{i, p_3(j-1)+k} = \left\{
	\begin{array}{ll}
	1, & \text{w.p 1}, \quad \text{if }\left(u^{(1)}_i, u^{(2)}_j, u^{(3)}_k\right)= (1, 1, 1);\\
	1, & \text{w.p 1/2}, \quad \text{if } \left(u^{(1)}_i, u^{(2)}_j, u^{(3)}_k\right) \neq (1, 1, 1);\\
	-1, & \text{w.p 1.2}, \quad \text{if } \left(u^{(1)}_i, u^{(2)}_j, u^{(3)}_k\right) \neq (1, 1, 1).
	\end{array}\right.
	\end{equation}
	$A_2$ and $A_3$ have the similar form. Therefore, $\tilde{A}_1$ are all 1 in the block of $(D_1\cap C)\times ((D_2\cap C)\otimes (D_3\cap C))$, and are with i.i.d. Rademacher entries outside the block. Since $C$ is uniformly randomly selected from $V_1$, $|V_1| = N/6, |D_k|=N/6, |C| = \kappa_N$, we know $m_1 = |D_1\cap C|, m_2 = |D_1\cap C|, m_3 = |D_1\cap C|$ satisfy hypergeometric distribution with parameter $(\kappa_N, N/2, N/6)$. Based on the concentration inequality of hypergeometric distribution (Theorem 1 in \cite{hush2005concentration}),
	\begin{equation}\label{ineq:V_1-cap-C}
	\frac{\kappa_N}{4}\leq m_k = |D_k\cap C| \leq \frac{\kappa_N}{2}, \quad k = 1,2,3
	\end{equation}
	with probability at least $1- C\exp(-c\kappa_N)$. Now the rest of the proof is similar to Theorem 3 in \cite{cai2016rate}. By \eqref{eq:tilde_A}, we have
	\begin{equation*}
	\begin{split}
	\left((U^{(1)})^\top \tilde{A}_1\right) \in \mathbb{R}^{N^2/36}, \quad \left((U^{(1)})^\top A_1\right)_{j} \left\{\begin{array}{ll}
	= \sqrt{m_1}, & j \in ((D_2\cap C)\otimes (D_3\cap C));\\
	\sim \frac{W}{\sqrt{m_1}}, & \text{otherwise},
	\end{array}
	\right.
	\end{split}
	\end{equation*}
	where $W$ has the same distribution as the sum of $m_1$ i.i.d. Rademacher random variables. Conditioning on $C$ satisfying \eqref{ineq:V_1-cap-C}, similarly as the derivation for Equation (1.15) in the Appendix of \cite{cai2016rate}, we can derive
	\begin{equation}\label{ineq:prop1-to-check1}
	\begin{split}
	\sigma_1^2\left((U^{(1)})^\top A_1\right) \geq \frac{N^2}{36} + \frac{m_1m_2m_3}{4} \quad \text{ with probability at least } 1 - C\exp(-cN);
	\end{split}
	\end{equation}
	\begin{equation}\label{ineq:prop1-to-check2}
	\begin{split}
	\sigma_2^2(A_1) \leq \frac{N^2}{36} + \frac{m_1m_2m_3}{8} \quad \text{ with probability at least } 1 - C\exp(-cN);
	\end{split}
	\end{equation}
	\begin{equation}\label{ineq:prop1-to-check3}
	\|(u^{(1)}_{\perp})^\top A_1 P_{(U^{(1)})^\top A_1} \| \leq C\sqrt{N} \quad \text{ with probability at least } 1 - C\exp(-cN).
	\end{equation}
	
	Under the circumstance that \eqref{ineq:V_1-cap-C}, \eqref{ineq:prop1-to-check1}, \eqref{ineq:prop1-to-check1}, and \eqref{ineq:prop1-to-check1} all hold, by Proposition 1 in \cite{cai2016rate}, we have
	\begin{equation*}
	\begin{split}
	\|\sin\Theta(\hat{u}_1^\top, U^{(1)})\| \leq & \frac{\sigma_2(A_1)	\|(u^{(1)}_{\perp})^\top A_1 P_{(U^{(1)})^\top A_1} \|}{\sigma_1^2\left((U^{(1)})^\top A_1\right) - \sigma_2^2(A_1)}\\
	\overset{\eqref{ineq:prop1-to-check1}\eqref{ineq:prop1-to-check1} \eqref{ineq:prop1-to-check1}}{\leq} & C\frac{\sqrt{N(N^2+m_1m_2m_3)}}{m_1m_2m_3} 
	\overset{\eqref{ineq:V_1-cap-C}}{\leq} C\frac{N^{3/2} + N^{1/2}\kappa_N^{3/2}}{\kappa_N^3}.
	\end{split}
	\end{equation*}
	Note that $\liminf_{N\to \infty} \kappa_N/\sqrt{N} = \infty$, $\lim_{N\to \infty} P( \text{\eqref{ineq:V_1-cap-C}, \eqref{ineq:prop1-to-check1}, \eqref{ineq:prop1-to-check1}, and \eqref{ineq:prop1-to-check1} all hold})=1$, we have
	$$\left\|\sin\Theta(\hat{u}_1^\top, (1_C)_{D_1})\right\| = \|\sin\Theta(\hat{u}_1^\top, U^{(1)})\|  \overset{d}{\to} 0, \quad \text{as } N\to \infty.$$ 
	The proofs for $k=2,3$ essentially follow. Therefore, we have finished the proof of this proposition. \quad $\square$
	
\section{Appendix: Technical Lemmas}

We collect all technical lemmas that has been used in the theoretical proofs throughout the paper in this section. 

The following lemma shows the equivalence between two widely considered Schatten $q$-norm distances for singular subspaces.
	\begin{Lemma}\label{lemma:norms}
		For any $U_1, U_2\in \mathbb{O}_{p,r}$ and all $1\leq q\leq +\infty$,
		$$
		\frac{1}{4}\|U_1U_1^\top-U_2U_2^\top\|_q\leq \|\sin\Theta(U_1,U_2)\|_q\leq \|U_1U_1^\top-U_2U_2^\top\|_q.
		$$
	\end{Lemma}

We will use the following properties of tensor algebra in the technical analysis of this paper.
\begin{Lemma}[Properties in Tensor Algebra]\label{lm:tensor-algebra}~
	\begin{itemize}
	\item Suppose $\X\in \mathbb{R}^{p_1\times p_2\times p_3}$, $U_k\in \mathbb{R}^{p_k\times r_k}$ for $k=1,2,3$. Then we have the following identity related to tensor matricizations,
	\begin{equation}\label{eq:tensor-algebra-1}
	\mathcal{M}_k\left(\X\times_{k+1} U_{k+1}^\top \times_{k+2} U_{k+2}^\top \right) = \mathcal{M}_k(\X) (U_{k+1}\otimes U_{k+2}),\quad k=1,2,3.
	\end{equation}
	\item Suppose we further have $\tilde{U}_k \in \mathbb{R}^{p_k\times \tilde{r}_k}$ for $k=1,2,3$, then
	\begin{equation}\label{eq:tensor-algebra-2}
	\left(\tilde{U}_1\otimes \tilde{U}_2\right)^\top = (\tilde{U}_1^\top) \otimes (\tilde{U}_2^\top),\quad \left(\tilde{U}_2\otimes \tilde{U}_3\right)^\top \left(U_2\otimes U_3\right) = \left(\tilde{U}_2^\top U_2\right)\otimes \left(\tilde{U}_3^\top U_3\right).
	\end{equation}
	\begin{equation}\label{eq:tensor-algebra-3}
	\begin{split}
	& \|U_2\otimes U_3\| = \|U_2\|\cdot \|U_3\|, \quad \|U_2\otimes U_3\|_{\rm F} = \|U_2\|_{\rm F}\cdot \|U_3\|_{\rm F},\\
	& \sigma_{\min}(U_2\otimes U_3) = \sigma_{\min}(U_2)\sigma_{\min}(U_3).
	\end{split}
	\end{equation}
	
	\item (Properties related to projections) Suppose $U_2\in \mathbb{O}_{p_2, r_2}, U_3\in \mathbb{O}_{p_3, r_3}$, and $U_{2\perp}\in \mathbb{O}_{p_2, p_2-r_2}, U_{3\perp}\in \mathbb{O}_{p_3, p_3-r_3}$ are their orthogonal complement, respectively. Then $P_{U_2\otimes U_3} = P_{U_2}\otimes P_{U_3}$, and we have the following decomposition
	\begin{equation}
	\begin{split}
	I_{p_2p_3} = & P_{I_{p_2}\otimes U_3} + P_{I_{p_3}\otimes {U_{3\perp}}} = P_{U_2\otimes I_{p_3}} + P_{U_{3\perp} \otimes I_{p_2}}\\
	= & P_{U_2\otimes U_3} + P_{U_{2\perp}\otimes U_3} + P_{U_2\otimes U_{3\perp}} + P_{U_{2\perp}\otimes U_{3\perp}}.
	\end{split}
	\end{equation}
	\end{itemize}
\end{Lemma}

The following lemma characterizes the maximum of norms for i.i.d. Gaussian tensors after any projections. 
\begin{Lemma}\label{lm:random_tensor_max_projection}
	For i.i.d. Gaussian tensor $\Z \in \mathbb{R}^{p_1 \times p_2 \times p_3}$, $\Z\overset{iid}{\sim} N(0, 1)$, we have the following tail bound for the projections,
	\begin{equation*}
	\begin{split}
	& P\left(\max_{\substack{V_2 \in \mathbb{R}^{p_2\times r_2}, V_3\in \mathbb{R}^{p_3\times r_3}\\\|V_2\| \leq 1, \|V_3\|\leq 1 }} \left\|\mathcal{M}_1\left(\Z \times_2 V_2^\top \times_3 V_3^\top\right)\right\|  \geq C\sqrt{p_1} + C\sqrt{r_2r_3} + C\sqrt{1+t}\left(\sqrt{p_2r_2} + \sqrt{p_3r_3}\right) \right)\\ 
	\leq & C\exp(-Ct(p_2r_2+p_3r_3))
	\end{split}
	\end{equation*}
	for any $t > 0$. Similar results also hold for $\mathcal{M}_2\left(\Z \times_1 V_1^\top \times_3 V_3^\top\right)$ and $\mathcal{M}_3\left(\Z \times_1 V_1^\top \times_2 V_2^\top \right)$.
	Meanwhile, there exists uniform $C>0$ such that
	\begin{equation}\label{ineq:random_tensor_target2}
	\begin{split}
	&P\left(\max_{\substack{V_1, V_2, V_3 \in \mathbb{R}^{p\times r}\\\max\{\|V_1\|, \|V_2\|, \|V_3\|\} \leq 1}} \left\|\Z \times_1 V_1^\top \times_2 V_2^\top \times_3 V_3^\top \right\|_{\rm F}^2 \geq Cr_1r_2r_3 + C(1+t)(p_1r_1+p_2r_2+p_3r_3)\right)\\
	\leq & \exp\left(-Ct(p_1r_1+p_2r_2+p_3r_3)\right)
	\end{split}
	\end{equation}
	for any $t > 0$.
\end{Lemma}

In the perturbation bound analysis in this paper, we also need the following technical result to bound the spectral and Frobenius norm for the projections.
\begin{Lemma}\label{lm:r-pinciple-compoenent}
	Suppose $X, Z \in \mathbb{R}^{p_1\times p_2}$, $\rank(X) = r$. If the singular value decomposition of $X$ and $Y$ are written as 
	$$Y = X+Z = \hat{U}\hat{\Sigma} \hat{V}^\top = \begin{bmatrix}
	\hat{U}_1 & \hat{U}_2
	\end{bmatrix} \cdot \begin{bmatrix}
	\hat{\Sigma}_1 & \\
	 & \hat{\Sigma}_2
	\end{bmatrix}\cdot \begin{bmatrix}
	\hat{V}_1^\top ~ & ~ \hat{V}_2^\top
	\end{bmatrix},$$
	where $\hat{U}_1 \in \mathbb{O}_{p_1, r}, \hat{V}_1\in \mathbb{O}_{p_2,r}$ correspond to the leading $r$ left and right singular vectors; and $\hat{U}_2 \in \mathbb{O}_{p_1, p_2-r}, \hat{V}_1\in \mathbb{O}_{p_2, p_2-r}$ correspond to their orthonormal complement. Then
	\begin{equation*}
	\left\|P_{\hat{U}_2}X\right\| \leq 2\|Z\|,\quad \left\|P_{\hat{U}_2}X\right\|_{\rm F} \leq \min\left\{2\sqrt{r}\|Z\|, 2\|Z\|_{\rm F}\right\}.
	\end{equation*}
\end{Lemma}

The following lemma provides a detailed analysis for $\varepsilon$-net for the class of regular matrices under various norms and for the low-rank matrices under spectral norm.
\begin{Lemma}[$\varepsilon$-net for Regular and Low-rank Matrices]\label{lm:covering_set_spectral} ~
	\begin{itemize}
		\item Suppose $\|\cdot\|_\bullet$ is any matrix norm, $\mathcal{X}_{p_1, p_2} = \{X\in \mathbb{R}^{p_1\times p_2}: \|X\|\leq 1\}$ is the unit ball around the center in $\|\cdot\|_\bullet$ norm. Then there exists an $\varepsilon$-net $\bar{\mathcal{X}}_{p_1, p_2}$ in $\|\cdot\|_\bullet$ norm with cardinality at most $((2+\varepsilon)/\varepsilon)^{p_1p_2}$ for $\mathcal{X}_{p_1, p_2}$. To be specific, there exists $X^{(1)},\ldots, X^{(N)}$ with $N \leq ((2+\varepsilon)/\varepsilon)^{p_1p_2}$, such that for all $X\in \mathcal{X}_{p_1, p_2}$, there exists $i\in \{1,\ldots, N\}$ satisfying $\|X^{(i)} - X\|\leq \varepsilon$.
		
		\item Let $\mathcal{X}_{p_1, p_2, r} = \{X\in \mathbb{R}^{p_1\times p_2}: \rank(X) \leq r, \|X\|\leq 1\}$ be the class of low-rank matrices under spectral norm. Then there exists an $\varepsilon$-net $\bar{\mathcal{X}}_r$ for $\mathcal{X}_{p_1, p_2, r}$ with cardinality at most $((4+\varepsilon)/\varepsilon)^{(p_1+p_2)r}$. Specifically, there exists $X^{(1)},\ldots, X^{(N)}$ with $N \leq ((4+\varepsilon)/\varepsilon)^{(p_1+p_2)r}$, such that for all $X\in \mathcal{X}_{p_1, p_2, r}$, there exists $i\in \{1,\ldots, N\}$ satisfying $\|X^{(i)} - X\|\leq \varepsilon$.
	\end{itemize}
\end{Lemma}

The next lemma characterizes the tail probability for i.i.d. Gaussian vector after multiplication of any fixed matrix.
\begin{Lemma}\label{lm:gaussian-vector-projection}
	Suppose $u\in \mathbb{R}^{p}$ such that $u \overset{iid}{\sim}N(0, 1)$, $A\in \mathbb{R}^{n\times p}$ is a fixed matrix. Then,
	\begin{equation*}
	\begin{split}
	& P\left(\|Au\|_2^2 - \|A\|_{\rm F}^2 \leq -2\|A^\top A\|_{\rm F}\sqrt{t}\right) \leq \exp(-t);\\ 
	& P\left(\|Au\|_2^2 - \|A\|_{\rm F}^2 \geq 2\|A^\top A\|_{\rm F}\sqrt{t} + 2\|A\|^2t\right) \leq \exp(-t).
	\end{split}
	\end{equation*}
\end{Lemma}

\section{Proof of Technical Lemmas}

\subsection{Proof of Lemma~\ref{lemma:detection}}
	Without loss of generality, assume that $p \equiv 0 ({\rm mod}\ 2)$.
	Hereafter, set $N=3p$ and $\kappa_N=20k$ with $k=\floor{p^{(1-\tau)/2}}$.
	Our main technique is based on a reduction scheme which maps any adjacency tensor $\A\in\{0,1\}^{N\times N\times N}$ to a random tensor $\Y\in\mathbb{R}^{p\times p\times p}$ in $O(N^3)$ number of flops. The technique was invented in \cite{ma2015computational}, adapted from a bottom-left trick in \cite{berthet2013computational}. Some other related methods can be found in \cite{cai2015computational} and \cite{wang2014statistical}. For the completeness and readability of our paper, we provide a detailed application of this technique to the tensor settings.
	
	To this end, for any $M\geq 3$ and $0<\mu\leq \frac{1}{2M}$, define two random variables
	$$
	\xi^+:=(Z+\mu){\bf 1}(|Z|\leq M)\quad {\rm and}\quad \xi^-:=(\tilde{Z}-\mu){\bf 1}(|\tilde{Z}|\leq M)
	$$
	where $Z$ and $\tilde{Z}$ denote independent standard normal random variables. The randomized mapping from $\A\in\{0,1\}^{N\times N\times N}$ to a random matrix $\Y\in\mathbb{R}^{p\times p\times p}$ is essentially one step of Gaussianization. For simplicity, denote $V_1:=\{1,2,\ldots, \frac{p}{2}\}\cup \{\frac{3p}{2}+1,\ldots,2p\}$,
	\begin{align*}
	V_2:=\big\{\frac{p}{2}+1,\ldots,p\big\}\cup \big\{2p+1,\ldots,\frac{5p}{2}\big\},
	\end{align*}
	and
	$$
	V_3:=\big\{p+1,\ldots,\frac{3p}{2}\big\}\cup \big\{\frac{5p}{2}+1,\ldots,3p\big\}.
	$$
	Therefore, $V_1, V_2, V_3$ are disjoint and $V_1\cup V_2\cup V_3=[N]$. Given an adjacency tensor $\A\in\{0,1\}^{N\times N\times N}$, let $\A_0=\A_{V_1,V_2,V_3}\in\mathbb{R}^{p\times p\times p}$ be a corner block of $\A$. Conditioned on $\A_0$, we generate a random tensor $\Y\in\mathbb{R}^{p\times p\times p}$ such that
	$$
	Y_{a,b,c}=\big(1-(A_0)_{a,b,c}\big){\Xi}^-_{a,b,c}+(A_0)_{a,b,c}{\Xi}^+_{a,b,c},\quad \forall a,b,c\in[p]
	$$
	where ${\bf \Xi}^-\in\mathbb{R}^{p\times p\times p}$ has i.i.d. entries with the same distribution as $\xi^-$ and ${\bf \Xi}^+\in\mathbb{R}^{p\times p\times p}$ has i.i.d. entries with the same distribution as $\xi^+$.
	Clearly, this process defines a deterministic map for any fixed ${\bf \Xi}^-, {\bf \Xi}^+\in\mathbb{R}^{p\times p\times p}$
	\begin{align*}
	\mathcal{T}: \{0,1\}^{N\times N\times N}\times \mathbb{R}^{p\times p\times p}\times\mathbb{R}^{p\times p\times p}&\mapsto \mathbb{R}^{p\times p\times p}\\
	(\A,{\bf \Xi}^-,{\bf \Xi}^+)\mapsto \Y.
	\end{align*}
	Let $\mathcal{L}(\X)$ denote the law of a random tensor $\X$. The total variation distance between two probability distributions $\mathbb{P}_1$ and $\mathbb{P}_2$ is denoted by ${\rm d_{TV}}(\mathbb{P}_1,\mathbb{P}_2)$. The following lemma is analogous to  \cite[Lemma~2]{ma2015computational} and the proof is skipped here.
	\begin{Lemma}\label{f12lemma}
		Let $M\geq 4$, $\mu\leq \frac{1}{2M}$, and $\eta$ be a Bernoulli random variable. Suppose $\xi$ is a random variable such that $(\xi|\eta=1)=\xi^+$ and $(\xi|\eta=0)=\xi^-$. 
		\begin{enumerate}[label=(\arabic*)]
			\item If $\mathbb{P}(\eta=1)=1$, then ${\rm d_{TV}}\big(\mathcal{L}(\xi),\mathcal{N}(\mu,1)\big)\leq e^{(1-M^2)/2}$;
			\item If $\mathbb{P}(\eta=0)=\mathbb{P}(\eta=1)=\frac{1}{2}$, then ${\rm d_{TV}}(\mathcal{L}(\xi),\mathcal{N}(0,1))\leq e^{-M^2/2}$.
		\end{enumerate}
	\end{Lemma}
	Our next step is to show that (by choosing $M=\sqrt{8\log 3p}$ and $\mu=(2M)^{-1}$), the law of $\Y=\mathcal{T}(\A, {\bf \Xi}^-, {\bf \Xi}^+)$ is asymptotically equivalent to a mixture over $\{\mathbb{P}_{\X}: \X\in\mathcal{M}_0(\bp,k,\br,\lambda)\}$ for $\lambda=\frac{p^{3(1-\tau)/4}}{2\sqrt{8\log 3p}}$ if $G\sim H_0$. On the other hand, if $G\sim H_1$, the law of $\Y=\mathcal{T}(\A, {\bf \Xi}^-, {\bf \Xi}^+)$ is asymptotically equivalent to a mixture over $\{\mathbb{P}_{\X}: \X\in\mathcal{M}_1(\bp,k,\br,\lambda)\}$. For an adjacency tensor $\A\in\{0,1\}^{N\times N\times N}$, we have $\mathcal{T}(\A,{\bf \Xi}^-, {\bf \Xi}^+)\in\mathbb{R}^{p\times p\times p}$. Recall that $\Y\in\mathbb{R}^{p_1\times p_2\times p_3}$ and we define an embedding $\ell: \mathbb{R}^{p\times p\times p}\mapsto \mathbb{R}^{p_1\times p_2\times p_3}$,
	\begin{equation}
	\ell(\A)_{ijk}=
	\begin{cases}
	A_{ijk}&{\rm if}\ (i,j,k)\in [p]\times [p]\times [p]; \\
	0& {\rm otherwise}.
	\end{cases}
	\end{equation}
	Lemma~\ref{lemma:tv} is similar to \cite[Lemma~4]{ma2015computational}. 
	We postpone the proof of Lemma~\ref{lemma:tv} to the Appendix.
	\begin{Lemma}\label{lemma:tv}
		Let $\A\in\mathbb{R}^{N\times N\times N}$ be the adjacency tensor of a hypergraph $G$ sampled from either $H_0$ or $H_1$ and $\Y=\ell\circ\mathcal{T}(\A, {\bf \Xi}^-, {\bf \Xi}^+)$. Suppose that $M=\sqrt{8\log N}$ and $\mu=\frac{1}{2M}$. For each $i=0,1$, if $G\sim H_i$, there exists a prior distribution $\pi_i$ on $\mathcal{M}_i(\bp,k,\br,\lambda)$ with $\lambda=\frac{p^{3/4(1-\tau)}}{2\sqrt{8\log 3p}}$ such that 
		$$
		{\rm d_{TV}}\big(\mathcal{L}(\Y),\mathbb{P}_{\pi_i}\big)\leq\frac{\sqrt{e}}{27N}+6k\big(0.86\big)^{2.5k}.
		$$
		where $\mathbb{P}_{\pi_i}=\int_{\mathcal{M}_i(\bp,k,\br,\lambda)}\mathbb{P}_\X(\cdot)\pi_i(d\X)$.
	\end{Lemma}
	Now, on the contradictory, suppose that the claim of Lemma~\ref{lemma:detection} does not hold. It means that there exists a sequence of polynomial-time tests $\{\phi_{p_t}\}$ with a sub-sequence $(p_t)_{t=1}^{\infty}$ of positive integers such that 
	$$
	\lim_{t\to\infty}\mathcal{R}_{\bp,\br,\lambda}(\phi_{p_t})=\lim_{t\to\infty}\Big\{\underset{\X\in\mathcal{M}_0(\bp,k,\br,\lambda)}{\sup}\mathbb{P}_{\X}\big\{\phi_{p_t}(\Y)=1\big\}+\underset{\X\in\mathcal{M}_1(\bp,k,\br,\lambda)}{\sup}\mathbb{P}_\X\big\{\phi_{p_t}(\Y)=0\big\}\Big\}<\frac{1}{2}.
	$$
	Define the test $\psi_{N_t}(\A)=\phi_{p_t}\big(\ell\circ \mathcal{T}(\A,{\bf \Xi}^-, {\bf \Xi}^+)\big)$ and we obtain a sequence of polynomial-time tests $\{\psi_{N_t}\}$ for problem (\ref{eq:hypergraph}) with $N_t=3p_t$ for $t=1,\ldots,\infty$. It suffices to compute 
	$$
	\mathcal{R}_{N_t,\kappa_{N_t}}(\psi_{N_t})=\mathbb{P}_{H_0}\big\{\psi_{N_t}(\A)=1\big\}+\mathbb{P}_{H_1}\big\{\psi_{N_t}(\A)=0\big\}
	$$
	with $\kappa_{N_t}=20\floor{p_t^{(1-\tau)/2}}$. Note that $\lim_{t\to\infty}\frac{\log \kappa_{N_t}}{\log\sqrt{N_t}}\leq 1-\tau$. By definition of ${\rm d_{TV}}$ and Lemma~\ref{lemma:tv}, under $H_0$,
	\begin{align*}
	\Big|\mathbb{P}_{H_0}\big\{\psi_{N_t}(\A)=1\big\}&-\mathbb{P}_{\pi_0}\big\{\phi_{p_t}(\Y)=1\big\}\Big|\\
	=&\Big|\mathbb{P}_{H_0}\big\{\phi_{p_t}\big(\ell\circ\mathcal{T}(\A, {\bf \Xi}^-, {\bf \Xi}^+)\big)=1\big\}-\mathbb{P}_{\pi_0}\big\{\phi_{p_t}(\Y)=1\big\}\Big|\\
	\leq&{\rm d_{TV}}\big(\mathcal{L}\big(\ell\circ\mathcal{T}(\A,{\bf \Xi}^-, {\bf \Xi}^+)\big),\mathbb{P}_{\pi_0}\big)\leq \frac{\sqrt{e}}{27N_t}+6k_t(0.86)^{2.5k_t}
	\end{align*}
	where $k_t=\floor{p_t^{(1-\tau)/2}}$ and we used the fact that the mixture $\mathbb{P}_{\pi_0}$ over $\mathcal{M}_0(\bp,k,\br,\lambda)$ is also a mixture over $\mathcal{F}_{\bp,\br}(\lambda)$. In a similar fashion,
	\begin{align*}
	\Big|\mathbb{P}_{H_1}\big\{\psi_{N_t}(\A)=0\big\}&-\mathbb{P}_{\pi_1}\big\{\phi_{p_t}(\Y)=0\big\}\Big|\leq
	\frac{\sqrt{e}}{27N_t}+6{k_t}(0.86)^{2.5k_t}
	\end{align*}
	As a result, 
	\begin{eqnarray*}
		\mathcal{R}_{N_t,\kappa_{N_t}}(\psi_{N_t})=\mathbb{P}_{H_0}\big\{\psi_{N_t}(\A)=1\big\}+\mathbb{P}_{H_1}\big\{\psi_{N_t}(\A)=0\big\}\\
		\leq\mathbb{P}_{\pi_0}\big\{\phi_{p_t}(\Y)=1\big\}+\mathbb{P}_{\pi_1}\big\{\phi_{p_t}(\Y)=0\big\}+\frac{2\sqrt{e}}{27N_t}+12k_t(0.86)^{2.5k_t}\\
		\leq\underset{\X\in\mathcal{M}_0(\bp,k,\br,\lambda)}{\sup}\mathbb{P}_{\X}\big\{\phi_{p_t}(\Y)=1\big\}+\underset{\X\in\mathcal{M}_1(\bp,k,\br,\lambda)}{\sup}\mathbb{P}_\X\big\{\phi_{p_t}(\Y)=0\big\}+\frac{2\sqrt{e}}{27N_t}+12k_t(0.86)^{2.5k_t}
	\end{eqnarray*}
	Therefore,
	$$
	\underset{t\to \infty}{\lim}\ \mathcal{R}_{N_t,\kappa_{N_t}}(\psi_{N_t})<\frac{1}{2},
	$$
	which contradicts the hypothesis {\bf H($\tau$)}. \quad $\square$

\subsection{Proof of Lemma~\ref{lemma:norms}} Let $\sigma_1\geq \sigma_2\geq \ldots\geq \sigma_r$ denote the singular values of $U_1^\top U_2$. It is easy to check that the singular values of $U_{1\perp}^\top U_2$ are $\sqrt{1-\sigma_r^2}\geq \ldots\geq \sqrt{1-\sigma_{2}^2}\geq \sqrt{1-\sigma_1^2}$, in view of the fact
$$
(U_1^\top U_2)^\top(U_1^\top U_2)+(U_{1\perp}^\top U_2)^\top(U_{1\perp}^\top U_2)=U_2^\top U_2=I_r.
$$
Recall that for all $1\leq q\leq +\infty$,
$$
\|\sin\Theta(U_1,U_2)\|_q=\Big(\sum_{i=1}^r\big(\sin(\cos^{-1}\sigma_i)\big)^q\Big)^{1/q}=\Big(\sum_{i=1}^r\big(1-\sigma_i^2\big)^{q/2}\Big)^{1/q}.
$$
The following fact is straightforward:
\begin{align*}
\|U_2U_2^\top-U_1U_1^\top \|_q\geq \|U_{1\perp}^\top U_2U_2^\top\|_q=\|U_{1\perp}^\top U_2\|_{q}=\Big(\sum_{i=1}^r\big(1-\sigma_i^2\big)^{q/2}\Big)^{1/q}
\end{align*}
which concludes $\|U_2U_2^\top-U_1U_1^\top \|_q\geq \|\sin\Theta(U_1,U_2)\|_q$. On the other hand,
\begin{align*}
\|U_2U_2^\top-U_1U_1^\top \|_q\leq& \|\Proj_{U_1}(U_2U_2^\top-U_1U_1^\top)\Proj_{U_1}\|_q+\|\Proj_{U_1}(U_2U_2^\top)\Proj_{U_1}^\perp\|_q\\
+&\|\Proj_{U_1}^\perp (U_2U_2^\top)\Proj_{U_1}\|_q+\|\Proj_{U_1}^\perp(U_2U_2^\top)\Proj_{U_1}^\perp\|_q\\
\leq& \big\|U_1\big(U_1^\top U_2 U_2^\top U_1-I_r\big)U_1^\top\big\|_q+\|U_2^\top U_{1\perp }\|_q+\|U_{1\perp}^\top U_2\|_q+\|U_{1\perp}^\top U_2 U_2^\top U_{1\perp }\|_q\\
\leq& 4\Big(\sum_{i=1}^r\big(1-\sigma_i^2\big)^{q/2}\Big)^{1/q}\leq 4\|\sin\Theta(U_1,U_2)\|_q
\end{align*}
where we used the fact $1-\sigma_{i}^2\leq \sqrt{1-\sigma_i^2}$ for all $1\leq i\leq r$. \quad $\square$

\subsection{Proof of Lemma \ref{lm:tensor-algebra}} 
\begin{itemize}
	\item First, we shall note that both $\mathcal{M}_k\left(\X\times_{k+1} U_{k+1}^\top \times_{k+2} U_{k+2}^\top \right)$ and $\mathcal{M}_k(\X) (U_{k+1}\otimes U_{k+2})$ are of dimension $p_k$-by-$(r_{k+1}r_{k+2})$. To prove they are equal, we just need to compare each of their entries. We focus on $k=1$ as the $k=2,3$ essentially follows. For any $1\leq i_1 \leq p_1, 1\leq i_2 \leq r_2, 1\leq i_3 \leq r_3$, one has
	\begin{equation*}
	\begin{split}
	& \left[\mathcal{M}_1\left(\X\times_{2} U_{2}^\top \times_{3} U_{3}^\top \right)\right]_{i_1, (i_2-1)r_3 + i_3} = \left(\X\times_{2} U_{2}^\top \times_{3} U_{3}^\top \right)_{i_1,i_2,i_3}\\
	= & \sum_{j_2=1}^{p_2}\sum_{j_3=1}^{p_3} X_{i_1, j_2, j_3} (U_2)_{j_2, i_2} (U_3)_{j_3, i_3} \\
	= & \sum_{j_2=1}^{p_2}\sum_{j_3=1}^{p_3}\left(\mathcal{M}_1(\X)\right)_{i_1, (j_2-1)p_3+j_3} \cdot \left(U_2\otimes U_3\right)_{(j_2-1)p_3+j_3, (i_2-1)r_3+i_3}\\
	= & \left(\mathcal{M}_1(\X) \cdot (U_2\otimes U_3)\right)_{i_1, (i_2-1)r_3+i_3}.
	\end{split}
	\end{equation*}
	This shows \eqref{eq:tensor-algebra-1}. 
	\item The proof for \eqref{eq:tensor-algebra-2} is essentially the same as \eqref{eq:tensor-algebra-1} as we only need to check each entries of the terms in \eqref{eq:tensor-algebra-2} are equal. For \eqref{eq:tensor-algebra-3}, let 
	$$U_2 = \sum_{i} \sigma_{2i}\cdot \alpha_{2i}  \beta_{2i}^\top,\quad U_3 = \sum_{j} \sigma_{3j}\cdot \alpha_{3j} \beta_{3j}^\top $$
	be the singular value decompositions. Then it is not hard to see the singular value decomposition of $U_2\otimes U_3$ can be written as
	$$U_2\otimes U_3 = \sum_{i, j} \sigma_{2i}\sigma_{3j} \cdot (\alpha_{i2}\otimes\alpha_{j3}) (\beta_{i2}\otimes\beta_{j3})^\top, $$
	so that the singular values of $U_\otimes U_3$ are $\{\sigma_i\cdot \sigma_j\}$. Then 
	$$\|U_2\otimes U_3\| = \max_{i,j}\sigma_{2i}\sigma_{3j} = \left(\max_i \sigma_{2i}\right)\cdot \left(\max_j \sigma_{3j}\right) = \|U_2\|\cdot \|U_3\|,$$
	$$\|U_2\otimes U_3\|_{\rm F}^2 = \sum_{i,j}\sigma_{2i}^2\sigma_{3j}^2 = \left(\sum_i \sigma_{2i}^2\right)\cdot \left(\sum_j \sigma_{3j}^2\right) = \|U_2\|_{\rm F}^2\cdot \|U_3\|_{\rm F}^2,$$
	and
	$$\sigma_{\min}(U_2\otimes U_3) = \min_{i,j}\sigma_{2i}\sigma_{3j} = \left(\min_i \sigma_{2i}\right)\cdot \left(\min_j \sigma_{3j}\right) = \sigma_{\min}(U_2)\cdot \sigma_{\min}(U_3). $$
	\item 
	\begin{equation*}
	\begin{split}
	& P_{I_{p_2}\otimes U_3} + P_{I_{p_2}\otimes U_{3\perp}} = (I_{p_2}\otimes U_3) (I_{p_2}\otimes U_3)^\top + (I_{p_2}\otimes U_{3\perp})(I_{p_2}\otimes U_{3\perp})^\top\\
	\overset{\eqref{eq:tensor-algebra-2}}{=} & (I_{p_2}I_{p_2}^\top)\otimes (U_3U_3^\top) + (I_{p_2}I_{p_2}^\top)\otimes (U_{3\perp}U_{3\perp}^\top) = I_{p_2}\otimes (U_3U_3^\top + U_{3\perp}U_{3\perp}^\top) \\
	= & I_{p_2}\otimes I_{p_3} = I_{p_2p_3}.
	\end{split}
	\end{equation*}
	The other identity can be shown similarly. \quad $\square$
\end{itemize}

\subsection{Proof of Lemma \ref{lm:random_tensor_max_projection}} The key idea for the proof of this lemma is via $\varepsilon$-net. By Lemma \ref{lm:covering_set_spectral}, for $k=1,2,3$, there exist $\varepsilon$-nets: $V_k^{(1)}, \ldots, V_k^{(N_k)}$ for $\{V_k\in \mathbb{R}^{p_k\times r_k}: \|V_k\|\leq 1\}$, $|N_k| \leq ((4+\varepsilon)/\varepsilon)^{p_kr_k}$, such that 
$$\text{For any } V\in \mathbb{R}^{p_k\times r_k} \text{ satisfying }\|V\|\leq 1, \text{there exists $V_k^{(j)}$ such that } \|V_k^{(j)} - V\| \leq \varepsilon.$$ 
For fixed $V_2^{(i)}$ and $V_3^{(j)}$, we consider
$$Z^{(ij)}_1 = \mathcal{M}_1\left(\Z\times_2 (V_2^{(i)})^\top\times_3 (V_3^{(j)})^\top\right) \in \mathbb{R}^{p_1\times (r_2r_3)}.$$
Clearly, each row of $Z^{(ij)}_1$ follows a joint Gaussian distribution: $N\left(0, \left(V_2^{(i)\top} V_2^{(i)}\right)\otimes \left(V_3^{(j)\top} V_3^{(j)}\right) \right)$, and $\left\|\left(V_2^{(i)\top} V_2^{(i)}\right)\otimes \left(V_3^{(j)\top} V_3^{(j)}\right)\right\| \leq 1$. Then by random matrix theory (e.g. \cite{vershynin2010introduction}),
\begin{equation*}
P\left(\|Z_1^{(ij)}\| \leq \sqrt{p_1} + \sqrt{r_2r_3} + t\right) \geq 1 - 2\exp(-t^2/2).
\end{equation*}
Then we further have
\begin{equation}\label{ineq:Z^ij}
P\left(\max_{i,j} \|Z_1^{(ij)}\| \leq \sqrt{p_1} + \sqrt{r_2r_3} + x\right) \geq 1 - 2((4+\varepsilon)/\varepsilon)^{p_2r_2 + p_3r_3}\exp(-x^2/2),
\end{equation}
for all $x>0$. Now, we assume
\begin{equation*}
\begin{split}
V_2^\ast , V_3^\ast = & \argmax_{\substack{V_2\in \mathbb{R}^{p_2\times r_2}, V_3\in \mathbb{R}^{p_3\times r_3}\\\|V_2\|\leq 1, \|V_3\|\leq 1}} \left\|\mathcal{M}_1\left(\Z \times_2 V_2^\top \times_3 V_3^\top\right)\right\|,\\
M = & \max_{\substack{V_2\in \mathbb{R}^{p_2\times r_2}, V_3\in \mathbb{R}^{p_3\times r_3}\\\|V_2\|\leq 1, \|V_3\|\leq 1}} \left\|\mathcal{M}_1\left(\Z \times_2 V_2^\top \times_3 V_3^\top\right)\right\|. 
\end{split}
\end{equation*}
By definition of the $\varepsilon$-net, we can find $1\leq i \leq N_2$ and $1\leq j\leq N_3$ such that $\|V_2^{(i)} - V_2^\ast\|\leq \varepsilon$ and $\|V_3^{(i)} - V_3^\ast\|\leq \varepsilon$. In this case under \eqref{ineq:Z^ij},
\begin{equation*}
\begin{split}
M = & \left\|\mathcal{M}_1\left(\Z\times_2 (V_2^\ast)^\top \times_3 (V_3^\ast)^\top\right)\right\|\\
\leq & \left\|\mathcal{M}_1\left(\Z\times_2 (V_2^{(i)})^\top \times_3 (V_3^{(j)})^\top\right)\right\| + \left\|\mathcal{M}_1\left(\Z\times_2 (V^\ast - V_2^{(i)})^\top \times_3 (V_3^{(j)})^\top\right)\right\|\\
& + \left\|\mathcal{M}_1\left(\Z\times_2 (V_2^\ast)^\top \times_3 (V_3^\ast - V_3^{(j)})^\top\right)\right\|\\
\leq & \sqrt{p_1} + \sqrt{r_2r_3} + x + \varepsilon M + \varepsilon M,
\end{split}
\end{equation*}
Therefore, we have
\begin{equation*}
\begin{split}
P\left(M \leq \frac{\sqrt{p_1} + \sqrt{r_2r_3}+x}{1-2\varepsilon}\right) \geq 1 - 2((4+\varepsilon)/\varepsilon)^{p_2r_2 + p_3r_3}\exp(-x^2/2).
\end{split}
\end{equation*}
By setting $\varepsilon = 1/3$, $x^2 = 2\log(13)(p_2r_2 + p_3r_3)(1+t)$, we have proved the first part of the lemma. 

The proof for the second part is similar. For any given $V_k\in \mathbb{R}^{p_k\times r_k}$ satisfying $\|V_1\|, \|V_2\|, \|V_3\|\leq 1$, we have
$\|V_1\otimes V_2\otimes V_3\| \leq 1$. By Lemma \ref{lm:gaussian-vector-projection}, we know 
\begin{equation*}
\begin{split}
P\Big( & \left\|\Z\times_1 V_1^\top \times_2 V_2^\top \times_3 V_3^\top\right\|_{\rm F}^2 - \|V_1\otimes V_2\otimes V_3\|_{\rm F}^2 \\
& \geq 2\sqrt{t\|(V_1^\top V_1) \otimes (V_2^\top V_2) \otimes (V_3^\top V_3)} + 2t \|V_1\otimes V_2\otimes V_3\|^2 \Big) \leq \exp(-t).
\end{split}
\end{equation*}
Since $\|V_1\otimes V_2\otimes V_3\| \leq 1$, $\|V_1\otimes V_2\otimes V_3\|_{\rm F}^2 = \|V_1\|_{\rm F}^2\|V_2\|_{\rm F}^2\|V_3\|_{\rm F}^2 \leq r_1r_2r_3$, then
\begin{equation*}
\begin{split}
& \|(V_1^\top V_1)\otimes (V_2^\top V_2)\otimes (V_3^\top V_3)\|_{\rm F}^2 = \|V_1^\top V_1\|_{\rm F}^2\|V_2^\top V_2\|_{\rm F}^2\|V_3^\top V_3\|_{\rm F}^2\\
= & \left(\sum_{i=1}^{r_1} \sigma_i^4(V_1)\right)\left(\sum_{i=1}^{r_2} \sigma_i^4(V_1)\right)\left(\sum_{i=1}^{r_3} \sigma_i^4(V_1)\right) \leq r_1r_2r_3,
\end{split}
\end{equation*} 
we have for any fixed $V_1, V_2, V_3$ and $x>0$ that
\begin{equation*}
P\left(\left\|\Z\times_1 V_1^\top \times_2 V_2^\top \times_3 V_3^\top\right\|_{\rm F}^2 \geq r_1r_2r_3 + 2\sqrt{r_1r_2r_3 x} + 2x\right) \leq \exp(-x).
\end{equation*}
By geometric inequality, $2\sqrt{r_1 r_2r_3x} \leq r_1r_2r_3 + x$, then we further have
\begin{equation*}
P\left(\left\|\Z\times_1 V_1^\top \times_2 V_2^\top \times_3 V_3^\top\right\|_{\rm F}^2 \geq 2r_1r_2r_3 + 3x\right) \leq \exp(-x).
\end{equation*}
The rest proof for this lemma is similar to the first part. By Lemma \ref{lm:covering_set_spectral}, one can find three $\varepsilon$-nets: $V_k^{(1)}, \ldots, V_k^{(N_k)}$ for $\{V_k\in \mathbb{R}^{p_k\times r_k}: \|V_k\|\leq 1\}$ such that $|N_k|\leq ((4+2\varepsilon)/\varepsilon)^{p_kr_k}$, $k=1,2,3$. Then by probability union bound,
\begin{equation}\label{ineq:jointly_V1_V2_V3}
\begin{split}
& \max_{V_1^{(a)}, V_2^{(b)}, V_3^{(c)}}P\left(\left\|\Z\times_1V_1^\top \times_2 V_2^\top \times_3 V_3^\top \right\|_{\rm F}^2 \geq 2r_1r_2r_3 + 3x\right) \\
\leq & \exp(-x)\cdot ((4+\varepsilon)/\varepsilon)^{p_1r_1 + p_2r_2 + p_3r_3}.
\end{split}
\end{equation}
When the inequality above holds, we suppose 
$$(V_1^\ast, V_2^\ast, V_3^\ast) = \argmax_{\substack{V_k\in \mathbb{R}^{p_k\times r_k}\\\|V_k\|\leq 1}} \left\|\Z\times_1 V_1^\top \times_2 V_2^\top \times_3 V_3^\top\right\|, \quad \text{and}\quad T = \left\|\Z\times_1 (V_1^\ast)^\top \times_2 (V_2^\ast)^\top \times_3 (V_3^\ast)^\top\right\|.$$ 
Then we can find $V_1^{(a)}, V_2^{(b)}, V_3^{(c)}$ in the corresponding $\varepsilon$-nets such that
\begin{equation*}
\|V_1^\ast - V_1^{(a)}\| \leq \varepsilon, \quad \|V_2^\ast - V_2^{(b)}\|\leq \varepsilon,\quad \|V_3^\ast - V_3^{(c)}\|\leq \varepsilon.
\end{equation*}
Then
\begin{equation*}
\begin{split}
T = & \left\|\Z\times_1(V_1^\ast)^\top \times_2 (V_2^\ast)^\top \times_3 (V_3^\ast)^\top \right\|\\ 
\leq & \left\|\Z\times_1 (V_1^{(a)})^\top \times_2 (V_2^{(b)})^\top \times_3 (V_3^{(c)})^\top \right\| + \left\|\Z\times_1 (V_1^{(a)} - V_1^\ast)^\top \times_2 (V_2^\ast)^\top \times_3 (V_3^\ast)^\top \right\|\\
& + \left\|\Z\times_1 (V_1^{(a)})^\top \times_2 (V_2^{(b)} - V_2^\ast)^\top  \times_3 V_3^\top \right\| + \left\|\Z\times_1 (V_1^{(a)})^\top \times_2 (V_2^{(b)})^\top \times_3 (V_3^{(c)} - V_3^\ast)^\top \right\|\\
\leq & 2r_1r_2r_3 + 3t + \left(\|V_1^\ast - V_1^{(a)}\|+\|V_2^\ast - V_2^{(b)}\|+\|V_3^\ast - V_3^{(c)}\|\right)\cdot T,
\end{split}
\end{equation*}
which implies $T \leq (2r_1r_2r_3+3x)/(1-3\varepsilon)$ provided that $\varepsilon<1/3$ and \eqref{ineq:jointly_V1_V2_V3} holds. Let $\varepsilon=1/9$, 
$x = (1+t)\log(37)\cdot \left(p_1r_1 + p_2r_2 + p_3r_3\right)$ for some large constant $C>0$, by \eqref{ineq:jointly_V1_V2_V3} again we have
\begin{equation}
\begin{split}
& \mathbb{P}\left(T \geq Cr_1r_2r_3 + C(1+t)(p_1r_1+p_2r_2+p_3r_3)\right)\\ 
\geq & \exp(-Ct(p_1r_1+p_2r_2+p_3r_3))
\end{split}
\end{equation}
for some uniform constant $C>0$. 
thus we have finished the proof for \eqref{ineq:random_tensor_target2}. \quad $\square$

\subsection{Proof of Lemma \ref{lm:r-pinciple-compoenent}} 
\begin{equation*}
\begin{split}
\left\|\Proj_{\hat{U}_{2}}X\right\| \leq & \left\|P_{\hat{U}_{2}}(X+Z)\right\| + \|Z\| = \sigma_{r+1}(Y) + \|Z\| = \min_{\substack{\tilde{X}\in \mathbb{R}^{p_1\times p_2}\\\rank(\tilde{X})\leq r}} \|Y - \tilde{X}\| + \|Z\| \\
\leq & \|Y-X\| + \|Z\| =  2\|Z\|.
\end{split}
\end{equation*}
Since $\rank\left(\Proj_{U_{2}}X\right) \leq  \rank(X) \leq r$, it is clear that
\begin{equation*}
\left\|\Proj_{U_{2}}X\right\|_{\rm F} \leq 2\sqrt{r}\|Z\|;
\end{equation*}
meanwhile,
\begin{equation*}
\begin{split}
\left\|\Proj_{\hat{U}_{2}}X\right\|_{\rm F} \leq & \left\|P_{\hat{U}_{2}}(X+Z)\right\|_{\rm F} + \|Z\|_{\rm F} = \left(\sum_{i=r+1}^{p_1\wedge p_2}\sigma_{i}^2(Y)\right)^{1/2} + \|Z\|_{\rm F}\\
\leq & \min_{\substack{\tilde{X}\in \mathbb{R}^{p_1\times p_2}\\\rank(\tilde{X})\leq r}} \|Y - \tilde{X}\|_{\rm F} + \|Z\|_{\rm F} \leq \|Y-X\|_{\rm F} + \|Z\|_{\rm F} \leq 2\|Z\|_{\rm F},
\end{split}
\end{equation*}
which has proved this lemma. \quad $\square$

\subsection{Proof of Lemma \ref{lm:covering_set_spectral}} 
\begin{itemize}
	\item We first consider the $\varepsilon$-net for $\mathcal{X}_{p_1, p_2}$. Note that $\mathcal{X}_{p_1, p_2}$ is a convex set in $\mathbb{R}^{p_1\times p_2}$, we sequentially pick matrices from $\mathcal{X}_{p_1, p_2}$, say $X^{(1)}, X^{(2)}, \ldots$ satisfying the following criterion: for each time $t$, the picked matrix satisfies $\min_{t'\leq t} \|X^{(t)}-X^{(t-1)}\|_\bullet \geq \varepsilon$, i.e., the distances from $X^{(t)}$ to all the other selected matrices are at least $\varepsilon$. We stop the selection process until it is not possible to select the next matrix satisfying such criterion. 
	
	Suppose now $X^{(1)},\ldots, X^{(N)}$ are all we have selected. Since it is not possible to select another matrix from $\mathcal{X}_{p_1, p_2}$ which meets the criterion, all matrices in $\mathcal{X}_{p_1, p_2}$ must be within $\varepsilon$ of some selected matrix in $\{X^{(1)},\ldots, X^{(N)}\}$, thus
	$$\mathcal{X}_{p_1, p_2} \subseteq \cup_{i=1}^N B(X^{(i)}, \varepsilon).$$
	Here $B(X^{(i)}, \varepsilon) = \{X\in \mathbb{R}^{p_1\times p_2}: \|X - X^{(i)}\|_\bullet \leq \varepsilon\}$ is the closed ball with center $X^{(i)}$ and radius $\varepsilon$, Therefore, $\{X^{(1)},\ldots, X^{(N)}\}$ is a $\varepsilon$-net. 
	
	On the other hand, for any $1\leq i < j \leq N$, $\|X^{(i)}- X^{(j)}\|_\bullet \geq \varepsilon$, so
	$$\{X\in \mathbb{R}^{p_1\times p_2}: \|X\|_\bullet \leq 1+\varepsilon/2\} \supseteq \cup_{i=1}^N B(X^{(i)}, \varepsilon/2),  $$
	and $B(X^{(i)}, \varepsilon/2) \cap B(X^{(j)}, \varepsilon/2)$ contains at most one matrix for any $1\leq i < j \leq N$. Therefore,
	\begin{equation}
	\begin{split}
	& (1+\varepsilon/2)^{p_1p_2} {\rm vol}(\mathcal{X}_{p_1, p_2}) = {\rm vol}(\{X\in \mathbb{R}^{p_1\times p_2}: \|X\|_\bullet \leq 1+\varepsilon/2\}) \\
	\leq & \sum_{i=1}^N {\rm vol}(B^{(i)}, \varepsilon/2) = N (\varepsilon/2)^{p_1p_2} {\rm vol}(\mathcal{X}_{p_1, p_2}),
	\end{split}
	\end{equation}
	which implies $N\leq ((2+\varepsilon)/\varepsilon)^{p_1p_2}$.
	\item By the first part of this lemma, there exist $(\varepsilon/2)$-nets $\bar{\mathcal{X}}_{p_1, r}$ and $\bar{\mathcal{X}}_{r, p_2}$ for $\{\|X\in\mathbb{R}^{p_1\times r}: \|X\| \leq 1\}$ and $\{\|X\in\mathbb{R}^{r\times p_2}: \|X\| \leq 1\}$, such that 
	$$\left|\bar{\mathcal{X}}_{p_1, r}\right| \leq \left(\frac{4+\varepsilon}{\varepsilon}\right)^{p_1r},\quad \left|\bar{\mathcal{X}}_{r, p_2}\right| \leq \left(\frac{4+\varepsilon}{\varepsilon}\right)^{p_2r}.$$
	Next, we argue that
	$$\bar{\mathcal{F}}_{p_1, p_2, r} := \left\{X\cdot Y: X\in \bar{\mathcal{X}}_{p_1, r}, Y\in \bar{\mathcal{X}}_{r, p_2}\right\} $$
	is an $\varepsilon$-net for $\mathcal{F}_{p_1, p_2, r}$ in the spectral norm. Actually for any $X\in \mathcal{F}_{p_1, p_2, r}$, we can find $A, B$ such that $X = A\cdot B$, $A\in \mathbb{R}^{p_1\times r}, \|A\|\leq 1; B\in \mathbb{R}^{r\times p_2}, \|B\|\leq 1$. Then we can find $A^\ast \in \bar{\mathcal{X}}_{p_1, r}$ and $B^\ast \in \bar{\mathcal{X}}_{r, p_2}$ such that $\|A-A^\ast\| \leq \varepsilon/2, \|B-B^\ast\|\leq \varepsilon/2$, thus $A^\ast B^\ast \in \bar{\mathcal{F}}_{p_1, p_2, r}$ satisfies
	\begin{equation*}
	\begin{split}
	\|X- A^\ast B^\ast \| = & \left\|(AB-AB^\ast) + (AB^\ast - A^\ast B^\ast)\right\|\\
	\leq & \|A\|\cdot\|B-B^\ast\| + \|A-A^\ast \|\cdot \|B^\ast\| \leq 1\cdot \varepsilon/2 + 1\cdot \varepsilon/2  = \varepsilon.
	\end{split}
	\end{equation*}
	Note that $\left|\bar{\mathcal{F}}_{p_1, p_2, r}\right|\leq \left|\bar{\mathcal{X}}_{p_1, r}\right|\cdot \left|\bar{\mathcal{X}}_{r, p_2}\right| \leq \left((4+\varepsilon)/\varepsilon\right)^{r(p_1+p_2)}$, this has finished the proof of this lemma. \quad $\square$
\end{itemize}

\subsection{Proof of Lemma \ref{lm:gaussian-vector-projection}} Suppose $A = U\Sigma V^\top$ is the singular value decomposition of $A$. Since $U, V$ are orthogonal and $u\overset{iid}{\sim}N(0, 1)$, $\|Au\|_2^2$ has the same distribution as $\sum_{i=1}^{p\wedge n} \sigma_i(A)^2 u_i^2$. By the exponential probability for general chi-square distribution (Lemma 1 in \cite{laurent2000adaptive}), we have
\begin{equation*}
\begin{split}
& P\left(\sum_{i=1}^{p\wedge n}\sigma_i^2(A)u_i^2 - \sum_{i=1}^{p\wedge n}\sigma_i^2(A) \leq -2\sqrt{t\sum_{i=1}^{p\wedge n}\sigma_i^4(A)}\right) \leq \exp(-x);\\ 
& P\left(\sum_{i=1}^{p\wedge n}\sigma_i^2(A)u_i^2 - \sum_{i=1}^{p\wedge n}\sigma_i^2(A) \geq 2\sqrt{t\sum_{i=1}^{p\wedge n}\sigma_i^4(A)} + 2t\max \sigma_i^2(A)\right) \leq \exp(-x),
\end{split}
\end{equation*}
which has finished the proof for Lemma \ref{lm:random_tensor_max_projection} since $\|A^\top A\|_{\rm F}^2 = \sum_{i=1}^{p\wedge n} \sigma_i^4(A)$, and $\|A\|=\max_i \sigma_i(A)$.\quad $\square$

\subsection{Proof of Lemma~\ref{lemma:tv}}
Clearly, it suffices to prove the claim for $i=0$, i.e., under $H_0$. Let $G=(V,E)\sim H_0$ with $V=\{1,2,\ldots,N\}$ and $\A$ denote its adjacency tensor, meaning that there is a clique of size $\kappa_N$ planted in the subset $\{1,2,\ldots,\floor{N/2}\}$.
Recall that $N=3p$ for an even integer $p$.
The vertices set of the planted clique is denoted by $C\subset \big\{1,\ldots,\frac{3p}{2}\big\}$ with $|C|=\kappa_N=20k$ where $k=\floor{p^{(1-\tau)/2}}$. Recall $V_1, V_2, V_3$ and define
$$
C_j:= C\cap V_j,\quad j=1,2,3
$$
which represents the subsets of clique vertices in $V_1, V_2, V_3$. 
If $\Y=\mathcal{T}(\A, {\bf \Xi}^-, {\bf \Xi}^+)\in\mathbb{R}^{p\times p\times p}$, it is clear that, under $H_0$, $\X=\mathbb{E}(\Y|C)$ is a sparse tensor with supports $S_1(\X)=C_1\subset [p/2], S_2(\X)=C_2-\frac{p}{2}\subset [p/2]$ and $S_3(\X)=C_3-p\subset [p/2]$.
We show that the sizes of $S_1(\X), S_2(\X), S_3(\X)$ are lower bounded by $k$ with high probability. 
\begin{Lemma}\label{s123lemma}
	There exists an event $\mathcal{E}$ on which $\min\{|S_1(\X)|,|S_2(\X)|,|S_3(\X)|\}\geq k$ and 
	$$
	\mathbb{P}({\mathcal{E}})\geq 1-6k(0.86)^{2.5k}.
	$$
\end{Lemma}
For any fixed realization $G\sim H_0$ with set of clique vertices $C=C_1\cup C_2\cup C_3$ (with corresponding supports $S_k:=S_k(\X), k=1,2,3$), we generate a Gaussian random tensor $\tilde{\Y}\in\mathbb{R}^{p\times p\times p}$ with independent entries such that
$$
\tilde{Y}(a,b,c)\sim\mathcal{N}(\mu,1)\quad \textrm{if } (a,b,c)\in S_1\times S_2\times S_3;\quad \tilde{Y}(a,b,c)\sim\mathcal{N}(0,1)\quad \textrm{otherwise},
$$
where $S_1=C_1, S_2=C_2-\frac{p}{2}$ and $S_3=C_3-p$. By Lemma~\ref{f12lemma}, we have 
$$
{\rm d_{TV}}\Big(\mathcal{L}\big(Y(a,b,c)\big| C\big),\mathcal{L}\big(\tilde{Y}(a,b,c)\big| C\big)\Big)\leq e^{(1-M^2)/2},\quad \forall\ a,b,c\in \{1,2,\ldots,p\}.
$$
As a result, since $M= \sqrt{8\log N}$ and $p=N/3$,
\begin{align*}
{\rm d_{TV}}\big(\mathcal{L}(\Y), \mathcal{L}(\tilde{\Y})\big)=&\mathbb{E}_C{\rm d_{TV}}\big(\mathcal{L}(\Y|C), \mathcal{L}(\tilde{\Y}|C)\big)\\
=&\mathbb{E}_C\sum_{a,b,c=1}^{p}{\rm d_{TV}}\Big(\mathcal{L}\big(Y(a,b,c)\big| C\big),\mathcal{L}\big(\tilde{Y}(a,b,c)\big| C\big)\Big)\leq p^3e^{(1-M^2)/2}\leq \frac{\sqrt{e}}{27N}.
\end{align*}
Now we show that $\mathcal{L}(\tilde{\Y}|\mathcal{E})$ is a mixture over $\big\{\mathbb{P}_{\X}, \X\in\mathcal{M}_0(\bp,k,\br,\lambda)\big\}$ with $\lambda=\frac{p^{3(1-\tau)/4}}{2\sqrt{8\log 3p}}$. Indeed, for any fixed $C$, let $\tilde\X=\mathbb{E}(\tilde\Y|C)$. Then,
$$
\tilde{X}(a,b,c)=\mathbb{E}\Big(\tilde{Y}(a,b,c)|C\Big)= \mu,\quad \forall (a,b,c)\in S_1\times S_2\times S_3.
$$
Recall that on $\mathcal{E}$, $\min\{|S_1|, |S_2|, |S_3|\}\geq k$. Therefore, $\tilde\X$ is of rank $1$ and on $\mathcal{E}$,
\begin{eqnarray*}
	\min\big\{\sigma_{\min}\big(\mathcal{M}_1(\tilde{\X})\big), \sigma_{\min}\big(\mathcal{M}_2(\tilde{\X})\big), \sigma_{\min}\big(\mathcal{M}_3(\tilde{\X})\big)\big\}\\
	\geq \mu\sqrt{|S_1||S_2||S_3|}
	\geq \mu k ^{3/2}\geq \mu \floor{p^{3(1-\tau)/4}}=\frac{\floor{p^{3(1-\tau)/4}}}{2\sqrt{8\log 3p}}
\end{eqnarray*}
since $\mu=\frac{1}{2M}$. The above fact indicates that $\tilde\X\in\mathcal{M}_0(\bp,k,\br,\lambda)$. In other words, under $H_0$, for any $C$ conditioned on $\mathcal{E}$, there exists $\X(C)\in\mathcal{M}_0(\bp,k,\br,\lambda)$ such that $\mathcal{L}(\tilde{\Y}|C)=\mathbb{P}_{\X(C)}$. Define the probability distribution $\pi_0=\mathcal{L}\big(\X(C)|\mathcal{E}\big)$ supported on $\mathcal{M}_0(\bp,k,\br,\lambda)$. Then $\mathcal{L}(\tilde{\Y}|\mathcal{E})=\mathbb{P}_{\pi_0}$ and 
\begin{align*}
{\rm d_{TV}}\big(\mathcal{L}(\Y),\mathbb{P}_{\pi_0}\big)\leq& {\rm d_{TV}}\big(\mathcal{L}(\Y),\mathcal{L}(\tilde{\Y})\big)+\textrm{TV}\big(\mathcal{L}(\tilde{\Y}),\mathbb{P}_{\pi_0}\big)\\
\leq& \frac{\sqrt{e}}{27N}+\mathbb{P}(\mathcal{E}^c)\leq \frac{\sqrt{e}}{27N}+6k\big(0.86\big)^{2.5k}.
\end{align*}

\subsection{Proof of Lemma~\ref{s123lemma}}
Recall that $\kappa\leq \sqrt{N/2}$ and N=3p. Let $N_1:=N/2=\frac{3p}{2}$.
Since $C$ is uniformly chosen from $\{1,2,\ldots,\frac{3p}{2}\}$ and $C_1\subset \big\{1,2,\ldots,\frac{p}{2}\big\}$, we have
\begin{align*}
\mathbb{P}\Big(|C_1|\leq &\frac{\kappa}{8}\Big)\leq \frac{\sum_{s=0}^{\kappa/8}{p/2\choose s}{p\choose \kappa-s}}{{N_1\choose \kappa}}\leq \frac{\kappa+1}{8}\frac{{p/2\choose \kappa/8}{p\choose 7\kappa/8}}{{N_1\choose \kappa}}\\
=&\frac{\kappa+1}{8}{\kappa\choose \kappa/8}\frac{(p)(p-1)\ldots(p-7\kappa/8+1)(p/2)(p/2-1)\ldots(p/2-\kappa/8+1)}{(3p/2)(3p/2-1)\ldots (3p/2-\kappa+1)}\\
\leq &\frac{\kappa+1}{8}(8e)^{\kappa/8}\big(\frac{2}{3}\big)^{\kappa}=\frac{\kappa}{8}\Big(8e\cdot2^8/3^8\Big)^{\kappa/8}\leq \frac{\kappa+1}{8}(0.86)^{\kappa/8}
\end{align*}
where we used the fact ${p/2\choose s}{p\choose \kappa-s}$ increases for $0\leq s\leq \kappa/8$ and inequality ${n\choose k}\leq (ne/k)^k$. Therefore, with probability at least $1-\frac{\kappa+1}{4}(0.86)^{\kappa/8}$,
$$
\frac{\kappa}{8}\leq |C_1|, |C_2|, |C_3|\leq \frac{7\kappa}{8}.
$$
Recall that $\kappa=20k$ and we conclude the proof.

\end{document}